\documentclass[12pt,reqno]{amsart}  
\usepackage{amsmath,amsthm,amssymb,mathtools,stmaryrd,enumerate,hyperref,xcolor,ytableau,tikz}
\mathtoolsset{showonlyrefs}

\newcommand{\bfa}{\boldsymbol{a}}
\newcommand{\bfb}{\boldsymbol{b}}
\newcommand{\bfc}{\boldsymbol{c}}
\newcommand{\bfd}{\boldsymbol{d}}
\newcommand{\bfk}{\boldsymbol{k}}
\newcommand{\bft}{\boldsymbol{t}}
\newcommand{\fh}{\mathfrak{h}}
\newcommand{\EE}{\mathcal E}
\newcommand{\DD}{\mathfrak{D}}
\newcommand{\calT}{\mathcal{T}}
\newcommand{\ZZ}{\mathbb{Z}}
\newcommand{\QQ}{\mathbb{Q}}
\newcommand{\Zgeq}{\mathbb{Z}_{+}}
\newcommand{\NN}{\mathbb{N}}

\hypersetup{pdfborder=0 0 0, pdftitle={Multiple q-zeta values as polynomials in partition multiplicities}, pdfauthor={Jan-Willem M. van Ittersum}, pdfkeywords={multiple q-zeta values, partitions, quasi-shuffle products, partial domino tilings, q-series identities}}
\usetikzlibrary{arrows}
\usetikzlibrary{calc}
\usetikzlibrary{decorations}
\usetikzlibrary{decorations.pathreplacing}
\tikzset{>=stealth'}
\usepackage{booktabs}
\usepackage{array}
\usepackage{longtable}
\makeatletter\def\p@enumi#1{\textup{(#1)}}\makeatother

\title[Multiple $q$-zeta values as polynomials in partition multiplicities]{Multiple $q$-zeta values \\as polynomials in partition multiplicities}

\author{Jan-Willem M.~van Ittersum}%
\email{\href{mailto:j.w.m.vanittersum@uva.nl}{j.w.m.vanittersum@uva.nl}} \address{Korteweg--de Vries Institute for Mathematics, University of Amsterdam, Postbus 94248, 1090 GE  Amsterdam, The Netherlands}



\newtheorem{theorem}{Theorem}[section]
\newtheorem{proposition}[theorem]{Proposition}

\newtheorem{lemma}[theorem]{Lemma}
\newtheorem{corollary}[theorem]{Corollary}

\numberwithin{equation}{section}
\theoremstyle{definition}
\newtheorem{definition}[theorem]{Definition}
\newtheorem{remark}[theorem]{Remark}
\newtheorem{question}[theorem]{Question}

\newcommand{\N}{\mathbb{N}}
\newcommand{\Q}{\mathbb{Q}}
\newcommand{\df}{\coloneqq}
\newcommand{\im}{\operatorname{im}}

\definecolor{dutchorange}{RGB}{220,50,0}
\definecolor{darkblue}{RGB}{0,71,171}

\newcommand{\dep}{\operatorname{dep}}
\newcommand{\wt}{\operatorname{wt}}

\renewcommand{\=}{\: =\: }
\newcommand{\defis}{\: :=\: }
\newcommand{\+}{\,+\,}
\newcommand{\meno}{\,-\,}
\newcommand{\msp}{\mspace{1mu}}

\newcommand{\A}{\mathcal{A}}
\newcommand{\QB}{\Q\langle \B \rangle}
\newcommand{\B}{\mathcal{B}}

\newcommand{\gs}{\mathcal{G}}
\newcommand{\Pp}{\mathsf{P}}
\newcommand{\calU}{\mathcal{U}}

\begin{document}
\begin{abstract}
Writing multiple~$q$-zeta values as sums over partitions, a conjecture of Bachmann, proved by Hirose, Maesaka and Watanabe by means of a recursion, says that polynomials in the parts and multiplicities can be replaced by polynomials in the multiplicities alone. We make this replacement explicit, derive from it bounds on weight and depth that make Bachmann's conjecture effective, and provide a closed formula for the Drop1 operator of Hirose, Maesaka, Seki and Watanabe.
\end{abstract}
\maketitle

\section{Introduction}
The Eisenstein series $P,Q$ and $R\mspace{1mu}$, given by
\begin{align}
P \df 1 - 24 \sum_{m,n\geq 1} m\, q^{mn}, \qquad 
Q \df 1 + 240 \sum_{m,n\geq 1} m^3\, q^{mn}, \qquad 
R \df 1 - 504 \sum_{m,n\geq 1} m^5\, q^{mn},
\end{align}
satisfy Ramanujan's differential equations
\begin{align}\label{eq:ramanujandiffeq}
q\frac{d}{dq}P \= \frac{P^2-Q}{12}\,, \qquad q\frac{d}{dq} Q \= \frac{PQ-R}{3}\,, \qquad q\frac{d}{dq} R \= \frac{PR-Q^2}{2}\,.
\end{align}
Hence, the algebra $\Q[P,Q,R]$ of quasimodular forms is closed under the differential operator~$q\frac{d}{dq}\mspace{1mu}$. More generally, for any polynomial $p\in\Q[x,y]$ all of whose monomials have odd total degree, the series
\begin{equation}\label{eq:qmf}  \sum_{m,n\geq 1} p(m,n)\,q^{mn} \end{equation}
is a quasimodular form (a derivative of an Eisenstein series). In other words, for such~$p$ the series in~\eqref{eq:qmf} can always be expressed as a \emph{polynomial} in series of the form $\sum_{m,n\geq 1} \tilde{p}(m)\,q^{mn} $ for $\tilde{p}\in \Q[x]\mspace{1mu}$, that is, as a linear combination of series of the form \[
\sum_{\substack{m_1,\ldots,m_r>0 \\0<n_1<\ldots<n_r}} \prod_{i=1}^r p_i(m_i)\,q^{m_in_i} 
\]
with $r\geq 0$ and $p_1,\ldots,p_r\in \Q[x]\mspace{1mu}$. 

The latter observation holds for every $p\in\Q[x,y]\mspace{1mu}$, as the simplest case of Theorem~\ref{thm:phi} below; in this paper we make it explicit. 
For example, for $p(x,y)=xy$ the $q$-series in~\eqref{eq:qmf} can be written as
\begin{align}\label{eq:inparts}
    \sum_{m,n\geq 1} {\color{dutchorange}m}\mspace{1mu} {\color{darkblue}n} \, q^{mn} &\= \frac{1}{2}\sum_{m,n\geq 1} ({\color{darkblue}
n}^2+{\color{darkblue}n}) \, q^{mn} \meno \sum_{\substack{m_1,m_2 \geq 1\\0<n_1<n_2}} {\color{darkblue}n_1}\, q^{m_1 n_1 + m_2 n_2}.
\end{align}
\begin{figure}[t]
\centering
\begin{tikzpicture}[inner sep=0in,outer sep=0in,scale=0.34]
\ytableausetup{smalltableaux}
\node at (-1.5,0) (n) {\ydiagram{6,6,6,2,2}};
\node at ($(n.east)+(11,0)$) (m) {\ydiagram{5,5,3,3,3,3}};
\node (n2M) at ($(n.north east)+(0,-3)$) {\ };
\node (m2M) at ($(m.north east)+(0,-2)$) {\ };
\node (n1M) at ($(n2M)+(-4,-0.1)$) {\ };
\node (m1M) at ($(m2M)+(-2,-0.1)$) {\ };
\node (n0M) at ($(n1M)+(0,-2)$) {\ };
\node (m0M) at ($(m1M)+(0,-4)$) {\ };
\draw[thick,darkblue,decorate,decoration={brace,amplitude=4pt}] ($(n.north west)+(0,0.2)$)--($(n.north east)+(0,0.2)$);
\node at ($0.5*(n.north west)+0.5*(n.north east)+(0,1.2)$) {$\textcolor{darkblue}{n_2}$};
\draw[thick,darkblue,decorate,decoration={brace,amplitude=4pt}] ($(n.south west)+(2,-0.2)$)--($(n.south west)+(0,-0.2)$);
\node at ($0.5*(n0M)+0.5*(n.south west)+(0,-0.2) + (0,-1.2)$) {$\textcolor{darkblue}{n_1}$};
\draw[thick,dutchorange,decorate,decoration={brace,amplitude=4pt}] ($(m.north west)+(0,0.2)$)--($(m.north east)+(0,0.2)$);
\draw[thick,dutchorange,decorate,decoration={brace,amplitude=4pt}] ($(m.south west)+(3,-0.2)$)--($(m.south west)+(0,-0.2)$);
\node at ($0.5*(m.north west)+0.5*(m.north east)+(0,1.2)$) {$\textcolor{dutchorange}{m_1+m_2}$};
\node at ($0.5*(m0M)+0.5*(m.south west)+(0,-0.2) + (0,-1.2)$) {$\textcolor{dutchorange}{m_2}$};

\draw[thick,dutchorange,decorate,decoration={brace,amplitude=4pt}] ($(n.north east)+(0.2,0)$)--($(n2M)+(0.2,0)$);
\node at ($0.5*(n.north east)+0.5*(n2M)+(1.6,0)$) {$\textcolor{dutchorange}{m_2}$};
\draw[thick,dutchorange,decorate,decoration={brace,amplitude=4pt}] ($(n1M)+(0.2,0)$)--($(n0M)+(0.2,0.1)$);
\node at ($0.5*(n1M)+0.5*(n0M)+(1.6,0)$) {$\textcolor{dutchorange}{m_1}$};
\draw[thick,darkblue,decorate,decoration={brace,amplitude=4pt}] ($(m.north east)+(0.2,0)$)--($(m2M)+(0.2,0)$);
\node at ($0.5*(m.north east)+0.5*(m2M)+(1.6,0)$) {$\textcolor{darkblue}{n_1}$};
\draw[thick,darkblue,decorate,decoration={brace,amplitude=4pt}] ($(m1M)+(0.2,0)$)--($(m0M)+(0.2,0.1)$);
\node at ($0.5*(m1M)+0.5*(m0M)+(3,0)$) {$\textcolor{darkblue}{n_2-n_1}$};
\draw[thick,|->] ($(n.east)+(3.8,0)$)--($(m.west)+(-2.5,0)$);
\end{tikzpicture}
\caption{Conjugating a partition with two distinct parts $n_1<n_2$ of multiplicities~$m_1$ and~$m_2\mspace{1mu}$: the parts become $m_2$ and $m_1+m_2\mspace{1mu}$, of multiplicities $n_2-n_1$ and~$n_1\mspace{1mu}$.}\label{fig:conjugate}
\end{figure}
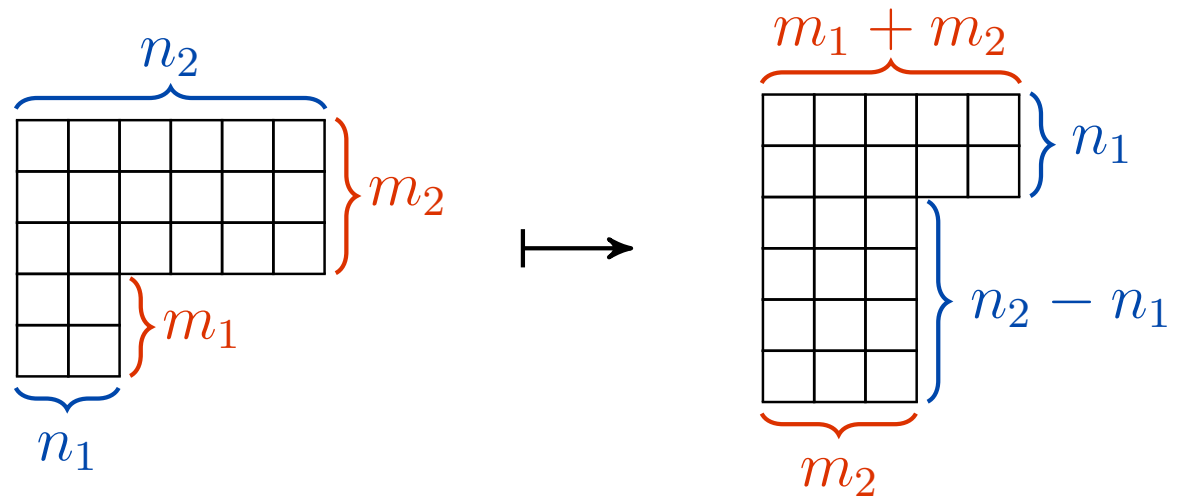%
Note that in the last sum, the coefficient of $q^N$ is given by the sum over all partitions of $N$ with parts~$n_1$ and~$n_2$ of multiplicities~$m_1$ and~$m_2$ respectively, whereas the other sums can be interpreted as sums over partitions consisting of only one part~$n$ of multiplicity~$m\mspace{1mu}$. Transposing (the Young diagram of) these partitions, as in Figure~\ref{fig:conjugate}, yields
\begin{align}\label{eq:derg1}
    \sum_{m,n\geq 1} {\color{dutchorange}m}\mspace{1mu} {\color{darkblue}n} \, q^{mn} &\= \frac{1}{2}\sum_{m,n\geq 1} ({\color{dutchorange}m}^2+{\color{dutchorange}m}) \, q^{mn} \meno  \sum_{\substack{m_1,m_2 \geq 1\\0<n_1<n_2}}  {\color{dutchorange}m_2}\, q^{m_1 n_1 + m_2 n_2}.
\end{align}
So \eqref{eq:inparts} and \eqref{eq:derg1} write this series with polynomials in just the {\color{darkblue}parts} or just the {\color{dutchorange}multiplicities}.

Fix $\epsilon\in\{0,1\}$ and let \[\Psi\=\Psi^{\epsilon}\,:\, \Q\,\oplus\,\bigoplus_{r\geq 1} x_r^{\epsilon}\,\Q[x_1,\ldots,x_r,y_1,\ldots,y_r] \,\to\, \Q\llbracket q\rrbracket \]
be the $\Q$-linear map uniquely determined by $1\mapsto 1$ and
\begin{equation} \label{eq:phi} p(x_1,\ldots,x_r,y_1,\ldots,y_r)\,\mapsto \sum_{\substack{m_1,\ldots,m_r\geq 1-\epsilon\\0<n_1<\ldots<n_r}} p(m_1,\ldots,m_r,n_1,\ldots,n_r)\, q^{m_1 n_1+\ldots+m_r n_r}.\end{equation}
For $\epsilon=1\mspace{1mu}$, the divisibility of~$p$ by~$x_r$ forces $m_r\geq1\mspace{1mu}$, so that the coefficient of every power of~$q$ is a finite sum. Write $\Psi^{\epsilon}_X\df\Psi_X$ for the restriction of $\Psi^{\epsilon}$ to $\Q\oplus\bigoplus_{r\geq 1} x_r^{\epsilon}\,\Q[x_1,\ldots,x_r]$ and $\Psi^{0}_Y$ for the restriction of $\Psi^{0}$ to $\bigoplus_{r\geq 0} \Q[y_1,\ldots,y_r]\mspace{1mu}$. Note that the coefficient of $q^N$ in~\eqref{eq:phi} is the sum of the polynomial~$p$ over all partitions of size~$N$ with (at most) $r$ different parts $n_1,\dots,n_r$ with multiplicities $m_1,\dots,m_r\msp$. Hence, by using the conjugation of partitions, we see $\im \Psi^{0}_X = \im \Psi^{0}_Y\mspace{1mu}$. Further, it is known that $\im\Psi^{1}\subseteq\im\Psi^{0}$.

Recently, Hirose, Maesaka and Watanabe \cite{HMW26} proved the following two conjectures of Bachmann, stated in \cite[Conjecture~4.3]{Ba} and \cite[Conjecture~4.1]{Br2}; both are recorded in \cite[Conjecture~1, Figure~1]{Br}, and partial results were obtained in \cite{Ba,V,Bu1,Br3}. 

\begin{theorem}\label{thm:phi} We have
\begin{align}
   \im\Psi^{0} &\= \im\Psi^{0}_X\,, \qquad   \im\Psi^{1}_X \= \Psi^{0}_X\bigl(\textstyle\bigoplus_{r\geq0}x_1\cdots x_r\,\Q[x_1,\ldots,x_r]\bigr)\,.\label{eq:Zq}
\end{align}
\end{theorem}

The content of the first equation is that, under $\Psi\mspace{1mu}$, polynomials in the parts \emph{and} the multiplicities give nothing beyond polynomials in the multiplicities alone (or, by conjugation, in the parts alone). The second equality can be interpreted as the statement that, for polynomials in the multiplicities alone, extending the sum in~\eqref{eq:phi} so as to allow the multiplicities to vanish gives nothing beyond restricting to polynomials divisible by all the multiplicities. Both equalities are a model-independent form of a more explicit statement of \cite{HMW26}, which we recall below.

The $q$-series in the image of $\Psi^{0}$ have been studied by various authors due to their connection to multiple zeta values (e.g., \cite{Ba,BK,BvI,Br,Bu1,EMS,T,Zh,Zu}). \emph{Multiple zeta values} are real numbers defined by 
\begin{align}\label{eq:defmzv}
    \zeta(k_1,\dots,k_r) \df \sum_{0<n_1 < \cdots < n_r} \frac{1}{n_1^{k_1} \cdots n_r^{k_r}}\qquad (k_1,\dots,k_{r-1}\geq 1,\ k_r\geq 2)\mspace{1mu}.
\end{align}
We write $\mathcal{Z}$ for the $\Q$-vector space spanned by all multiple zeta values. In~\cite{BvI} it has been shown that every $q$-series in the image of $\Psi^{0}$ is a $q$-analogue of multiple zeta values and we therefore write $ \mathcal{Z}_q \df \im \Psi^{0}. $ More precisely, \cite[Theorem~1.2]{BvI} gives for $f\in\mathcal{Z}_q$ an exponent $k\geq0$ for which the limit $\lim_{q\to 1} (1-q)^k f$ exists and lies in $\mathcal{Z}\mspace{1mu}$, and every element of $\mathcal{Z}$ arises in this way.

Hirose, Maesaka and Watanabe introduced the following model of $\mathcal Z_q\mspace{1mu}$: for $\bfc=(c_1,\dots,c_d)$ with $c_j\geq0$ and $c_d\neq0$ put
\begin{equation}\label{eq:HMW-model}
  \zeta^\dagger_q(\bfc)\df\sum_{\substack{0<n_1\leq\dots\leq n_d\\
  n_j<n_{j+1}\text{ if }c_j\neq0}}\ \prod_{j\,:\,c_j\neq0}\frac{q^{n_j}}{(1-q^{n_j})^{c_j}}\ \in \, \mathcal Z_q\mspace{1mu}.
\end{equation}
These series span~$\im\Psi^{0}$, and those with all $c_j\geq1$ already span~$\im\Psi^{0}_X$ (see Section~\ref{sec:models}). Writing $\bfc=(0^{\ell_1-1},k_1,\dots,0^{\ell_r-1},k_r)\mspace{1mu}$, their theorem states that
\begin{equation}\label{eq:HMWintro}
  \zeta^\dagger_q(0^{\ell_1-1},k_1,\dots,0^{\ell_r-1},k_r)
  \= Z^\dagger_q\bigl(\EE_{0}(\ell_1,k_1,\dots,\ell_r,k_r)\bigr)\mspace{1mu},
\end{equation}
where $\EE_{0}(\ell_1,k_1,\dots,\ell_r,k_r)$ (the $\mathfrak E_q(\ell_1,k_1,\dots,\ell_r,k_r)$ of \cite{HMW26}), defined by the recursion recalled in Definition~\ref{def:Eq}, is an explicit element  of the free algebra $\fh\df\ZZ\langle x,y\rangle$ in two non-commuting letters and $Z^\dagger_q$ is the linear map sending the word $yx^{c_1-1}\cdots yx^{c_s-1}$ to $\zeta^\dagger_q(c_1,\dots,c_s)\mspace{1mu}$. So \eqref{eq:HMWintro} is an explicit form of the equality $\im\Psi^{0}=\im\Psi^{0}_X\mspace{1mu}$. A second family~$\EE_{1}(\bfc)\mspace{1mu}$, defined by the same recursion with one exponent changed, does the same for the other equality of Theorem~\ref{thm:phi} (see Theorem~\ref{thm:HMW}). What we give in this paper is a solution of that recursion for $\EE_{0}$ and $\EE_{1}$ in closed form. 

To state it, fix $\epsilon\in\{0,1\}$ and let $t_1,t_2,\ldots$ be commuting variables, central in $\fh\llbracket t_1,t_2,\ldots\rrbracket\mspace{1mu}$. For a finite set $S=\{s_1<\dots<s_{2m}\}\subset\Zgeq$ of even cardinality let \[
  \gs^{\epsilon}_S\df\gs_S\df\sum_{c_1,\ldots,c_{2m}\geq1}\EE_{\epsilon}(c_1,\ldots,c_{2m})\,t_{s_1}^{c_1}\cdots t_{s_{2m}}^{c_{2m}}
\]
be the generating series of the words, so that $\gs_\emptyset=1\mspace{1mu}$, and write $\gs^{\epsilon}_r\df\gs_r\df\gs_{[2r]}$ with $[2r]\df\{1,2,\ldots,2r\}\mspace{1mu}$. We call $D_i(S)\df\{s_i,s_{i+1}\}$ a \emph{domino} of $S\mspace{1mu}$, even--odd or odd--even as $i$ is even or odd, and we write $D_i\df\{i,i+1\}\mspace{1mu}$. Let $\calT(S)$ be the set of \emph{partial domino tilings} of~$S\mspace{1mu}$, that is, disjoint unions of dominoes, no two of them odd--even and occupying four consecutive positions of~$S\mspace{1mu}$, and for $A\in\calT(S)$ let $\kappa(A)$ be the number of its dominoes and $\mathrm{eo}(A)$ the number of even--odd dominoes among them (see Section~\ref{subsec:notation}).

Writing $\hat t_s\df t_s/(1-x\msp t_s)\mspace{1mu}$, we attach to every domino $D=\{u,v\}$ the variable $\eta_D\df x(x-1)\,\hat t_u\hat t_v$ and to a disjoint union~$A$ of dominoes the product $\eta_A\df\prod_D\eta_D$ over them, so that $\eta_\emptyset=1\mspace{1mu}$. We define \[
  \Pp_S\df\prod_{\substack{0<i<2m\\ i\text{ odd}}}\bigl(1-\eta_{D_i(S)}\bigr)\mspace{1mu},
  \qquad
  \Pp^{*}_S\df\prod_{\substack{0<i<2m\\ i\text{ even}}}\bigl(1-\eta_{D_i(S)}\bigr)\mspace{1mu}.
\]
In terms of the series~$\gs_S$ the recursion defining~$\EE_{0}$ and $\EE_{1}$ (Definition~\ref{def:Eq}) becomes the following identity, which we prove in Proposition~\ref{prop:rec}: for every nonempty such set~$S\mspace{1mu}$,
\begin{equation}\label{eq:recw}
  \gs_S\=\sum_{A\in\calT(S)}\gs_{S\setminus A}\;y\;\frac{w_{A,S}}{x-1}\,,
\end{equation}
where the \emph{weights} of the recursion are
\begin{equation}\label{eq:weights}
  w^{\epsilon}_{A,S}\df w_{A,S}\df(-1)^{\mathrm{eo}(A)}\,x^{(\epsilon-1)\kappa(A)}\,\frac{\Pp_{S\setminus A}}{\Pp^{*}_S}\,\eta_{A}-\delta_{A,\emptyset}\qquad(A\in\calT(S))\mspace{1mu}.
\end{equation}
In particular, $w_{\emptyset,S}=\Pp_S/\Pp^{*}_S-1$ and $w_{\emptyset,\emptyset}=0\mspace{1mu}$.

Let $\calU(r)$ be the set of $U\subseteq[2r]$ for which $[2r]\setminus U$ is a disjoint union of dominoes of~$[2r]$ (this set has $F_{2r+1}$ elements, with $F_n$ the $n$th Fibonacci number). Let $\mathsf C^{\epsilon}_r\df\mathsf C_r$ be the matrix with rows and columns indexed by~$\calU(r)$ and entries \[
  \mathsf C^{\epsilon}_{S\setminus A,S}\df\mathsf C_{S\setminus A,S}\df\frac{w_{A,S}}{1-x}\qquad(A\in\calT(S))\mspace{1mu},
\]
and set $\mathsf C_{T,S}\df0$ for all other~$T\mspace{1mu}$. By Lemma~\ref{lem:Cring} these entries are power series in the~$t_s$ with coefficients in~$\ZZ[x]$ and without constant term, so that $1+y\,\mathsf C_r$ is invertible over $\fh\llbracket t_1,\dots,t_{2r}\rrbracket\mspace{1mu}$. As the column $S=\emptyset$ of~$\mathsf C_r$ vanishes, the identities~\eqref{eq:recw} for $S\in\calU(r)\mspace{1mu}$, together with $\gs_\emptyset=1\mspace{1mu}$, say that the row vector $\gs^{\epsilon}\df\gs\df(\gs_U)_{U\in\calU(r)}$ satisfies \[
  \gs\,(1+y\,\mathsf C_r)\=e_\emptyset\,,
\]
with $e_\emptyset$ the unit row vector at~$\emptyset\mspace{1mu}$. Solving for the entry at~$[2r]$ gives the following theorem.

\begin{theorem}\label{thm:intro}
For all $r\geq1$ and $\epsilon\in\{0,1\}\mspace{1mu}$, \[
  \gs^{\epsilon}_r\=\bigl[(1+y\,\mathsf C^{\epsilon}_r)^{-1}\bigr]_{\emptyset,[2r]}\,.
\]
\end{theorem}

In depth one $\calU(1)=\{\emptyset,[2]\}$ and
\begin{equation}\label{eq:depthone}
  \mathsf C_1\=\begin{pmatrix}0&-x^{\epsilon}\hat t_1\hat t_2\\[2pt]0&x\,\hat t_1\hat t_2\end{pmatrix}\,,
  \qquad\text{so that}\qquad
  \gs_1   \=y\,x^{\epsilon}\hat t_1\hat t_2\,\bigl(1+y\,x\,\hat t_1\hat t_2\bigr)^{-1}\,.
\end{equation}
Expanding the inverse and taking the coefficient of $t_1^2t_2^2$ gives $\EE_{0}(2,2)=yx^{2}-y^{2}x\mspace{1mu}$, and applying~$Z^\dagger_q$ to it yields \[
  \zeta^\dagger_q(0,2)\=\zeta^\dagger_q(3)\meno\zeta^\dagger_q(1,2)\mspace{1mu},
\]
 which is \eqref{eq:derg1}. Namely, the left-hand side is $\sum_{m,n\geq1}mn\,q^{mn}$ and the two terms on the right are its two partition expressions in the multiplicities alone.

Such a rewriting exists for every polynomial in the parts and the multiplicities: Theorem~\ref{thm:phi} amounts to the existence of $\Q$-linear maps \[
  \mathcal L^{0}\colon\bigoplus_{r\geq0}\Q[x_1,\ldots,x_r,y_1,\ldots,y_r]\to
  \bigoplus_{r\geq0}\Q[x_1,\ldots,x_r]\,,\qquad\text{with}\qquad \Psi^{0}\=\Psi^{0}_X\circ\mathcal L^{0}\,,
\]
and \[
  \mathcal L^{1}\colon\bigoplus_{r\geq0}x_r\Q[x_1,\ldots,x_r]\to\bigoplus_{r\geq0}x_1\cdots x_r\Q[x_1,\ldots,x_r]\,,\qquad\text{with}\qquad \Psi^{1}_X\=\Psi^{0}_X\circ\mathcal L^{1}\,.
\]
As $\Psi^{0}_X$ has a large kernel, such maps are not unique; the words~$\EE_{\epsilon}$ select one (cf.\ Section~\ref{sec:Lmaps}), and Theorem~\ref{thm:intro} makes it explicit. For example, for $p=x_1y_1$ the identity $\Psi^{0}(p)=\Psi^{0}_X(\mathcal L^{0}p)$ is \eqref{eq:derg1}.

By Proposition~\ref{prop:star}, the map~$\Psi^{\epsilon}$ is a homomorphism of algebras for a stuffle product~$\ast^{\epsilon}$ on polynomials which interleaves the parts of two partitions in all possible ways and adds the multiplicities of coinciding parts. The map~$\mathcal L^{\epsilon}$ is not a homomorphism of algebras: one checks, e.g., that $\mathcal L^{0}(y_1\ast^{0}y_1)\neq \mathcal L^{0}(y_1)\ast^{0}\mathcal L^{0}(y_1)$ and $\mathcal L^{1}(x_2\ast^{1}x_2) \neq \mathcal L^{1}(x_2)\ast^{1}\mathcal L^{1}(x_2)\mspace{1mu}$. But the subspace $\bigoplus_{r\geq0}(x_1\cdots x_r)^{\epsilon}\Q[x_1,\ldots,x_r]$ is a subalgebra, and $\mathcal L^{\epsilon}$ is linear over it.
\begin{proposition}\label{prop:Wlinear}
Let $\epsilon\in\{0,1\}$ and $f\in\bigoplus_{r\geq0}(x_1\cdots x_r)^{\epsilon}\Q[x_1,\ldots,x_r]\mspace{1mu}$. Then $\mathcal L^{\epsilon}(p\ast^{\epsilon} f)=\mathcal L^{\epsilon}(p)\ast^{\epsilon} f$ for every~$p$ in the domain of~$\mathcal L^{\epsilon}$.
\end{proposition}

\begin{question}\label{q:formal}
Is there a proof of Theorem~\ref{thm:phi} and Proposition~\ref{prop:Wlinear} within the formal multiple Eisenstein series of \cite[Section~2]{BvI2}, in which swap invariance (the conjugation of partitions) and the stuffle product are imposed on formal symbols, or equivalently within Burmester's balanced setup \cite{Bu2} (see also \cite[Section~6]{BvI2}), without passing to $q$-series?
\end{question}

For $\bfd=(d_1,\dots,d_s)$ write $|\bfd|\df d_1+\dots+d_s\mspace{1mu}$. Each letter~$y$ produced by the recursion for $\EE_{\epsilon}$ lowers both $|\boldsymbol\ell|$ and~$|\bfk|$, which bounds the depth of the words. We record the bound because it is the form in which Bachmann stated his conjecture \cite[Conjecture~4.3]{Ba}.

\begin{theorem}\label{thm:depth}
For $\boldsymbol\ell,\bfk\in\Zgeq^{r}$ with $r\geq1\mspace{1mu}$, put $w\df|\boldsymbol\ell|+|\bfk|$ and $s_0\df\min(|\boldsymbol\ell|,|\bfk|)\mspace{1mu}$. Then
\begin{align*}
  \zeta^\dagger_q(0^{\ell_1-1},k_1,\dots,0^{\ell_r-1},k_r)\ &\in\ \bigl\langle\,\zeta^\dagger_q(\bfd)\bigm|1\leq s\leq s_0,\ \bfd\in\Zgeq^{s}\mspace{1mu},\ |\bfd|\leq w-r\,\bigr\rangle_{\ZZ}\,,\\[4pt]
  \zeta^{\mathrm{BZ}}_q(1^{\ell_1-1},k_1+1,\dots,1^{\ell_r-1},k_r+1)\ &\in\ \bigl\langle\,\zeta^{\mathrm{BZ}}_q(\bfd)\bigm|1\leq s\leq s_0,\ \bfd\in\ZZ_{\geq2}^{s}\mspace{1mu},\ |\bfd|\leq w\,\bigr\rangle_{\ZZ}\,,
\end{align*}
where $\zeta^{\mathrm{BZ}}_q$ denotes the Bradley--Zhao model of multiple~$q$-zeta values recalled in Section~\ref{sec:models}.
\end{theorem}

Theorem~\ref{thm:depth} is an effective version of Theorem~\ref{thm:phi}, since the series on the left span the left-hand sides of \eqref{eq:Zq} and the series on the right span the right-hand sides. The bound on the weight is implicit in \cite{HMW26}, where no bound on the depth is stated.

Let $\EE^{\mathrm{cl}}_{\epsilon}(\bfc)$ be the homogeneous part of $\EE_{\epsilon}(\bfc)$ of weight $|\bfc|-(1-\epsilon)r\mspace{1mu}$, the largest weight that occurs (see Corollary~\ref{cor:topweight} and Proposition~\ref{prop:classicalZ}). For $\epsilon=1$ these \emph{classical words} are the words~$\mathfrak D(\bfc)$ by which Hirose, Maesaka, Seki and Watanabe \cite{HMSW} write every multiple zeta value as a $\ZZ$-linear combination of those with all entries at least~$2\mspace{1mu}$,
\begin{equation}\label{eq:drop1}
  \zeta(1^{\ell_1-1},k_1+1,\dots,1^{\ell_r-1},k_r+1)\=Z\bigl(\EE^{\mathrm{cl}}_{1}(\ell_1,k_1,\dots,\ell_r,k_r)\bigr)\mspace{1mu},
\end{equation}
where $Z$ sends $yx^{l_1-1}\cdots yx^{l_s-1}$ to $\zeta(l_1,\dots,l_s)\mspace{1mu}$. The linear map~$\mathcal D$ sending $y^{\ell_1}x^{k_1}\cdots y^{\ell_r}x^{k_r}$ to $\EE^{\mathrm{cl}}_{1}(\ell_1,k_1,\dots,\ell_r,k_r)$ is their \emph{Drop1 operator}, which also describes the derivative of multiple Eisenstein series \cite[Theorem~D]{BKM26}. Bachmann, Kanno and Maesaka remark that a non-recursive expression for it seems difficult to find \cite[Section~4.1]{BKM26}. The following corollary of Theorem~\ref{thm:intro} provides one.

Put $\eta^{\mathrm{cl}}_D\df x^{2}\hat t_u\hat t_v$ for a domino $D=\{u,v\}\mspace{1mu}$, let $\Pp^{\mathrm{cl}}_S$ and $\Pp^{\mathrm{cl},*}_S$ be the products~$\Pp_S$ and $\Pp^{*}_S$ above with $\eta^{\mathrm{cl}}$ in place of~$\eta\mspace{1mu}$, and let $\mathsf C^{\mathrm{cl}}_r$ be the matrix with rows and columns indexed by $\calU(r)$ and entries \[
  \mathsf C^{\mathrm{cl}}_{S\setminus A,S}\df-\frac{w^{\mathrm{cl}}_{A,S}}{x}\,,\qquad
  w^{\mathrm{cl}}_{A,S}\df(-1)^{\mathrm{eo}(A)}\,x^{(\epsilon-1)\kappa(A)}\,\frac{\Pp^{\mathrm{cl}}_{S\setminus A}}{\Pp^{\mathrm{cl},*}_S}\,\eta^{\mathrm{cl}}_A-\delta_{A,\emptyset}\qquad(A\in\calT(S))\mspace{1mu},
\]
and set $\mathsf C^{\mathrm{cl}}_{T,S}\df0$ for all other~$T\mspace{1mu}$.

\begin{corollary}\label{cor:introcl}
For $r\geq1$ and $\epsilon\in\{0,1\}\mspace{1mu}$, \[
  \sum_{\bfc\in\Zgeq^{2r}}\EE^{\mathrm{cl}}_{\epsilon}(\bfc)\,\bft^{\bfc}\=\bigl[(1+y\,\mathsf C^{\mathrm{cl}}_r)^{-1}\bigr]_{\emptyset,[2r]}\,.
\]
\end{corollary}

The paper is organised as follows. In Section~\ref{sec:prelim} we recall the models of $q$-multiple zeta values and the recursion of Hirose, Maesaka and Watanabe, explain why Theorem~\ref{thm:phi} is equivalent to their theorem, prove Theorem~\ref{thm:depth}, recall the characterisation of the classical words by finite multiple zeta sums and define the maps~$\mathcal L^{\epsilon}$. In Section~\ref{sec:solving} we solve that recursion, which gives Theorem~\ref{thm:intro} and Corollary~\ref{cor:introcl}. The weights~$w_{A,S}$ are studied in Section~\ref{sec:formal}, where we obtain a formula for the individual words by expanding the inverse in Theorem~\ref{thm:intro}. Finally, in Section~\ref{sec:qseries} we prove Proposition~\ref{prop:Wlinear} and turn the solution into an infinite ordered product of matrices, whose expansions are identities between $q$-series, among them identities between the quasimodular forms of the first paragraph. First values of the words~$\EE_{\epsilon}$ and of the polynomials~$\mathcal L^{0}(x^{a}y^{b})$ are listed in Appendix~\ref{app:table}. 

\subsection*{Acknowledgements}
We thank Henrik Bachmann for interesting suggestions. The author's research is part of the project \emph{Modular forms and multiple zeta values in mathematical physics} with file number VI.Veni.242.073 of the research programme Veni, which is financed by the Dutch Research Council (NWO) under the grant \href{https://doi.org/10.61686/OJGCY95356}{\nolinkurl{doi.org/10.61686/OJGCY95356}}.

\subsection*{Use of AI}
During the preparation of this work, after the results had been obtained and a first version of the paper written, the author used Claude Code (Anthropic) to rewrite parts of the text and to write computer code verifying statements in the paper. After using this tool, the author reviewed and edited the content as needed and takes full responsibility for the content of the publication.

\section{Preliminaries}\label{sec:prelim}

Throughout, we write $\NN = \{0,1,2,\dots\}$ for the set of natural numbers and $\Zgeq=\{1,2,\dots\}$ for the positive integers. 

\subsection{\texorpdfstring{Models of $q$-MZV's}{Models of q-MZV's}}\label{sec:models}
In the $q$-MZV literature, authors often work with specific \emph{models} of $q$-multiple zeta values, that is, with distinguished spanning families of $\mathcal{Z}_q$ and of its subspaces. (They are spanning families rather than bases: the members satisfy many linear relations, such as the dualities recalled in Section~\ref{sec:HMWthm}.) We introduce the models needed in the sequel. Set $x_0=y_0=0$ and write
{\allowdisplaybreaks
\begin{align*}
\mathbb{I} &\defis \{(k_1,\ldots,k_r)\in \Zgeq^r \mid r\geq 0\}\,, \\
\mathbb{I}^{\mathrm{adm}} &\defis \{(k_1,\ldots,k_r)\in \mathbb{I} \mid r=0 \text{ or } k_r\geq 2\}\,, \\
\mathbb{I}^{\geq 2} &\defis \{(k_1,\ldots,k_r)\in \mathbb{I} \mid r=0 \text{ or } k_j\geq 2 \text{ for all } j\}\,, \\
\overline{\mathbb{I}} &\defis \{(c_1,\ldots,c_d)\in \N^d \mid d\geq 0\}\,,\\
\overline{\mathbb{I}}^{\mathrm{adm}} &\defis \{(c_1,\ldots,c_{d})\in \overline{\mathbb{I}} \mid d= 0 \text{ or } c_d\neq 0\}\,.
\end{align*}}%
(In \cite{HMW26}, the authors denote by $\overline{1}$ what we denote by $0\mspace{1mu}$, following Burmester \cite{Bu1}.)

Every $\boldsymbol{c}\in \overline{\mathbb{I}}^{\mathrm{adm}}$ can be written in exactly one way as
\begin{equation}\label{eq:ckl}
\boldsymbol{c} \= (0^{\ell_1-1},k_1,\ldots, 0^{\ell_r-1}, k_r)
\end{equation}
with $r\geq 0$ and $\boldsymbol{k},\boldsymbol{\ell} \in \Zgeq^r\mspace{1mu}$, where $(0^{\ell})$ denotes the sequence $(0,\ldots,0)$ of length~$\ell\mspace{1mu}$. 

\subsubsection*{Hirose--Maesaka--Watanabe model}  The polynomials
\begin{equation}\label{eq:pc}
  p_{\boldsymbol c}(\mathbf{x},\mathbf{y})\df\prod_{j=1}^r\binom{x_j+k_j-2}{k_j-1}\binom{y_j-y_{j-1}+\ell_j-2}{\ell_j-1}
\end{equation}
form a basis of $\bigoplus_{r\geq 0} \Q[x_1,\ldots,x_r,y_1,\ldots,y_r]\mspace{1mu}$, indexed by $\boldsymbol{c}=(0^{\ell_1-1},k_1,\ldots,0^{\ell_r-1},k_r)\in\overline{\mathbb{I}}^{\mathrm{adm}}$ as in~\eqref{eq:ckl}. Expanding each factor as a geometric series shows that $\Psi^{0}(p_{\boldsymbol c})=\zeta^\dagger_q(\boldsymbol{c})\mspace{1mu}$, where $\zeta^\dagger_q$ is given by~\eqref{eq:HMW-model}.
For $\boldsymbol c$ with all $c_j\geq1$ the polynomial~$p_{\boldsymbol c}$ involves the multiplicities alone,
\begin{equation}\label{eq:pcmult}
  p_{\boldsymbol c}\=\prod_{j=1}^{d}\binom{x_j+c_j-2}{c_j-1}\qquad(c_1,\dots,c_d\geq1)\mspace{1mu},
\end{equation}
and these $p_{\boldsymbol c}$ form a basis of $\bigoplus_{d\geq0}\Q[x_1,\dots,x_d]\mspace{1mu}$.

\subsubsection*{Bradley--Zhao model} The Bradley--Zhao model is given by \[
\zeta_q^{\mathrm{BZ}}(\boldsymbol{k})\defis  \sum_{0<n_1 < \ldots < n_r } \prod_{\substack{1\leq j \leq r}} \frac{q^{n_j(k_j-1)}}{(1-q^{n_j})^{k_j}}
\qquad(\boldsymbol{k}\in \mathbb{I}^{\mathrm{adm}})\mspace{1mu},
\]
the admissibility $k_r\geq2$ being needed for convergence in $\Q[\![q]\!]\mspace{1mu}$. Termwise this model differs from the previous by $q\mapsto q^{-1}$ and a sign, i.e.,
\begin{equation}\label{eq:BZdagger}
  (-1)^{k}\,\frac{q^{-n(k-1)}}{(1-q^{-n})^{k}}\=\frac{q^{n}}{(1-q^{n})^{k}}\,,
\end{equation}
where the right-hand side is the corresponding factor of $\zeta^\dagger_q\mspace{1mu}$. This observation is what relates the two models at the finite level in Section~\ref{sec:HMWthm}.

 We write \[
  \mathcal{Z}_{q,1}\defis\operatorname{Span}_\Q\bigl\{\zeta^{\mathrm{BZ}}_q(\boldsymbol k)\bigm|\boldsymbol k\in\mathbb{I}^{\mathrm{adm}}\bigr\}\,,
\]
which is the space denoted so in \cite[Lemma~2.4]{HMW26}. Expanding every factor by means of
\begin{equation}\label{eq:BZgeom}
  \frac{q^{n(k-1)}}{(1-q^{n})^{k}}\=\sum_{m\geq 0}\binom{m}{k-1}q^{mn}\qquad(k\geq 1)\mspace{1mu},
\end{equation}
where the term $m=0$ is present precisely when $k=1\mspace{1mu}$, we obtain for every $\boldsymbol{k}\in \mathbb{I}^{\mathrm{adm}}$
\begin{equation}\label{eq:BZPsi1}
  \zeta^{\mathrm{BZ}}_q(\boldsymbol k)\=\Psi^{1}_X\biggl(\prod_{j=1}^{r}\binom{x_j}{k_j-1}\biggr)\mspace{1mu},
\end{equation}
the right-hand side being defined because $k_r\geq2$ makes the product divisible by~$x_r\mspace{1mu}$. The polynomials $\prod_{j=1}^r\binom{x_j}{k_j-1}$ with $\boldsymbol{k}\in\mathbb{I}^{\mathrm{adm}}$ form a basis of $\bigoplus_{r\geq 0} x_r\Q[x_1,\ldots,x_r]\mspace{1mu}$, and those with $\boldsymbol{k}\in \mathbb{I}^{\geq 2}$ a basis of $\bigoplus_{r\geq 0} x_1\cdots x_r\Q[x_1,\ldots,x_r]\mspace{1mu}$. Hence
\begin{equation}\label{eq:Zq1model}
  \mathcal{Z}_{q,1}\=\im\Psi^{1}_X\,,
\end{equation}
which is the left-hand side of the second equality of \eqref{eq:Zq}. Moreover, for $\boldsymbol{k}\in \mathbb{I}^{\geq 2}$ all terms with some $m_j=0$ in~\eqref{eq:BZPsi1} vanish, so that
\begin{equation}\label{eq:BZgeq2}
  \zeta^{\mathrm{BZ}}_q(\boldsymbol k)\=\Psi^{0}_X\biggl(\prod_{j=1}^{r}\binom{x_j}{k_j-1}\biggr)
  \qquad(\boldsymbol k\in\mathbb{I}^{\geq2})\mspace{1mu},
\end{equation}
and the span of these is its right-hand side.

\subsubsection*{Words}\label{sec:setup} Let $\QB$ be the non-commutative polynomial algebra in the alphabet $\mathcal{B}=\{b_0,b_1,\ldots \}\mspace{1mu}$, and let $\mathcal{B}^\ast$ be the free monoid generated by $\mathcal{B}\mspace{1mu}$, whose elements we call \emph{words}.
Moreover, with $\fh=\ZZ\langle x,y\rangle$ as in the introduction, consider the non-commutative polynomial algebras \[
\mathfrak{h}^1\defis \ZZ+y\mathfrak{h}\,, \quad \mathfrak{h}^0\defis \ZZ+y\mathfrak{h}x\,, \quad \mathfrak{h}^{\geq 2} \defis \ZZ+yx\ZZ\langle x, yx\rangle\,,
\]
so that $\mathfrak{h}^{\geq2}\subset\mathfrak{h}^0\subset\mathfrak{h}^1$. Below we write $\mathfrak{h}^{1+\epsilon}$ for $\mathfrak{h}^1$ if $\epsilon=0$ and for $\mathfrak{h}^{\geq2}$ if $\epsilon=1\mspace{1mu}$, so as to treat the two cases at once. We consider the $\ZZ$-modules $\mathfrak{h}^1, \mathfrak{h}^0$ and $\mathfrak{h}^{\geq2}$ to be contained in $\QB$ via \[
yx^{k_1-1}\cdots yx^{k_r-1} \mapsto b_{k_1}\cdots b_{k_r}
\]
and define linear maps $Z_q^\dagger:\mathfrak{h}^1 \to \Q[\![q]\!]$ and $Z_q^{\mathrm{BZ}}: \mathfrak{h}^0 \to \Q[\![q]\!]$ by $Z_q^\dagger(1):=1=:Z_q^{\mathrm{BZ}}(1)$ and \[
Z_q^\dagger(b_{c_1}\cdots b_{c_r}) \defis \zeta_q^\dagger(c_1,\ldots,c_r)\mspace{1mu}, \qquad Z_q^{\mathrm{BZ}}(b_{k_1}\cdots b_{k_r}) \defis \zeta_q^{\mathrm{BZ}}(k_1,\ldots,k_r)\mspace{1mu}.
\]
Since $\Psi^{0}(p_{\boldsymbol c})=\zeta^\dagger_q(\boldsymbol c)$ and, by \eqref{eq:pcmult}, $p_{\boldsymbol c}$ involves the multiplicities alone when $\boldsymbol c\geq1\mspace{1mu}$, we have $Z^\dagger_q=\Psi^{0}_X\circ\rho$ for the isomorphism \[\rho\colon\fh^1\otimes\Q\to\bigoplus_{r\geq0}\Q[x_1,\ldots,x_r]\,, \qquad
b_{c_1}\cdots b_{c_r}\mapsto\prod_{j}\binom{x_j+c_j-2}{c_j-1}\,.\]

\subsection{The stuffle product on polynomials}\label{sec:products}
For $f\in\Q[x,y]$ and $p\in\Q[x_1,\ldots,x_{r-1},y_1,\ldots,y_{r-1}]$ put \[
  f\otimes p\df f(x_1,y_1)\,p(x_2,\ldots,x_r,y_2,\ldots,y_r)\;\in\;\Q[x_1,\ldots,x_{r},y_1,\ldots,y_{r}]\,,
\]
so that every element of $\Q[x_1,\ldots,x_r,y_1,\ldots,y_r]$ with $r\geq1$ is a sum of such products. For $f,g\in\Q[x,y]$ put
\[
  (f\diamond^{\epsilon}g)(x,y)\df(f\diamond g)(x,y)\df\sum_{a=1-\epsilon}^{x-1+\epsilon}f(a,y)\,g(x-a,y)\mspace{1mu},
\]
which is a polynomial in $x$ and~$y\mspace{1mu}$. This product is commutative and associative. Let $\ast^{\epsilon}\df\ast$ be the corresponding quasi-shuffle (or stuffle) product on $\bigoplus_{r\geq0}\Q[x_1,\ldots,x_r,y_1,\ldots,y_r]$ (see \cite{H} and \cite[Section~2]{HI}), that is, the bilinear product with unit~$1$ such that
\[
  (f\otimes p)\ast(g\otimes p')\=f\otimes\bigl(p\ast(g\otimes p')\bigr)\+g\otimes\bigl((f\otimes p)\ast p'\bigr)\+(f\diamond g)\otimes(p\ast p')
\]
for $f,g\in\Q[x,y]\mspace{1mu}$, $p\in\Q[x_1,\ldots,x_{r-1},y_1,\ldots,y_{r-1}]$ and $p'\in\Q[x_1,\ldots,x_{r'-1},y_1,\ldots,y_{r'-1}]\mspace{1mu}$. For instance, $y_1\ast^{0}y_1=2\,y_1y_2+(x_1-1)\,y_1^{2}\mspace{1mu}$.

\begin{proposition}\label{prop:star}
For $\epsilon\in\{0,1\}$ the product $\ast^{\epsilon}$ makes $\Q\oplus\bigoplus_{r\geq1}x_r^{\epsilon}\Q[x_1,\ldots,x_r,y_1,\ldots,y_r]$ a commutative and associative algebra with unit~$1\mspace{1mu}$, in which $\Q\oplus\bigoplus_{r\geq1}x_r^{\epsilon}\Q[x_1,\ldots,x_r]$ and $\Q\oplus\bigoplus_{r\geq1}x_1\cdots x_r\Q[x_1,\ldots,x_r]$ are subalgebras, and $\Psi^{\epsilon}$ is a homomorphism of algebras. On the latter, the products $\ast^{0}$ and~$\ast^{1}$ agree.
\end{proposition}

\begin{proof}
For $\epsilon=0$ this is in \cite[Section~3]{BvI}, where $\Psi^{0}$ with variables~$u_i$ in place of~$q^i$ is the injective map~$\Phi$ (up to the order of the variables and the roles of $x$ and~$y$).

The following argument covers both values of~$\epsilon\mspace{1mu}$. By \cite[Theorem~2.1]{HI}, with the monomials~$x^ky^\ell$ in the role of the letters, the product~$\ast$ is commutative and associative. By induction on the number of variables, $(f_1\otimes\dots\otimes f_r)\ast(g_1\otimes\dots\otimes g_{r'})$ is a sum of products $h_1\otimes\dots\otimes h_t$ in which $h_1,\dots,h_t$ interleave $f_1,\dots,f_r$ and $g_1,\dots,g_{r'}\mspace{1mu}$, with some pairs $f_i,g_{i'}$ replaced by $f_i\diamond g_{i'}\mspace{1mu}$; in particular~$h_t$ is $f_r\mspace{1mu}$, $g_{r'}$ or $f_r\diamond g_{r'}\mspace{1mu}$. For $f,g$ divisible by~$x$ the terms $a=0$ and $a=x$ of $f\diamond^{1}g$ vanish, so that $f\diamond^{1}g=f\diamond^{0}g\mspace{1mu}$, and this polynomial is divisible by~$x\mspace{1mu}$. Hence the three spaces in the statement are closed under~$\ast\mspace{1mu}$, and $\ast^{0}$ and~$\ast^{1}$ agree on the last of them.

For $p=f_1\otimes\dots\otimes f_r$ the sum~\eqref{eq:phi} reads \[
  \Psi^{\epsilon}(p)\=\sum_{0<n_1<\dots<n_r}\ \prod_{j=1}^{r}\psi_{f_j}(n_j)\mspace{1mu},\qquad
  \psi_f(n)\df\sum_{m\geq1-\epsilon}f(m,n)\,q^{mn}\,,
\]
and $\psi_f(n)\,\psi_g(n)=\psi_{f\diamond g}(n)\mspace{1mu}$, since the sum of $f(m,n)\,g(m',n)$ over $m,m'\geq1-\epsilon$ with $m+m'=M$ is $(f\diamond g)(M,n)\mspace{1mu}$. Split the sum for $\Psi^{\epsilon}(f\otimes p)\,\Psi^{\epsilon}(g\otimes p')$ according to whether the smallest part of the first factor is smaller than, larger than or equal to the smallest part of the second. By induction on the number of variables, applied to the sums over parts larger than a fixed integer, the three cases give the three terms in the definition of~$\ast\mspace{1mu}$, so that $\Psi^{\epsilon}(f\otimes p)\,\Psi^{\epsilon}(g\otimes p')=\Psi^{\epsilon}\bigl((f\otimes p)\ast(g\otimes p')\bigr)\mspace{1mu}$.
\end{proof}

\subsection{Subsets of integers}
\label{subsec:notation}
For $n\in\NN$ we write $[n]\defis \{1,\ldots,n\}\mspace{1mu}$. For $\bfc=(c_{1},\ldots,c_{n})\in\Zgeq^{n}\mspace{1mu}$, let \[
  \{\bfc=1\}\defis\{i\in[n]\mid c_{i}=1\}\,,\qquad
  \{\bfc>1\}\defis\{i\in[n]\mid c_{i}>1\}\=[n]\setminus \{\bfc=1\}\,.
\]
If $\bfc\in\ZZ^n$ and $A,B\subseteq[n]$ are disjoint, we write $\bfc_B\in\ZZ^n$ for the sequence obtained from $\bfc$ by subtracting~$1$ from the entries indexed by~$B\mspace{1mu}$; if all entries of $\bfc$ indexed by~$A$ are equal to~$1\mspace{1mu}$, we write $\bfc^{A}\in\ZZ^{n-|A|}$ for the sequence obtained from $\bfc$ by deleting them. We put $\bfc^{A}_{B}\defis(\bfc_B)^{A}$. Any term of a sum below in which $\bfc^{A}$ occurs while some entry of $\bfc$ indexed by~$A$ differs from~$1$ is to be read as zero.

 For a finite subset $B \subseteq \Zgeq\mspace{1mu}$, let $\kappa(B) \defis \frac{1}{2}|B|$ and let $\sigma_{B}$ be the unique order-preserving bijection $\NN\setminus B \to \NN\mspace{1mu}$, i.e., $ \sigma_{B}(i) = \bigl|[i]\setminus B\bigr|\mspace{1mu}$. Moreover, for $i\in \NN\mspace{1mu}$, we call $D_i\defis\{i,i+1\}$ a \emph{domino} and define the counting functions
{\allowdisplaybreaks
\begin{alignat*}{2}
\mathrm{eo}(B)&\defis |\{i \in 2\NN \mid D_i \subseteq B\}|\,,&\qquad \mathrm{oe}(B)&\defis |\{\,i \in 2\NN+1 \mid D_i \subseteq B\}|\,, \\
\alpha(B)&\defis |\{\,i\in2\NN\mid B\cap D_i\neq\emptyset\,\}|\,,&\qquad \beta(B)&\defis |\{\,i\in 2\NN+1\mid B\cap D_i\neq\emptyset\,\}|\,,
\end{alignat*}}%
that is, $\alpha(B)$ and $\beta(B)$ are the number of dominoes~$D_i$ needed to cover~$B\mspace{1mu}$, when restricting to even and odd~$i$ respectively. Moreover, $\alpha(B)=|B|-\mathrm{eo}(B)$ and $\beta(B)=|B|-\mathrm{oe}(B)\mspace{1mu}$. Note that $i=0$ is allowed here, so that the domino $D_0=\{0,1\}\mspace{1mu}$, which meets~$[n]$ in the single position~$1\mspace{1mu}$, is counted by $\alpha\mspace{1mu}$, e.g., $\alpha(\{1\})=1\mspace{1mu}$.
Further, we let $\boldsymbol{e}_n(i)=\boldsymbol{e}(i)$ be the $i$th standard basis vector in $\QQ^n$ and let $\boldsymbol{e}(B)\defis\sum_{i\in B}\boldsymbol{e}(i)$ for $B\subseteq[n]$ denote the \emph{characteristic vector} of $B\mspace{1mu}$.

 Let $S=\{s_1<\ldots<s_n\}\subset\Zgeq$ be finite. Recall that for $i\in[n-1]$ we call $D_i(S)\defis\{s_i,s_{i+1}\}$ a domino of~$S\mspace{1mu}$, which is \emph{even--odd} if $i$ is even and \emph{odd--even} if $i$ is odd, and we call two dominoes of~$S$ \emph{adjacent} if they occupy four consecutive positions of~$S\mspace{1mu}$.
 We define the following subsets of~$S\mspace{1mu}$:
{\allowdisplaybreaks
\begin{align}\nonumber
\mathrm{EO}(S)
  &\defis \Bigl\{\;\bigcup_{i \in R} D_i(S) \; \bigm|\; R \,\subseteq\, [n-1]\cap 2\ZZ\Bigr\} \,,\\\nonumber
\mathrm{OE}(S)
  &\defis \Bigl\{\;\bigcup_{i \in R} D_i(S) \; \bigm|\; R \,\subseteq\, [n-1] \cap (2\ZZ+1)\Bigr\}\,,\\\nonumber
\mathrm{OE}^{*}(S)
  &\defis \Bigl\{\;\bigcup_{i \in R} D_i(S) \; \bigm|\; R \,\subseteq\, [n-1] \cap (2\ZZ+1),\
       |i-i'| \neq 2  \text{ for all } i,i' \in R\Bigr\}\,,\\
\calT(S)
  &\defis \bigl\{\,A \cup A^* \ \bigm|\ 
       A \in \mathrm{EO}(S),\ A^* \in \mathrm{OE}^{*}(S),\
       A \cap A^* = \emptyset\bigr\}\,, \label{def:domino-tiling}
\end{align}}%
where $\calT(S)$ is the set of partial domino tilings of~$S\mspace{1mu}$, already defined in the introduction.
Observe that every $A\in\calT(S)$ is a disjoint union of dominoes of~$S$ in exactly one way. 
For a tiling $A\in\calT(S)$ we write $\mathrm{eo}(A)$ and $\mathrm{oe}(A)$ for the numbers of its even--odd and odd--even dominoes, as in the introduction, so that $\kappa(A)=\mathrm{eo}(A)+\mathrm{oe}(A)\mspace{1mu}$. For $S=[n]$ one has $D_i([n])=D_i\mspace{1mu}$, and for $A\in\calT([n])$ the number~$\mathrm{eo}(A)$ agrees with the display above: an even--odd domino $D_i\subseteq A$ is one of the dominoes of~$A\mspace{1mu}$, as otherwise $D_{i-1}$ and $D_{i+1}$ would be adjacent odd--even dominoes of~$A\mspace{1mu}$. For $\mathrm{oe}(A)$ this need not hold, as $A=D_2\cup D_4$ shows.

 We write $\calT(r)\defis\calT([2r])\mspace{1mu}$. For instance $\calT(1)=\{\emptyset,D_1\}$ and $\calT(2)=\{\emptyset,D_1,D_2,D_3\}\mspace{1mu}$, and $\calT(r')\subseteq\calT(r)$ for $r'<r\mspace{1mu}$. The numbers~$|\calT(r)|\mspace{1mu}$, equal to $2,4,10,24,58,140$ for $r=1,\dots,6\mspace{1mu}$, are twice the Pell numbers, that is $|\calT(r)|=2|\calT(r-1)|+|\calT(r-2)|\mspace{1mu}$.

\subsection{The Hirose--Maesaka--Watanabe Theorem}\label{sec:HMWthm}
 We recall the theorem that proves Bachmann's conjectures. The main construction is a pair of maps $\EE_{\epsilon}\colon \bigcup_{r\geq 0}\ZZ^{2r} \to \mathfrak{h}^{1+\epsilon}$ for $\epsilon\in\{0,1\}\mspace{1mu}$, written~$\mathfrak E^\epsilon_q$ in \cite{HMW26}:
\begin{definition}
\label{def:Eq}
For $r\geq 0, \bfc=(c_{1},\dots,c_{2r})\in\Zgeq^{2r}$ and $\epsilon\in\{0,1\}\mspace{1mu}$, define $\EE_{\epsilon}(\bfc)\in\fh^{1+\epsilon}$ recursively by $\EE_{\epsilon}(\emptyset)=1$ and, for $r\geq 1$
\begin{align*}
\EE_{\epsilon}(\bfc)
  \defis & -\!\!\!\sum_{\substack{\emptyset\neq B\subseteq\{\bfc>1\}}}\!\!\!\!\!
      (-1)^{|B|}\,\EE_{\epsilon}(\bfc_{B})\,x^{\alpha(B)}
   \+\!\!\sum_{\substack{B\subseteq \{\bfc >1\}}}\!\!
      (-1)^{|B|}\,\EE_{\epsilon}(\bfc_{B})
      \left(\sum_{h=\alpha(B)+1}^{\beta(B)}\!\! yx^{h-1} \meno\!\! \sum_{h=\beta(B)+1}^{\alpha(B)}\!\! yx^{h-1}\right)\\
  &\quad
   \+\sum_{\substack{\emptyset \neq A\in \calT(r): A\subseteq\{\bfc=1\}\\
                    B\subseteq \{\bfc>1\}}}
    \sum_{h=1}^{\kappa(A)}
      (-1)^{\mathrm{eo}(A)+|B|+\kappa(A)-h}
      \binom{\kappa(A)-1}{h-1}
      \EE_{\epsilon}(\bfc^A_B)\,
      yx^{\epsilon\kappa(A)+\beta(\sigma_{A}(B))+h-1}.
\end{align*}
Extend by $\EE_{\epsilon}(\bfc)=0$ if $c_i\leq 0$ for some $i\mspace{1mu}$. 
\end{definition}
Note that at least one of the two sums over $h$ in the second sum is empty. 

Next, we introduce the finite version of the Hirose--Maesaka--Watanabe model of multiple~$q$-zeta values. 
For $N\geq1$ and $\bfc=(c_1,\dots,c_d)\in\overline{\mathbb{I}}^{\mathrm{adm}}$ let \[
  \zeta^\dagger_{q,N}(\bfc)\df\sum_{\substack{0<n_1\leq\dots\leq n_d<N\\
  n_j<n_{j+1}\text{ if }c_j\neq0}}\ \prod_{j\,:\,c_j=0}\frac{1}{1-q^{N-n_j}}\
  \prod_{j\,:\,c_j\neq0}\frac{q^{n_j}}{(1-q^{n_j})^{c_j}}\,,
\]
which tends to $\zeta^\dagger_q(\bfc)$ as $N\to\infty\mspace{1mu}$, and for $\bfc$ with all $c_j\geq1$ let \[
  \zeta^{\diamond}_{q,N}(\bfc)\df\sum_{A\subseteq\{j\,:\,c_j=1\}}\
  \sum_{\substack{0<n_1\leq\dots\leq n_d<N\\ n_j<n_{j+1}\text{ if }j\notin A}}\
  \prod_{j\in A}\frac{1}{1-q^{N-n_j}}\ \prod_{j\notin A}\frac{q^{n_j}}{(1-q^{n_j})^{c_j}}\,.
\]
Following \cite{HMW26}, we treat the two cases at once: for $\bfc=(\ell_1,k_1,\dots,\ell_r,k_r)\in\Zgeq^{2r}$ put
\begin{equation}\label{eq:xidef}
  \xi^{\epsilon}_N(\bfc;q)\df\xi_N(\bfc;q)\df\begin{cases}
  \zeta^\dagger_{q,N}\bigl(0^{\ell_1-1},k_1,\dots,0^{\ell_r-1},k_r\bigr)&\text{if }\epsilon=0\mspace{1mu},\\[2pt]
  \zeta^{\diamond}_{q,N}\bigl(1^{\ell_1-1},k_1+1,\dots,1^{\ell_r-1},k_r+1\bigr)&\text{if }\epsilon=1\mspace{1mu}.
  \end{cases}
\end{equation}
Write $\xi^{\epsilon}(\bfc;q)\df\xi(\bfc;q)\df\lim_{N\to\infty}\xi_N(\bfc;q)\mspace{1mu}$.  Let $Z^\dagger_{q,N}\colon\fh^1\to\ZZ[\![q]\!]$ be the linear map with $Z^\dagger_{q,N}(1)=1$ and \[
  Z^\dagger_{q,N}\bigl(yx^{c_1-1}\cdots yx^{c_s-1}\bigr)\=\zeta^\dagger_{q,N}(c_1,\dots,c_s)\mspace{1mu}.
\]
 Then \cite[Theorem~3.7]{HMW26} gives the following.

\begin{theorem}\label{thm:HMW}
For $N\geq1\mspace{1mu}$, $r\geq0\mspace{1mu}$, $\epsilon\in\{0,1\}$ and $\bfc\in\Zgeq^{2r}\mspace{1mu}$,
\begin{equation}\label{eq:HMWfinite}
  \xi^{\epsilon}_{N}(\bfc;q)\=Z^\dagger_{q,N}\bigl(\EE_{\epsilon}(\bfc)\bigr)\mspace{1mu},
\end{equation}
and, in the limit,
\begin{equation}\label{eq:HMWinf}
  \xi^{\epsilon}(\bfc;q)\=Z^\dagger_{q}\bigl(\EE_{\epsilon}(\bfc)\bigr)\mspace{1mu}.
\end{equation}
\end{theorem}

\begin{remark}\label{rem:BZform}
For $\epsilon=0\mspace{1mu}$, \eqref{eq:HMWinf} reads $\zeta^\dagger_q(0^{\ell_1-1},k_1,\dots,0^{\ell_r-1},k_r)=Z^\dagger_q(\EE_{0}(\bfc))\mspace{1mu}$, which is the first identity of \cite[Theorem~1.1]{HMW26}. Their second identity is the case $\epsilon=1\mspace{1mu}$, transported to the Bradley--Zhao model: putting \[
  \DD(\bfc)\defis(-1)^{c_{1}+\cdots+c_{2r}}\,\theta\bigl(\EE_{1}(\bfc)\bigr)
  \,\in\,\fh^{\geq2}\,,
\]
with $\theta$ the $\ZZ$-algebra automorphism of $\fh$ sending $x\mapsto-x\mspace{1mu}$, $y\mapsto-y\mspace{1mu}$, one has \[
  \zeta_q^{\mathrm{BZ}}(1^{\ell_1-1},k_1+1,\dots,1^{\ell_r-1},k_r+1)
  \= Z_q^{\mathrm{BZ}}\bigl(\DD(\bfc)\bigr)\mspace{1mu}.
\]
This follows from \eqref{eq:HMWfinite} for $\epsilon=1$ after substituting $q\mapsto q^{-1}$, both sides of which are rational functions of~$q\mspace{1mu}$, and letting $N\to\infty\mspace{1mu}$: the signs from \eqref{eq:BZdagger} cancel against those in $\DD\mspace{1mu}$.
\end{remark}

\begin{proof}[Proof of Theorem~\ref{thm:phi}]The series~$\zeta^\dagger_q(\bfc)$ span~$\im\Psi^{0}$, and those with all $c_j\geq1$ span~$\im\Psi^{0}_X\mspace{1mu}$. By \eqref{eq:HMWinf} for $\epsilon=0\mspace{1mu}$, every $\zeta^\dagger_q(\bfc)$ is $Z^\dagger_q$ of an element of~$\fh^1$, hence lies in $\im\Psi^{0}_X$ by $Z^\dagger_q=\Psi^{0}_X\circ\rho\mspace{1mu}$. This is the first equality, the inclusion $\supseteq$ being trivial. For the second, $\im\Psi^{1}_X$ is spanned by the $\zeta^{\mathrm{BZ}}_q(\bfk)$ with $\bfk\in\mathbb I^{\mathrm{adm}}$ by \eqref{eq:Zq1model}, each of which is $Z^{\mathrm{BZ}}_q$ of an element of $\fh^{\geq2}$ by Remark~\ref{rem:BZform}, hence lies in $\Psi^{0}_X\bigl(\bigoplus_{r\geq0}x_1\cdots x_r\Q[x_1,\ldots,x_r]\bigr)$ by \eqref{eq:BZgeq2}; the reverse inclusion is \eqref{eq:BZgeq2} again.
\end{proof}

The map $u\mapsto(Z^\dagger_{q,N}(u))_{N\geq1}$ from $\fh^1$ to $\ZZ[\![q]\!]^{\NN}$ is injective by \cite[Lemma~2.6]{HMW26}, so that \eqref{eq:HMWfinite} characterises the elements~$\EE_{\epsilon}(\bfc)$ uniquely. The words are invariant under reversal of~$\bfc\mspace{1mu}$,
\begin{equation}\label{eq:reflection}
  \EE_{\epsilon}(c_1,\dots,c_{2r})\=\EE_{\epsilon}(c_{2r},\dots,c_1)\mspace{1mu},
\end{equation}
as noted in \cite[Remark~1.4]{HMW26}. By \eqref{eq:HMWinf} this gives the duality $\zeta^\dagger_q(0^{\ell_1-1},k_1,\dots,0^{\ell_r-1},k_r) =\zeta^\dagger_q(0^{k_r-1},\ell_r,\dots,0^{k_1-1},\ell_1)$ and its counterpart for the Bradley--Zhao model.

\subsection{Weight and depth}\label{sec:classical}
The \emph{weight} $\wt(w)$ of a word $w\in\fh$ is its number of letters, its \emph{depth} $\dep(w)$ the number of letters~$y\mspace{1mu}$; we call~$r$ the depth of $\bfc\in\Zgeq^{2r}\mspace{1mu}$.

\begin{corollary}\label{cor:topweight}
For $\bfc\in\Zgeq^{2r}$ every word~$w$ occurring in $\EE_{\epsilon}(\bfc)$ has
\begin{equation}\label{eq:weightbound}
  \wt(w)\ \leq\ |\bfc|-(1-\epsilon)r
\end{equation}
and
\begin{equation}\label{eq:depthbound}
  \dep(w)\ \leq\ \min(c_1+c_3+\dots+c_{2r-1},\;c_2+c_4+\dots+c_{2r})\mspace{1mu}.
\end{equation}
\end{corollary}

\begin{proof}
Every index at which $\EE_\epsilon$ occurs on the right-hand side of the recursion in Definition~\ref{def:Eq} has a smaller sum of entries than~$\bfc$ (the term $B=\emptyset$ of the second sum vanishes), so both bounds follow by induction on $|\bfc|=c_1+\dots+c_{2r}\mspace{1mu}$. For the weight, it suffices that each term, a multiple of $\EE_\epsilon(\bfc')\,u$ with $\bfc'\in\Zgeq^{2r'}$ and a word~$u\mspace{1mu}$, has $\wt(u)\leq|\bfc|-|\bfc'|-(1-\epsilon)(r-r')\mspace{1mu}$. In the first sum $\alpha(B)\leq|B|$ letters are added and $|\bfc|$ drops by~$|B|\mspace{1mu}$. In the second sum $h\leq\max(\alpha(B),\beta(B))\leq|B|$ letters are added. In the third sum $\epsilon\kappa(A)+\beta(\sigma_A(B))+h\leq(1+\epsilon)\kappa(A)+|B|$ letters are added, while $|\bfc|-|\bfc^A_B|-(1-\epsilon)\kappa(A)=(1+\epsilon)\kappa(A)+|B|\mspace{1mu}$.

For the depth, every term carrying a letter~$y$ lowers both $c_1+c_3+\dots+c_{2r-1}$ and $c_2+c_4+\dots+c_{2r}$ by at least~$1\mspace{1mu}$. In the second sum the sums over~$h$ are empty unless $\alpha(B)\neq\beta(B)\mspace{1mu}$, and as $\alpha(B)-\beta(B)=\mathrm{oe}(B)-\mathrm{eo}(B)$ this requires~$B$ to contain a domino~$D_i\mspace{1mu}$, whose two elements have different parities. In the third sum a nonempty tiling~$A$ is deleted, and every domino of~$A$ contains one position of each parity.
\end{proof}

Recall from the introduction that the classical word~$\EE^{\mathrm{cl}}_{\epsilon}(\bfc)$ is the homogeneous part of $\EE_{\epsilon}(\bfc)$ of weight $|\bfc|-(1-\epsilon)r\mspace{1mu}$. Let $Z_N\colon\fh^1\to\Q$ be the linear map with $Z_N(1)=1$ and \[
  Z_N(b_{k_1}\cdots b_{k_s})\df\sum_{0<n_1<\dots<n_s<N}\frac{1}{n_1^{k_1}\cdots n_s^{k_s}}\,,
\]
and let $\xi^{\mathrm{cl}}_N(\bfc)$ be the finite sum~$\xi^{\epsilon}_N(\bfc;q)$ with every factor $q^{n}/(1-q^{n})^{c}$ replaced by $n^{-c}$ and every factor $1/(1-q^{N-n})$ by $1/(N-n)\mspace{1mu}$, so that \[
  \xi^{\mathrm{cl}}_N(\bfc)\=\lim_{q\to1}(1-q)^{|\bfc|-(1-\epsilon)r}\,\xi^{\epsilon}_N(\bfc;q)\mspace{1mu}.
\]
For $\epsilon=0$ this is the sum $\zeta_N\binom{\boldsymbol\ell}{\bfk}$ of \cite[Remark~1.3]{HMW26}, for $\epsilon=1$ the diamond multiple zeta value $\zeta^{\diamond}_N(1^{\ell_1-1},k_1+1,\dots,1^{\ell_r-1},k_r+1)$ of \cite[Definition~2.3]{HMSW}.

\begin{proposition}[{\cite[Corollary~3.8]{HMW26}}]\label{prop:classicalZ}
For $\bfc\in\Zgeq^{2r}\mspace{1mu}$, $\epsilon\in\{0,1\}$ and $N\geq1$ we have $Z_N(\EE^{\mathrm{cl}}_{\epsilon}(\bfc))=\xi^{\mathrm{cl}}_N(\bfc)\mspace{1mu}$.
\end{proposition}

Since $u\mapsto(Z_N(u))_{N\geq1}$ is injective on~$\fh^1$ by \cite[Lemma~2.1]{HMSW}, Proposition~\ref{prop:classicalZ} characterises the classical words. It identifies them with the words~$\mathfrak E^{\epsilon}(\bfc)$ of \cite[Corollary~3.8]{HMW26}, which for $\epsilon=1$ are the words~$\mathfrak D(\bfc)$ of \cite[Theorem~2.5]{HMSW}. So $\EE^{\mathrm{cl}}_{1}$ is the Drop1 operator $\mathcal D\colon\fh^0\to\fh^{\geq2}$ of the introduction, and letting $N\to\infty$ gives~\eqref{eq:drop1} (see \cite[Theorem~1.8]{HMSW}).

\begin{corollary}\label{cor:mzvdepth}
Every multiple zeta value of weight~$w$ and depth~$d$ is a $\ZZ$-linear combination of multiple zeta values $\zeta(l_1,\dots,l_s)$ with $l_1,\dots,l_s\geq2\mspace{1mu}$, $l_1+\dots+l_s=w$ and $s\leq\min(d,w-d)\mspace{1mu}$.
\end{corollary}

\begin{proof}
Write the index as $(1^{\ell_1-1},k_1+1,\dots,1^{\ell_r-1},k_r+1)$ and put $\bfc=(\ell_1,k_1,\dots,\ell_r,k_r)\mspace{1mu}$, so that $d=\ell_1+\dots+\ell_r$ and $w=|\bfc|\mspace{1mu}$. By \eqref{eq:drop1} the value is $Z(\EE^{\mathrm{cl}}_{1}(\bfc))\mspace{1mu}$, and $\EE^{\mathrm{cl}}_{1}(\bfc)\in\fh^{\geq2}$ is homogeneous of weight~$w$ and satisfies \eqref{eq:depthbound}.
\end{proof}

\begin{proof}[Proof of Theorem~\ref{thm:depth}]
Let $\bfc=(\ell_1,k_1,\dots,\ell_r,k_r)\mspace{1mu}$, so that $|\bfc|=w\mspace{1mu}$. By Corollary~\ref{cor:topweight} the words of $\EE_{\epsilon}(\bfc)$ have weight at most $w-(1-\epsilon)r$ and depth at most~$s_0\mspace{1mu}$, and $\EE_{\epsilon}(\bfc)$ has no constant term, since for $N=1$ both sides of \eqref{eq:HMWfinite} vanish except for $Z^\dagger_{q,1}(1)=1\mspace{1mu}$. For $\epsilon=0\mspace{1mu}$, \eqref{eq:HMWinf} gives $\zeta^\dagger_q(0^{\ell_1-1},k_1,\dots,0^{\ell_r-1},k_r)=Z^\dagger_q(\EE_{0}(\bfc))$ with $Z^\dagger_q(b_{d_1}\cdots b_{d_s})=\zeta^\dagger_q(d_1,\dots,d_s)\mspace{1mu}$, which is the first statement. For $\epsilon=1\mspace{1mu}$, Remark~\ref{rem:BZform} gives $\zeta^{\mathrm{BZ}}_q(1^{\ell_1-1},k_1+1,\dots,1^{\ell_r-1},k_r+1)=Z^{\mathrm{BZ}}_q(\DD(\bfc))$ with $\DD(\bfc)\in\fh^{\geq2}$, and $\theta$ preserves weight and depth, which is the second.
\end{proof}

\subsection{\texorpdfstring{The maps $\mathcal L^{\epsilon}$}{The maps L}}\label{sec:Lmaps}
The maps~$\mathcal L^{\epsilon}$ of the introduction are defined as follows. The map~$\mathcal L^{0}$ sends the polynomial~$p_{\bfc}\mspace{1mu}$, with $\bfc=(0^{\ell_1-1},k_1,\dots,0^{\ell_r-1},k_r)\mspace{1mu}$, to $\rho\bigl(\EE_{0}(\ell_1,k_1,\dots,\ell_r,k_r)\bigr)\mspace{1mu}$. By \eqref{eq:HMWinf} and the identities $\Psi^{0}(p_{\bfc})=\zeta^\dagger_q(\bfc)$ and $Z^\dagger_q=\Psi^{0}_X\circ\rho$ it then follows that $\Psi^{0}=\Psi^{0}_X\circ\mathcal L^{0}$. For example, one has $\mathcal L^{0}(y_1)=x_1\mspace{1mu}$, since $p_{(0,1)}=y_1$ and $\EE_{0}(2,1)=b_2\mspace{1mu}$, and $\mathcal L^{0}(x_1y_1)=\binom{x_1+1}{2}\oplus(-x_2)$ by \eqref{eq:derg1}, where $f\oplus g$ denotes the element of $\Q[x_1]\oplus\Q[x_1,x_2]$ with components~$f$ and~$g\mspace{1mu}$.

The map~$\mathcal L^{1}$ is defined by sending $\prod_j\binom{x_j}{k'_j-1}\mspace{1mu}$, with $\bfk'=(1^{\ell_1-1},k_1+1,\dots,1^{\ell_r-1},k_r+1)\mspace{1mu}$, to the image of $\DD(\ell_1,k_1,\dots,\ell_r,k_r)\in\fh^{\geq2}$ under $b_{d_1}\cdots b_{d_s}\mapsto\prod_j\binom{x_j}{d_j-1}\mspace{1mu}$. By \eqref{eq:BZPsi1}, Remark~\ref{rem:BZform} and \eqref{eq:BZgeq2} it then follows that $\Psi^{1}_X=\Psi^{0}_X\circ\mathcal L^{1}$.  For example, one has $\mathcal L^{1}(x_2)=\binom{x_1}{2}\mspace{1mu}$, as $\DD(2,1)=\DD(1,2)=b_3$ by \eqref{eq:reflection} and Lemma~\ref{lem:Lprops}\ref{it:Lprops-id}. This is the identity $\zeta^{\mathrm{BZ}}_q(1,2)=\zeta^{\mathrm{BZ}}_q(3)\mspace{1mu}$, a $q$-analogue of $\zeta(2,1)=\zeta(3)\mspace{1mu}$.

Conjugation of partitions, as in Figure~\ref{fig:conjugate}, is the substitution \[
  \iota\colon\quad x_j\mapsto y_{r+1-j}-y_{r-j}\,,\qquad y_j\mapsto x_{r+1-j}+\dots+x_r\qquad(1\leq j\leq r,\ y_0\df0)
\]
on $\Q[x_1,\ldots,x_r,y_1,\ldots,y_r]\mspace{1mu}$. One has $\Psi^{0}\circ\iota=\Psi^{0}$, since conjugation is a bijection on the set of partitions of a fixed integer with a given number of different parts.

\begin{lemma}\label{lem:Lprops}\mbox{}
\begin{enumerate}[{\upshape(i)}]\itemsep4pt
\item\label{it:Lprops-id} The map~$\mathcal L^{\epsilon}$ is the identity on $\bigoplus_{r\geq0}(x_1\cdots x_r)^{\epsilon}\Q[x_1,\ldots,x_r]\mspace{1mu}$.
\item\label{it:Lprops-iota} One has $\mathcal L^{0}\circ\iota=\mathcal L^{0}$.
\item\label{it:Lprops-y} For every $p\in\Q[y_1,\dots,y_r]$ one has \[
  \mathcal L^{0}(p)(x_1,\dots,x_r)\=p\bigl(x_r,\,x_{r-1}+x_r,\,\dots,\,x_1+\dots+x_r\bigr)\mspace{1mu}.
\]
\end{enumerate}
\end{lemma}

\begin{proof}
For $c_1,\dots,c_d\geq1$ one has \[
  \xi^{0}_N(1,c_1,\dots,1,c_d;q)\=\zeta^\dagger_{q,N}(c_1,\dots,c_d)\=Z^\dagger_{q,N}(b_{c_1}\cdots b_{c_d})\mspace{1mu},
\]
so that \eqref{eq:HMWfinite} and the injectivity recalled after Theorem~\ref{thm:HMW} give $\EE_{0}(1,c_1,\dots,1,c_d)=b_{c_1}\cdots b_{c_d}\mspace{1mu}$. By \eqref{eq:pcmult} this means that $\mathcal L^{0}$ fixes $p_{(c_1,\dots,c_d)}\mspace{1mu}$. In the same way, the second identity of \cite[Theorem~1.2]{HMW26}, which has no entry equal to~$1$ when all $k'_j\geq2\mspace{1mu}$, and the injectivity in \cite[Lemma~2.6]{HMW26} give $\DD(1,k'_1-1,\dots,1,k'_d-1)=b_{k'_1}\cdots b_{k'_d}\mspace{1mu}$, so that $\mathcal L^{1}$ fixes $\prod_j\binom{x_j}{k'_j-1}\mspace{1mu}$. This proves~\ref{it:Lprops-id}.

The substitution~$\iota$ sends~$x_j$ to $y_{r+1-j}-y_{r-j}$ and $y_j-y_{j-1}$ to $x_{r+1-j}\mspace{1mu}$. Hence it sends the polynomial~\eqref{eq:pc} with index $(0^{\ell_1-1},k_1,\dots,0^{\ell_r-1},k_r)$ to the one with index $(0^{k_r-1},\ell_r,\dots,0^{k_1-1},\ell_1)\mspace{1mu}$, so that $\mathcal L^{0}\circ\iota=\mathcal L^{0}$ by \eqref{eq:reflection}. For $p\in\Q[y_1,\dots,y_r]$ the polynomial~$\iota(p)$ lies in $\Q[x_1,\ldots,x_r]\mspace{1mu}$, and therefore $\mathcal L^{0}(p)=\mathcal L^{0}(\iota(p))=\iota(p)$ by~\ref{it:Lprops-id}.
\end{proof}

By part~\ref{it:Lprops-y}, the restriction of $\mathcal L^{0}$ to the polynomials in the parts alone is the conjugation that gives $\im\Psi^{0}_X=\im\Psi^{0}_Y$ in the introduction. For $\epsilon=1$ the reversal of indices in~\eqref{eq:reflection} gives the counterpart for the Bradley--Zhao model of the duality in Section~\ref{sec:HMWthm}.

\section{Solving the recursion}\label{sec:solving}
In this section we solve the recursion of Definition~\ref{def:Eq}. First we write the recursion for the generating series of the words in Proposition~\ref{prop:rec}, and then we prove Theorem~\ref{thm:intro}.

\subsection{The recursion for the generating series}\label{sec:closed}
Let $S=\{s_1<\dots<s_n\}\subset\Zgeq$ be finite of even cardinality. For $\bfc\in\NN^{S}$ write $\bft^{\bfc}\df\prod_{s\in S}t_s^{c_s}$ and $\EE_{\epsilon}(\bfc)\df\EE_{\epsilon}(c_{s_1},\dots,c_{s_n})\mspace{1mu}$, so that $\gs_S=\sum_{\bfc\in\Zgeq^{S}}\EE_{\epsilon}(\bfc)\,\bft^{\bfc}$, and compute~$\mathrm{eo}\mspace{1mu}$, $\alpha$ and~$\beta$ of a subset of~$S$ from the positions of its elements in~$S\mspace{1mu}$. Put $\Pi_S\df\prod_{s\in S}(1-x\mspace{1mu} t_s)\mspace{1mu}$. We require the following identities.

\begin{lemma}\label{lem:alphabeta}
Let $S\subset\Zgeq$ be finite of even cardinality. Then \[
  \sum_{B\subseteq S}(-1)^{|B|}x^{\alpha(B)}\,\bft^{\boldsymbol e(B)}\=\Pi_S\,\Pp^{*}_S\,,
  \qquad
  \sum_{B\subseteq S}(-1)^{|B|}x^{\beta(B)}\,\bft^{\boldsymbol e(B)}\=\Pi_S\,\Pp_S\,.
\]
\end{lemma}

\begin{proof}
By definition~$\alpha(B)$ is the number of dominoes~$D_i$ with $i$ even meeting the image of~$B\mspace{1mu}$. These dominoes and the two positions that none of them covers, namely $1$ and~$n\mspace{1mu}$, partition~$[n]\mspace{1mu}$, so that the sum factors accordingly: a domino $D_i=\{i,i+1\}$ contributes $1-x(t_{s_i}+t_{s_{i+1}})+x\mspace{1mu} t_{s_i}t_{s_{i+1}}$ according as $B$ contains neither, one or both of its two elements, and an uncovered position~$i$ contributes $1-x\mspace{1mu} t_{s_i}\mspace{1mu}$. Now \[
  \frac{1-x(u+v)+x\mspace{1mu} uv}{(1-xu)(1-xv)}
  \= 1\meno\frac{u}{1-xu}\,\frac{v}{1-xv}\,x(x-1)\mspace{1mu},
\]
 so dividing this product by $\Pi_S$ yields~$\Pp^{*}_S\mspace{1mu}$. The second identity is proven analogously.
\end{proof}

\begin{lemma}\label{lem:Cring}
For $A\in\calT(S)$ the quotient $w^{\epsilon}_{A,S}/(x-1)$ lies in $\ZZ[x]\llbracket t_s\mid s\in S\rrbracket$ and has no constant term. In particular $1+y\,\mathsf C^{\epsilon}_r$ is invertible over $\fh\llbracket t_1,\dots,t_{2r}\rrbracket\mspace{1mu}$.
\end{lemma}

\begin{proof}
For $A\neq\emptyset$ the factor $x^{(\epsilon-1)\kappa(A)}$ combines with $\eta_A=(x(x-1))^{\kappa(A)}\prod_{s\in A}\hat t_s$ into \[x^{\epsilon\kappa(A)}(x-1)^{\kappa(A)}\prod_{s\in A}\hat t_s\,,\] which is divisible by $(x-1)^{\kappa(A)}$ and carries no negative power of~$x$ for $\epsilon\in \{0,1\}\mspace{1mu}$. For $A=\emptyset$ the weight $w_{\emptyset,S}=(\Pp_S-\Pp^{*}_S)/\Pp^{*}_S$ is divisible by $x-1\mspace{1mu}$, as both $\Pp_S$ and $\Pp^{*}_S$ are $\equiv 1$ modulo~${x-1}\mspace{1mu}$. In either case the quotient by $x-1$ has no constant term, since $\hat t_s$ is divisible by~$t_s\mspace{1mu}$.
\end{proof}

\begin{proposition}\label{prop:rec}
For $\epsilon\in\{0,1\}$ the identity~\eqref{eq:recw} holds for every nonempty finite $S\subset\Zgeq$ of even cardinality.
\end{proposition}

\begin{proof}
Multiply the recursion in Definition~\ref{def:Eq} by $\bft^{\bfc}$ and sum over $\bfc\in\Zgeq^{S}\mspace{1mu}$. For $A\in\calT(S)\mspace{1mu}$, including $A=\emptyset\mspace{1mu}$, and $B\subseteq S\setminus A\mspace{1mu}$, the map $\bfc\mapsto\bfc^{A}_{B}$ is a bijection from the $\bfc$ with $A\subseteq\{\bfc=1\}$ and $B\subseteq\{\bfc>1\}$ onto $\Zgeq^{S\setminus A}\mspace{1mu}$, with $\bft^{\bfc}=\bft^{\boldsymbol e(A)}\,\bft^{\boldsymbol e(B)}\,\bft^{\bfc^{A}_{B}}$; all other pairs $(A,B)$ contribute nothing. The sums over~$h$ are \[
  \sum_{h=\alpha(B)+1}^{\beta(B)}\!\!x^{h-1}\meno\!\!\sum_{h=\beta(B)+1}^{\alpha(B)}\!\!x^{h-1}\=\frac{x^{\beta(B)}-x^{\alpha(B)}}{x-1}\,,\qquad
  \sum_{h=1}^{\kappa(A)}(-1)^{\kappa(A)-h}\binom{\kappa(A)-1}{h-1}x^{h-1}\=(x-1)^{\kappa(A)-1}\,.
\]
The sums over~$B$ follow from Lemma~\ref{lem:alphabeta}, which we apply to~$S$ in the first two sums and to $S\setminus A$ in the third. The part with $x^{\beta(B)}$ of the second sum is $\gs_S\,y\,\Pi_S\Pp_S/(x-1)\mspace{1mu}$, which is the term $A=\emptyset$ of the third sum extended to all $A\in\calT(S)\mspace{1mu}$. Hence \[
  \gs_S\,\Pi_S\Pp^{*}_S\=-\gs_S\,y\,\frac{\Pi_S\Pp^{*}_S}{x-1}\+\sum_{A\in\calT(S)}(-1)^{\mathrm{eo}(A)}\,\bft^{\boldsymbol e(A)}\,\gs_{S\setminus A}\;y\,x^{\epsilon\kappa(A)}(x-1)^{\kappa(A)-1}\,\Pi_{S\setminus A}\Pp_{S\setminus A}\,.
\]
Multiply on the right by $(x-1)(\Pi_S\Pp^{*}_S)^{-1}$, which commutes with $x$ and the~$t_s\mspace{1mu}$. As $\Pi_S=\Pi_A\Pi_{S\setminus A}\mspace{1mu}$, one has $\bft^{\boldsymbol e(A)}\Pi_{S\setminus A}\Pi_S^{-1}=\prod_{s\in A}\hat t_s\mspace{1mu}$, and $x^{\epsilon\kappa(A)}(x-1)^{\kappa(A)}\prod_{s\in A}\hat t_s=x^{(\epsilon-1)\kappa(A)}\eta_A\mspace{1mu}$. This gives \[
  \gs_S\,(x-1)\=-\gs_S\,y\+\sum_{A\in\calT(S)}\gs_{S\setminus A}\;y\;(-1)^{\mathrm{eo}(A)}\,x^{(\epsilon-1)\kappa(A)}\,\frac{\Pp_{S\setminus A}}{\Pp^{*}_S}\,\eta_A\,,
\]
that is, $\gs_S\,(x-1)=\sum_{A\in\calT(S)}\gs_{S\setminus A}\,y\,w_{A,S}\mspace{1mu}$. Observe that the term $-\gs_S\,y$ comes from the $-\delta_{A,\emptyset}$ in~$w_{\emptyset,S}\mspace{1mu}$. As $\fh\llbracket t_s\mid s\in S\rrbracket$ has no zero divisors, dividing on the right by $x-1$ (cf.\ Lemma~\ref{lem:Cring}) gives~\eqref{eq:recw}.
\end{proof}

\subsection{The recursion solved}\label{sec:solved}
\begin{proof}[Proof of Theorem~\ref{thm:intro}]
By Proposition~\ref{prop:rec} we have $\gs\,(1+y\,\mathsf C_r)=e_\emptyset\mspace{1mu}$. Observe here that for $S\in\calU(r)$ and $A\in\calT(S)$ one has $S\setminus A\in\calU(r)\mspace{1mu}$, since the two elements of a domino of~$S$ enclose an interval of $[2r]\setminus S$ of even length. As $1+y\,\mathsf C_r$ is invertible by Lemma~\ref{lem:Cring}, this gives
\begin{equation}\label{eq:solved}
  \gs_U\=\bigl[(1+y\,\mathsf C_r)^{-1}\bigr]_{\emptyset,U}\qquad(U\in\calU(r))\mspace{1mu}.
\end{equation}
The case $U=[2r]$ is the theorem.
\end{proof}

\begin{remark}\label{rem:otherentries}
Only the row $T=\emptyset$ of $(1+y\,\mathsf C_r)^{-1}$ enters Theorem~\ref{thm:intro}. Expanding the inverse as a geometric series gives, for all $T,U\in\calU(r)\mspace{1mu}$, \[
  \bigl[(1+y\,\mathsf C_r)^{-1}\bigr]_{T,U}\=\sum_{s\geq0}\ \sum_{T=T_0\subseteq T_1\subseteq\dots\subseteq T_s=U}\bigl(-y\,\mathsf C_{T_0,T_1}\bigr)\cdots\bigl(-y\,\mathsf C_{T_{s-1},T_s}\bigr)\mspace{1mu},
\]
where the inner sum runs over the sequences in~$\calU(r)$ with $T_i\setminus T_{i-1}\in\calT(T_i)$ for all~$i$ (\emph{chains}). In this expression, each coefficient of a monomial in the~$t_s$ is a finite sum, since the entries of $\mathsf C_r$ have no constant term. The expression vanishes unless $T\subseteq U\mspace{1mu}$.
\end{remark}

\begin{proof}[Proof of Corollary~\ref{cor:introcl}]
Grade $\fh\llbracket t_1,\dots,t_{2r}\rrbracket$ by $\deg x=\deg y=1$ and $\deg t_s=-1\mspace{1mu}$. An element has degree at most~$d$ if its coefficient of $\bft^{\bfc}$ has weight at most $d+|\bfc|$ for every~$\bfc\mspace{1mu}$. Its component of degree~$d$ is then well defined, and components of maximal degree multiply. Now $\hat t_s=\sum_{n\geq1}x^{n-1}t_s^{n}$ is homogeneous of degree~$-1\mspace{1mu}$, and $\eta_D=\eta^{\mathrm{cl}}_D-x\hat t_u\hat t_v$ has degree~$0$ with component~$\eta^{\mathrm{cl}}_D\mspace{1mu}$. Hence $w_{A,S}$ has degree at most $(\epsilon-1)\kappa(A)$ with component $w^{\mathrm{cl}}_{A,S}\mspace{1mu}$. By Lemma~\ref{lem:Cring} we have $w_{A,S}=(1-x)\,\mathsf C_{S\setminus A,S}$ with $\mathsf C_{S\setminus A,S}\in\ZZ[x]\llbracket t_s\mid s\in S\rrbracket\mspace{1mu}$, so $\mathsf C_{S\setminus A,S}$ has degree at most $(\epsilon-1)\kappa(A)-1$ with component $-w^{\mathrm{cl}}_{A,S}/x=\mathsf C^{\mathrm{cl}}_{S\setminus A,S}\mspace{1mu}$. Along a sequence $\emptyset=T_0\subseteq\dots\subseteq T_s=[2r]$ as in Remark~\ref{rem:otherentries} the numbers $\kappa(T_i\setminus T_{i-1})$ add up to~$r$ and $s$ letters~$y$ are added, so every term of the expansion of $\gs_r$ there has degree at most $(\epsilon-1)r\mspace{1mu}$, and the component of degree $(\epsilon-1)r$ of~$\gs_r$ is the same expansion with $\mathsf C^{\mathrm{cl}}_r$ in place of~$\mathsf C_r\mspace{1mu}$, which is the right-hand side. A word~$w$ in the coefficient of $\bft^{\bfc}$ has degree $\wt(w)-|\bfc|\mspace{1mu}$, so this component is $\sum_{\bfc}\EE^{\mathrm{cl}}_{\epsilon}(\bfc)\,\bft^{\bfc}$.
\end{proof}

Note that the proof gives the weight bound~\eqref{eq:weightbound} once more.

\section{From the weights to the words}\label{sec:formal}
In this section we first collect properties of the weights and then expand the inverse in Theorem~\ref{thm:intro}, to obtain a formula for the individual words in Corollary~\ref{cor:main}. In the subsequent Section~\ref{sec:qseries}, we interpret the results in this section on the level of $q$-series by applying the evaluation maps~$Z^\dagger_{q}\mspace{1mu}$.

\subsection{Expanding the inverse}\label{sec:weights}
In depth two, $\calU(2)$ consists of the sets~$\emptyset\mspace{1mu}$, $\{1,2\}\mspace{1mu}$, $\{3,4\}\mspace{1mu}$, $\{1,4\}$ and~$[4]\mspace{1mu}$. For $\epsilon=1$ the weights are $w_{S,S}=\eta_S$ and $w_{\emptyset,S}=-\eta_S$ for $S\in\{\{1,2\},\{3,4\},\{1,4\}\}\mspace{1mu}$, and
\begin{align*}
  w_{\{1,2\},[4]}&\=\frac{\eta_{\{1,2\}}\,(1-\eta_{\{3,4\}})}{1-\eta_{\{2,3\}}}\,,&
  w_{\{3,4\},[4]}&\=\frac{\eta_{\{3,4\}}\,(1-\eta_{\{1,2\}})}{1-\eta_{\{2,3\}}}\,,\\[4pt]
  w_{\{2,3\},[4]}&\=-\frac{\eta_{\{2,3\}}\,(1-\eta_{\{1,4\}})}{1-\eta_{\{2,3\}}}\,,&
  w_{\emptyset,[4]}&\=\frac{(1-\eta_{\{1,2\}})(1-\eta_{\{3,4\}})}{1-\eta_{\{2,3\}}}-1\,.
\end{align*}
For $\epsilon=0$ the weights $w_{A,S}$ with $A\neq\emptyset$ are to be multiplied by~$x^{-1}$. Note that the weight $w_{[4],[4]}$ is not defined, because $[4]=\{1,2\}\cup\{3,4\}\mspace{1mu}$, two adjacent odd--even dominoes, is not a partial domino tiling. For $\epsilon=1$ the weights $w_{A,[4]}$ sum to~$0\mspace{1mu}$: multiplied by $1-\eta_{\{2,3\}}\mspace{1mu}$, this says \[(1-\eta_{\{1,2\}})(1-\eta_{\{3,4\}})+\eta_{\{3,4\}}(1-\eta_{\{1,2\}})+\eta_{\{1,2\}}(1-\eta_{\{3,4\}})-\eta_{\{2,3\}}(1-\eta_{\{1,4\}})=1-\eta_{\{2,3\}}\,.\] 
(This follows from $\eta_{\{1,4\}}\eta_{\{2,3\}}=\eta_{\{1,2\}}\eta_{\{3,4\}}\mspace{1mu}$.)
For $\epsilon=0$ the weights entering~\eqref{eq:recw} for a given~$S$ do not sum to~$0$ (already in depth one $w^{0}_{\emptyset,[2]}+w^{0}_{[2],[2]}=(x^{-1}-1)\,\eta_{[2]}$). For $\epsilon=1$ they always do. For $S\in\calU(r)$ this means that the column~$S$ of~$\mathsf C^{1}_r$ sums to~$0\mspace{1mu}$.

\begin{proposition}\label{prop:collapse}
For every finite $S\subset\Zgeq$ of even cardinality, $\sum_{A\in\calT(S)}w^{1}_{A,S}=0\mspace{1mu}$.
\end{proposition}

\begin{proof} 
Multiplying by $\Pp^{*}_S$ on the right, the claim becomes \[
  \sum_{A\in\calT(S)}(-1)^{\mathrm{eo}(A)}\,\eta_A\,\Pp_{S\setminus A}\=\Pp^{*}_S\,.
\]
Without loss of generality, assume $S=[n]\mspace{1mu}$. Then $s_i=i$ and $D_i(S)=D_i\mspace{1mu}$. For odd $i\leq n-1$ put $[i,n]\df\{i,i+1,\dots,n\}$ and \[
  F(i)\df\sum_{A\in\calT([i,n])}(-1)^{\mathrm{eo}(A)}\,\eta_A\,\Pp_{[i,n]\setminus A}\,,\qquad
  F^{*}(i)\df\sum_{\substack{A\in\calT([i,n])\\ D_i\not\subseteq A}}(-1)^{\mathrm{eo}(A)}\,\eta_A\,\Pp_{[i,n]\setminus A}\,,
\]
so that $F(1)$ is the left-hand side above. Put $F(n+1)=F^{*}(n+1)=1\mspace{1mu}$, and let $\Pp^{*}(i)$ be the product of the factors $1-\eta_{D_j}$ over the even~$j$ with $i<j<n\mspace{1mu}$, so that $\Pp^{*}(1)=\Pp^{*}_S$ and $\Pp^{*}(n-1)=\Pp^{*}(n+1)=1\mspace{1mu}$. We show $F(i)=\Pp^{*}(i)$ by descending induction on the odd number~$i\mspace{1mu}$. The case $i=n+1$ is the convention $F(n+1)=1=\Pp^{*}(n+1)\mspace{1mu}$.

Splitting off the domino~$D_i\mspace{1mu}$, of weight $+\eta_{D_i}\mspace{1mu}$, which excludes an odd--even domino at $i+2\mspace{1mu}$, gives $F^{*}(i)=F(i)-\eta_{D_i}F^{*}(i+2)$ and hence, on iterating,
\begin{equation}\label{eq:Fstar}
  F^{*}(i)\=\sum_{k\geq0}(-1)^{k}\,\eta_{D_i}\eta_{D_{i+2}}\cdots\eta_{D_{i+2k-2}}\,F(i+2k)\mspace{1mu}.
\end{equation}
In $F(i)$ either~$D_i$ belongs to the tiling and contributes $\eta_{D_i}F^{*}(i+2)\mspace{1mu}$, or $i$ is left uncovered. In the second case $\Pp$ pairs~$i$ with the first uncovered position to its right, say $i+1+2m$ with $0\leq m\leq\frac12(n-i-1)\mspace{1mu}$, and every position in between is covered. Such a partner exists, since the tail $i,\dots,n$ has even length and the covered positions come in dominoes, so that the number of uncovered positions is even and hence at least two. An interval of even length has exactly one tiling by dominoes, here $D_{i+1},D_{i+3},\dots,D_{i+2m-1}\mspace{1mu}$. These $m$ dominoes are all even--odd, since their indices have the parity of $i+1\mspace{1mu}$. As even--odd dominoes they may be adjacent to one another, and they leave the tail beyond $i+1+2m$ unrestricted. Hence
\begin{equation}\label{eq:split}
  F(i)\=\eta_{D_i}F^{*}(i+2)\+\sum_{m\geq0}(-1)^{m}
  \eta_{D_{i+1}}\eta_{D_{i+3}}\cdots\eta_{D_{i+2m-1}}
  \bigl(1-\eta_{\{i,i+1+2m\}}\bigr)F(i+2m+2)\mspace{1mu}.
\end{equation}
The two terms of \eqref{eq:split} carrying~$t_i$ cancel. Indeed, by \eqref{eq:Fstar} the first equals \[
  \sum_{k\geq0}(-1)^{k}\,\eta_{D_i}\eta_{D_{i+2}}\cdots\eta_{D_{i+2k}}\,F(i+2k+2)\mspace{1mu},
\]
and the contribution of $-\eta_{\{i,i+1+2m\}}$ to the second is the same sum with the opposite sign, since $\eta_{D_{i+1}}\cdots\eta_{D_{i+2m-1}}\,\eta_{\{i,i+1+2m\}} =\eta_{D_i}\cdots\eta_{D_{i+2m}}\mspace{1mu}$, as both products equal $\bigl(x(x-1)\bigr)^{m+1}\hat t_i\cdots\hat t_{i+2m+1}\mspace{1mu}$. What is left is \[
  F(i)\=\sum_{m\geq0}(-1)^{m}\,\eta_{D_{i+1}}\eta_{D_{i+3}}\cdots\eta_{D_{i+2m-1}}
  \,F(i+2m+2)\mspace{1mu}.
\]
By the inductive hypothesis $F(j)=\Pp^{*}(j)$ for odd $j>i\mspace{1mu}$, and $\Pp^{*}(j)=(1-\eta_{D_{j+1}})\Pp^{*}(j+2)$ for odd $j\leq n-3\mspace{1mu}$. Hence the $m$th summand equals \[
  (-1)^{m}\,\eta_{D_{i+1}}\cdots\eta_{D_{i+2m-1}}\,\Pp^{*}(i+2m)\meno(-1)^{m+1}\,\eta_{D_{i+1}}\cdots\eta_{D_{i+2m+1}}\,\Pp^{*}(i+2m+2)\mspace{1mu},
\]
except for the last one, $m=\frac12(n-i-1)\mspace{1mu}$, which equals the first term alone, as $\Pp^{*}(n-1)=\Pp^{*}(n+1)=1\mspace{1mu}$. The sum telescopes, and $F(i)=\Pp^{*}(i)\mspace{1mu}$.
\end{proof}

\begin{corollary}\label{cor:signcount}
For every $r\geq1\mspace{1mu}$, \[
\sum_{A\in\calT(r)}(-1)^{\mathrm{eo}(A)}X^{\kappa(A)}\=1+X\,.\]
\end{corollary}

\begin{proof}
Give all the $\hat t_s$ the same value~$\hat t\mspace{1mu}$, so that every $\eta_{D}$ becomes $\eta\df x(x-1)\hat t^{\,2}$. Proposition~\ref{prop:collapse} for $S=[2r]\mspace{1mu}$, in the form of the first display of its proof, then reads \[
  \sum_{A}(-1)^{\mathrm{eo}(A)}\eta^{\kappa(A)}(1-\eta)^{r-\kappa(A)}\=(1-\eta)^{r-1}\,.
\]
The claimed identity follows by dividing by $(1-\eta)^{r}$ and putting $X=\eta/(1-\eta)\mspace{1mu}$.
\end{proof}

For $\epsilon=1\mspace{1mu}$, Proposition~\ref{prop:collapse} says that every column of $\mathsf C^{1}_r$ sums to~$0\mspace{1mu}$, so that every column of $1+y\,\mathsf C^{1}_r\mspace{1mu}$, and hence of its inverse, sums to~$1\mspace{1mu}$. By Remark~\ref{rem:otherentries} this means that for every $U\in\calU(r)$ \[
  \sum_{T\in\calU(r)}\ \sum_{s\geq0}\ \sum_{T=T_0\subseteq T_1\subseteq\dots\subseteq T_s=U}\bigl(-y\,\mathsf C^{1}_{T_0,T_1}\bigr)\cdots\bigl(-y\,\mathsf C^{1}_{T_{s-1},T_s}\bigr)\=1\,,
\]
where the inner sum runs over the sequences in~$\calU(r)$ with $T_i\setminus T_{i-1}\in\calT(T_i)$ for all~$i\mspace{1mu}$. 

\begin{definition}For $T,U\in\calU(r)\mspace{1mu}$, a \emph{chain} from $T$ to~$U$ is a sequence $T=T_0\subseteq T_1\subseteq\dots\subseteq T_s=U$ in~$\calU(r)$ with $T_i\setminus T_{i-1}\in\calT(T_i)$ for all $i\in[s]\mspace{1mu}$. It is \emph{strict} if all inclusions are strict.
\end{definition}

 By Theorem~\ref{thm:intro} and Remark~\ref{rem:otherentries} with $T=\emptyset$ and $U=[2r]$ we obtain
\begin{equation}\label{eq:paths}
  \gs_{[2r]}\=\sum_{s\geq0}(-1)^{s}\!\!\!
  \sum_{\emptyset=T_0\subseteq T_1\subseteq\dots\subseteq T_s=[2r]}\!\!\!
  y\,\mathsf C_{T_0,T_1}\;y\,\mathsf C_{T_1,T_2}\cdots y\,\mathsf C_{T_{s-1},T_s}\,,
\end{equation}
where the inner sum runs over the chains from $\emptyset$ to~$[2r]\mspace{1mu}$. The $s$th summand carries exactly~$s$ letters~$y\mspace{1mu}$, so it is the part of $\gs_{[2r]}$ of depth~$s\mspace{1mu}$. Hence the part of depth~$s$ of $\EE_{\epsilon}(\bfc)$ is the coefficient of $\bft^{\bfc}$ in $\bigl[(-y\,\mathsf C_r)^{s}\bigr]_{\emptyset,[2r]}\mspace{1mu}$. In Section~\ref{sec:words}, we expand \eqref{eq:paths} in the~$t_s\mspace{1mu}$, once in Section~\ref{sec:coeffs} we have identified the coefficients of the entries of~$\mathsf C_r\mspace{1mu}$.

 Resumming the steps with $T_{i-1}=T_i$ turns \eqref{eq:paths} into a sum over strict chains,
\begin{equation}\label{eq:chains}
  \gs_{[2r]}\=\sum_{s\geq 0}\sum_{\emptyset=U_0\subsetneq U_1\subsetneq\dots\subsetneq
  U_s=[2r]}\ \prod_{i=1}^{s}\bigl(-y\,\mathsf C_{U_{i-1},U_i}\bigr)
  \bigl(1+y\,\mathsf C_{U_i,U_i}\bigr)^{-1}\,.
\end{equation}

In depth one \eqref{eq:chains} is \eqref{eq:depthone}. In depth two the strict chains from $\emptyset$ to $[4]$ are $\emptyset\subsetneq U_1\subsetneq[4]$ with $U_1\in\{\{1,2\},\{3,4\},\{1,4\}\}\mspace{1mu}$, and \eqref{eq:chains} becomes
\begin{align*}
  \gs_{[4]}\=\Bigl[\;
  &y\,x^{\epsilon}\hat t_1\hat t_2\,\bigl(1+y\,x\,\hat t_1\hat t_2\bigr)^{-1}\,y\,x^{\epsilon}\hat t_3\hat t_4\,\bigl(1-x(x-1)\hat t_1\hat t_2\bigr)\\
  &\+ y\,x^{\epsilon}\hat t_3\hat t_4\,\bigl(1+y\,x\,\hat t_3\hat t_4\bigr)^{-1}\,y\,x^{\epsilon}\hat t_1\hat t_2\,\bigl(1-x(x-1)\hat t_3\hat t_4\bigr)\\
  &\meno y\,x^{\epsilon}\hat t_1\hat t_4\,\bigl(1+y\,x\,\hat t_1\hat t_4\bigr)^{-1}\,y\,x^{\epsilon}\hat t_2\hat t_3\,\bigl(1-x(x-1)\hat t_1\hat t_4\bigr)
  \;\Bigr]\\
  &\qquad\Bigl(1-x(x-1)\hat t_2\hat t_3\+y\,x\,\bigl(\hat t_1\hat t_2+\hat t_3\hat t_4-\hat t_2\hat t_3-x(x-1)\hat t_1\hat t_2\hat t_3\hat t_4\bigr)\Bigr)^{-1}.
\end{align*}

\subsection{\texorpdfstring{The coefficients of $\mathsf C_r$}{The coefficients of C}}\label{sec:coeffs}
To read the identity in Theorem~\ref{thm:intro} coefficientwise we need the coefficients of the inverse of $\Pi_S\,\Pp^{*}_S\mspace{1mu}$. They are products of the following polynomials, which invert the two-variable building block $1-x(u+v)+x\mspace{1mu} uv$ of Lemma~\ref{lem:alphabeta}. For $k,\ell\in\ZZ$ let
\begin{equation}\label{eq:gammagen}
  \sum_{k,\ell\geq0}\gamma(k,\ell;x)\,u^{k}v^{\ell}\=\bigl(1-x(u+v)+x\mspace{1mu} uv\bigr)^{-1}\,,
\end{equation}
that is, expanding the geometric series and collecting the terms with $m$ factors $-uv\,x\mspace{1mu}$, \[
\gamma(k,\ell;x) \= \sum_{m\geq 0} (-1)^m \binom{k+\ell-m}{k-m,\ell-m,m} x^{k+\ell-m}\,,
\]
where the multinomial coefficients are $0$ if $k<m$ or $\ell<m\mspace{1mu}$. In particular $\gamma(k,\ell;x)=0$ as soon as $k<0$ or $\ell<0\mspace{1mu}$, a convention used throughout below. Moreover $\gamma(k,0;x)=x^k=\gamma(0,k;x)$ for ${k\geq0}\mspace{1mu}$, so that $\gamma(0,0;x)=1\mspace{1mu}$. Now define polynomials parametrized by a vector ${\bfa=(a_1,\ldots,a_n)\in \ZZ^n}$ with $n$ even by
\begin{align*}
\gamma_{\bfa}(x)\defis \prod_{\substack{0\leq i \leq n \\ i \text{ even}}}\gamma(a_i,a_{i+1};x) \qquad (a_i=0 \text{ if } i\not \in [n])\mspace{1mu}.
\end{align*}
The product is taken over the blocks $\{1\},\{2,3\},\dots,\{n-2,n-1\},\{n\}$ into which the dominoes~$D_i$ with $i$ even divide~$[n]\mspace{1mu}$. These blocks govern~$\Pi_S\Pp^{*}_S$ in Lemma~\ref{lem:alphabeta}. For example, for $a_1,a_2,a_3,a_4\geq 0$ we have \[
\gamma_{(a_1,a_2)}(x)\=x^{a_1+a_2}, \qquad \gamma_{(a_1,a_2,a_3,a_4)}(x) \= \sum_{m\geq 0}(-1)^m \binom{a_2+a_3-m}{a_2-m,a_3-m,m} x^{a_1+a_2+a_3+a_4-m}.
\]

\begin{proposition}\label{prop:sumsofgamma} Let $n>0$ be even. Then 
\begin{enumerate}[{\upshape(i)}]\itemsep4pt
\item\label{it:gamma-inverse} $\sum_{\bfa\in\NN^{n}}\gamma_{\bfa}(x)\,\bft^{\bfa}$ is the inverse of $\Pi_{[n]}\,\Pp^{*}_{[n]}\mspace{1mu}$, equivalently, for every $\bfa \in \ZZ^n$, \[\sum_{\substack{B\subseteq[n]}} (-1)^{|B|} x^{\alpha(B)} \gamma_{\bfa-\boldsymbol{e}(B)}(x) \= \begin{cases} 1 &\text{if } \bfa=\boldsymbol{0}\mspace{1mu}, \\ 0 & \text{else.}\end{cases}  \]
\item\label{it:gamma-mod} $\Pp_{[n]}/\Pp^{*}_{[n]}\equiv1$ modulo $x-1\mspace{1mu}$, equivalently, for every $\bfa\in\ZZ^n$, \[\sum_{\substack{B\subseteq[n]}} (-1)^{|B|} x^{\beta(B)} \gamma_{\bfa-\boldsymbol{e}(B)}(x) \  \in \ \begin{cases} \{1\} & \text{if } \bfa=\boldsymbol{0}\mspace{1mu}, \\  (x-1)\,\ZZ[x] & \text{else.}\end{cases}  \]
\end{enumerate}
\end{proposition}
\begin{proof}
By the proof of the first identity of Lemma~\ref{lem:alphabeta}, we see that $\Pi_{[n]}\Pp^{*}_{[n]}$ is the product of $1-x\mspace{1mu} t_1\mspace{1mu}$, of $1-x(t_i+t_{i+1})+x\mspace{1mu} t_it_{i+1}$ for the dominoes $\{i,i+1\}$ with $i$ even, and of $1-x\mspace{1mu}t_n\mspace{1mu}$. The series $\sum_{\bfa}\gamma_{\bfa}\bft^{\bfa}$ factors similarly, with factors $\sum_{k}x^{k}t_1^{k}=(1-x\mspace{1mu} t_1)^{-1}$, $\sum_{k,\ell}\gamma(k,\ell;x)t_i^{k}t_{i+1}^{\ell}$ and $(1-x\mspace{1mu}t_n)^{-1}$, so~\ref{it:gamma-inverse} is \eqref{eq:gammagen}. Comparing the coefficient of $\bft^{\bfa}$ in the product of the two series gives the displayed identity, as the coefficient of $\bft^{\boldsymbol e(B)}$ in $\Pi_{[n]}\Pp^{*}_{[n]}$ is $(-1)^{|B|}x^{\alpha(B)}$ by Lemma~\ref{lem:alphabeta}. 

For~\ref{it:gamma-mod}, every $\eta_D$ is divisible by $x-1\mspace{1mu}$, so that $\Pp_{[n]}/\Pp^{*}_{[n]}$ has constant term~$1$ and is $\equiv1$ modulo $x-1\mspace{1mu}$. Since \[
  \frac{\Pp_{[n]}}{\Pp^{*}_{[n]}}\=\bigl(\Pi_{[n]}\Pp_{[n]}\bigr)\bigl(\Pi_{[n]}\Pp^{*}_{[n]}\bigr)^{-1},
\]
the second identity of Lemma~\ref{lem:alphabeta} and~\ref{it:gamma-inverse} give its coefficient of $\bft^{\bfa}$ as the displayed sum. 
\end{proof}

Given a subset $A\subseteq [n]$ of even cardinality, ${\bfa=(a_1,\ldots,a_n)\in \NN^n}$ and $\epsilon\in \{0,1\}\mspace{1mu}$, we define $p_{A,\bfa}(x)\df0$ if $A=\emptyset$ and $\bfa=\boldsymbol{0}\mspace{1mu}$, or if $a_i=0$ for some $i\in A\mspace{1mu}$, and in all other cases by \[
p^{\epsilon}_{A,\bfa}(x)\defis p_{A,\bfa}(x)\defis \sum_{\substack{B\subseteq [n]\\A\cap B=\emptyset}} (-1)^{\mathrm{eo}(A)+|B|}x^{\epsilon\kappa(A)+\beta(\sigma_{A}(B))}(x-1)^{\kappa(A)-1}\gamma_{\bfa-\boldsymbol{e}(A\cup B)}(x)\mspace{1mu}.
\]
For $A\neq\emptyset$ this is a polynomial. For $A=\emptyset$ the factor $(x-1)^{\kappa(A)-1}$ is to be read as division by $x-1\mspace{1mu}$, which leaves no remainder in $\ZZ[x]$ by Proposition~\ref{prop:sumsofgamma}\ref{it:gamma-mod}. For a finite set $S\subset\Zgeq$ of even cardinality, $A\subseteq S$ and $\bfa\in\NN^{S}$ we define $p_{A,\bfa}$ by the same formula, with $A\mspace{1mu}$, $B$ and~$\bfa$ read through the positions in~$S\mspace{1mu}$, as in Section~\ref{sec:closed}.

These polynomials are the coefficients of the entries of the matrix~$\mathsf C_r$ of Theorem~\ref{thm:intro}, and more generally of $\mathsf C_{S\setminus A,S}\df w_{A,S}/(1-x)$ for arbitrary finite $S\subset\Zgeq$ of even cardinality and $A\in\calT(S)\mspace{1mu}$.

\begin{lemma}\label{lem:pcoeff}
Let $S\subset\Zgeq$ be finite of even cardinality and $A\in\calT(S)\mspace{1mu}$. Then
\begin{enumerate}[{\upshape(i)}]\itemsep4pt
\item\label{it:pcoeff-coeff} for every $\bfa\in\NN^{S}$, the coefficient of $\bft^{\bfa}$ in $-\mathsf C^{\epsilon}_{S\setminus A,S}$ equals $p^{\epsilon}_{A,\bfa}(x)\mspace{1mu}$;
\item\label{it:pcoeff-series} $-\mathsf C^{\epsilon}_{S\setminus A,S}$ is a power series with coefficients in $\ZZ[x,x^{-1}]$ in the products~$\hat t_s\hat t_{s'}\mspace{1mu}$, where $\{s,s'\}$ runs over the pairs of elements that are consecutive in~$S$ or consecutive in~$S\setminus A\mspace{1mu}$.
\end{enumerate}
\end{lemma}

\begin{proof}
Write $T\defis S\setminus A\mspace{1mu}$, $\kappa\defis\kappa(A)\mspace{1mu}$, and recall from Lemma~\ref{lem:alphabeta} and Proposition~\ref{prop:sumsofgamma} that, for every finite $R\subset\Zgeq$ of even cardinality,
\begin{equation}\label{eq:threeseries}
  \Pi_R\,\Pp_R\=\sum_{B\subseteq R}(-1)^{|B|}x^{\beta(B)}\bft^{\boldsymbol e(B)},\qquad
  \bigl(\Pi_R\,\Pp^{*}_R\bigr)^{-1}\=\sum_{\bfb\in\NN^{R}}\gamma_{\bfb}(x)\,\bft^{\bfb},
\end{equation}
where $\beta$ and $\gamma_{\bfb}$ are computed from the positions in~$R\mspace{1mu}$. The second identity is Proposition~\ref{prop:sumsofgamma}\ref{it:gamma-inverse}, which says that $\sum_{\bfb}\gamma_{\bfb}\bft^{\bfb}$ is inverse to $\sum_B(-1)^{|B|}x^{\alpha(B)}\bft^{\boldsymbol e(B)}$.

If $A\neq\emptyset\mspace{1mu}$, then $-\mathsf C_{T,S}=w_{A,S}/(x-1)$ and, since $\eta_A=(x(x-1))^{\kappa}\prod_{s\in A}\hat t_s$ and $\Pi_A\,\Pi_T=\Pi_S\mspace{1mu}$, \[
  \frac{w_{A,S}}{x-1}
  \=(-1)^{\mathrm{eo}(A)}x^{\epsilon\kappa}(x-1)^{\kappa-1}\,\bft^{\boldsymbol e(A)}\,
  \bigl(\Pi_T\,\Pp_T\bigr)\bigl(\Pi_S\,\Pp^{*}_S\bigr)^{-1}\,.
\]
Inserting \eqref{eq:threeseries} for $R=T$ and $R=S$ and reading off the coefficient of $\bft^{\bfa}$ gives \[
  (-1)^{\mathrm{eo}(A)}x^{\epsilon\kappa}(x-1)^{\kappa-1}\sum_{B\subseteq T}(-1)^{|B|}x^{\beta_T(B)}\,
  \gamma_{\bfa-\boldsymbol e(A)-\boldsymbol e(B)}(x)\mspace{1mu},
\]
where $\beta_T(B)\mspace{1mu}$, the function~$\beta$ of $B$ inside $T=S\setminus A\mspace{1mu}$, is $\beta(\sigma_A(B))$ in the notation of Section~\ref{subsec:notation}. This is $p_{A,\bfa}(x)\mspace{1mu}$. If $a_i=0$ for some $i\in A\mspace{1mu}$, the coefficient vanishes because of the factor~$\bft^{\boldsymbol e(A)}$, in accordance with the convention.

If $A=\emptyset\mspace{1mu}$, then $-\mathsf C_{S,S}=w_{\emptyset,S}/(x-1)$ with $w_{\emptyset,S}+1=\Pp_S/\Pp^{*}_S=(\Pi_S\Pp_S)(\Pi_S\Pp^{*}_S)^{-1}$, whose coefficient of $\bft^{\bfa}$ is $\sum_{B\subseteq S}(-1)^{|B|}x^{\beta(B)}\gamma_{\bfa-\boldsymbol e(B)}(x)$ by \eqref{eq:threeseries}. By Proposition~\ref{prop:sumsofgamma}\ref{it:gamma-mod} this is $1$ for $\bfa=\boldsymbol 0\mspace{1mu}$, so that the constant term of $-\mathsf C_{S,S}$ vanishes, and for $\bfa\neq\boldsymbol0$ it is divisible by $x-1$ with quotient $p_{\emptyset,\bfa}(x)\mspace{1mu}$. This proves~\ref{it:pcoeff-coeff}.

For~\ref{it:pcoeff-series}, $\Pp_T$ is a polynomial in the products $\eta_{\{s,s'\}}=x(x-1)\hat t_s\hat t_{s'}$ over pairs consecutive in~$T\mspace{1mu}$, the products~$\Pp^{*}_S$ and $\Pp_S$ are polynomials in such products over pairs consecutive in~$S\mspace{1mu}$, and $\eta_A$ is the product of $x(x-1)\hat t_s\hat t_{s'}$ over the dominoes $\{s,s'\}$ of~$A\mspace{1mu}$, whose two elements are consecutive in~$S\mspace{1mu}$. As $\Pp^{*}_S$ has constant term~$1\mspace{1mu}$, its inverse is a power series in these products, and hence so are $w_{A,S}$ and $-\mathsf C_{T,S}=w_{A,S}/(x-1)\mspace{1mu}$.
\end{proof}

In depth one, by \eqref{eq:depthone}, $-\mathsf C_{\emptyset,[2]}=x^{\epsilon}\hat t_1\hat t_2$ and $-\mathsf C_{[2],[2]}=-x\,\hat t_1\hat t_2\mspace{1mu}$, so that for $\bfa\in\Zgeq^{2}$ \[
  p_{\{1,2\},\bfa}(x)\=x^{a_1+a_2+\epsilon-2}\,,\qquad p_{\emptyset,\bfa}(x)\=-x^{a_1+a_2-1}\,,
\]
while $p_{A,\bfa}=0$ if $a_1=0$ or $a_2=0\mspace{1mu}$. In depth two the polynomials~$\gamma$ enter. For instance, the weights of Section~\ref{sec:weights} give \[
  -\mathsf C_{[4],[4]}\big|_{t_1=t_4=0}\=\frac{x\,t_2t_3}{1-x(t_2+t_3)+x\,t_2t_3}\,,
\]
so that $p_{\emptyset,(0,a_2,a_3,0)}(x)=x\,\gamma(a_2-1,a_3-1;x)$ by \eqref{eq:gammagen}. For example, $p_{\emptyset,(0,2,2,0)}(x)=2x^{3}-x^{2}$.

Part~\ref{it:pcoeff-series} of Lemma~\ref{lem:pcoeff} forces many of the $p_{A,\bfa}$ to vanish.

\begin{corollary}\label{cor:vanishing}
Let $S\subset\Zgeq$ be finite of even cardinality, $A\in\calT(S)$ and $\bfa\in\NN^{S}$. If some $i\in S$ has $a_i>0$ and $a_j=0$ for every neighbour~$j$ of $i$ in~$S\mspace{1mu}$, then $p^{\epsilon}_{A,\bfa}(x)=0\mspace{1mu}$. Moreover $p^{\epsilon}_{\emptyset,\bfa}(x)=0$ whenever~$\bfa$ has at most one nonzero entry.
\end{corollary}

\begin{proof}
If $\bfa$ vanishes at an element of~$A\mspace{1mu}$, then $p_{A,\bfa}=0$ by definition. Otherwise no neighbour of $i$ lies in~$A\mspace{1mu}$, so $i\notin A$ and the neighbours of $i$ in $S\setminus A$ are those in~$S\mspace{1mu}$. By Lemma~\ref{lem:pcoeff}\ref{it:pcoeff-series}, $-\mathsf C_{S\setminus A,S}$ is a power series in the products~$\hat t_s\hat t_{s'}$ with $s'$ a neighbour of $s$ in $S$ or in $S\setminus A\mspace{1mu}$, and $\hat t_s$ is divisible by~$t_s\mspace{1mu}$. Hence every monomial of $-\mathsf C_{S\setminus A,S}$ containing~$t_i$ contains~$t_j$ for a neighbour~$j$ of $i$ in~$S\mspace{1mu}$, so the coefficient of $\bft^{\bfa}$ vanishes. The second claim is the first one together with $p_{\emptyset,\boldsymbol0}=0\mspace{1mu}$.
\end{proof}

\subsection{Words}\label{sec:words}
We now take the coefficient of $\bft^{\bfc}$ in~\eqref{eq:paths}. By Lemma~\ref{lem:pcoeff} the coefficients of the links $-\mathsf C_{T_{j-1},T_j}$ are the polynomials $p_{T_j\setminus T_{j-1},\bfa}\mspace{1mu}$. We record a chain together with a decomposition of~$\bfc$ over its steps in a matrix.

For a matrix $M=(m_{i,j})\in\NN^{2r\times s}$ and $j\in\{0,\dots,s\}$ let \[
  T_j(M)\defis\bigl\{\,i\in[2r]\ \big|\ m_{i,j'}>0\text{ for some }j'\leq j\,\bigr\}
\]
be the set of rows having a nonzero entry among the first~$j$ columns, so that $\emptyset=T_0(M)\subseteq T_1(M)\subseteq\dots\subseteq T_s(M)\mspace{1mu}$, and put $A_j(M)\defis T_j(M)\setminus T_{j-1}(M)\mspace{1mu}$, the set of rows whose first nonzero entry stands in column~$j\mspace{1mu}$. Write $\bfa_j(M)\in\NN^{T_j(M)}$ for the $j$th column of~$M$ with the entries in the rows outside~$T_j(M)\mspace{1mu}$, which are zero, deleted.

\begin{definition}\label{def:dominomatrix}
For $\bfc\in\Zgeq^{2r}$ and $s\geq0\mspace{1mu}$, let $\mathcal{M}_{2r\times s}(\bfc)$ be the set of matrices $M\in\NN^{2r\times s}$ without zero columns, whose $i$th row sums to $c_i$ for every $i\in[2r]\mspace{1mu}$, and for which $T_0(M)\subseteq\dots\subseteq T_s(M)$ is a chain.
\end{definition}

As $\bfc$ has no zero entry, $T_s(M)=[2r]$ holds automatically. The set $\mathcal{M}_{2r\times s}(\bfc)$ is finite and empty for $s>|\bfc|\mspace{1mu}$, since each column has sum at least~$1\mspace{1mu}$.

\begin{remark}\label{rem:polrecursion}
By Lemma~\ref{lem:pcoeff}, comparing the coefficients of $\bft^{\bfc}$ in~\eqref{eq:recw} for $S=[2r]$ gives the recursion \[
\EE_{\epsilon}(\bfc)
        \=
  \sum_{A\in\calT(r)}\ \sum_{\substack{\bfd\in\Zgeq^{2r}\\ \bfd\leq\bfc}}
      \,\EE_{\epsilon}(\bfd^A)\,
      y\,p^{\epsilon}_{A,\bfc-\bfd+\boldsymbol{e}(A)}(x)\mspace{1mu},
\]
in which a term vanishes unless~$\bfd$ equals~$1$ on~$A\mspace{1mu}$, and the term with $A=\emptyset$ and $\bfd=\bfc$ is absent, since $p_{\emptyset,\boldsymbol0}=0\mspace{1mu}$. Iterating it gives a slightly different proof of Corollary~\ref{cor:main} below.
\end{remark}

\begin{corollary}\label{cor:main} For $r\geq 0\mspace{1mu}$, $\bfc\in\Zgeq^{2r}$ and $\epsilon\in\{0,1\}\mspace{1mu}$,
\begin{align*}
\EE_{\epsilon}(\bfc)
        \= &
  \sum_{s\geq 0}\ \sum_{\substack{M\in \mathcal{M}_{2r\times s}(\bfc)}} y\,p^{\epsilon}_{A_1(M),\bfa_1(M)}(x) \cdots
      y\,p^{\epsilon}_{A_s(M),\bfa_s(M)}(x)\mspace{1mu}.
\end{align*}
\end{corollary}

\begin{proof}
The entry $\mathsf C_{T,S}$ involves only the variables~$t_s$ with $s\in S\mspace{1mu}$, and these commute with everything except~$y\mspace{1mu}$.  Hence the coefficient of $\bft^{\bfc}$ in the summand of \eqref{eq:paths} belonging to the chain $\emptyset=T_0\subseteq\dots\subseteq T_s=[2r]$ is \[
  \sum_{\bfa_1+\dots+\bfa_s=\bfc}\ \prod_{j=1}^{s}\,y\,[\bft^{\bfa_j}]\bigl(-\mathsf C_{T_{j-1},T_j}\bigr)
  \=\sum_{\bfa_1+\dots+\bfa_s=\bfc}\ \prod_{j=1}^{s}\,y\,p_{T_j\setminus T_{j-1},\bfa_j}(x)\mspace{1mu},
\]
by Lemma~\ref{lem:pcoeff}, where the sign~$(-1)^s$ is absorbed into the $s$ factors, the sum runs over the decompositions of $\bfc$ into vectors $\bfa_j\in\NN^{2r}$ supported on~$T_j\mspace{1mu}$, each $\bfa_j$ is identified with its restriction to $T_j\mspace{1mu}$, on which $p^{\epsilon}$ is evaluated, and the products are written from left to right. A term vanishes unless every $\bfa_j$ is positive on $T_j\setminus T_{j-1}\mspace{1mu}$, by the convention on $p_{A,\bfa}\mspace{1mu}$. It also vanishes if some $\bfa_j$ is zero, since then $T_j=T_{j-1}$ is forced and $p_{\emptyset,\boldsymbol0}=0\mspace{1mu}$. Given a chain and a decomposition with a nonvanishing term, let $M$ be the matrix with columns $\bfa_1,\dots,\bfa_s\mspace{1mu}$. Then $T_j(M)=T_j$ for all~$j\mspace{1mu}$. Indeed, the inclusion $\subseteq$ holds because~$\bfa_{j'}$ is supported on $T_{j'}\subseteq T_j$ for $j'\leq j\mspace{1mu}$, and $\supseteq$ because every $i\in T_j$ lies in $T_{j'}\setminus T_{j'-1}$ for some $j'\leq j\mspace{1mu}$, where the $i$th entry of $\bfa_{j'}$ is positive. So $M\in\mathcal{M}_{2r\times s}(\bfc)\mspace{1mu}$, with $A_j(M)=T_j\setminus T_{j-1}$ and $\bfa_j(M)$ the restriction of $\bfa_j$ to~$T_j\mspace{1mu}$. Conversely, every $M\in\mathcal{M}_{2r\times s}(\bfc)$ arises in this way from exactly one pair, namely from its own chain~$(T_j(M))_j$ and its columns. Summing over $s$ and over all chains in~\eqref{eq:paths} gives the corollary.
\end{proof}

\begin{remark}
By Corollary~\ref{cor:vanishing} the summand of $M$ vanishes unless every column of $M$ has sum at least~$2$ and no column has an isolated nonzero entry, that is, an entry $m_{i,j}>0$ whose neighbours $m_{i^-,j}$ and $m_{i^+,j}$ in the rows~$i^\pm$ adjacent to $i$ within~$T_j(M)$ are both zero (with $m_{i^\pm,j}=0$ if $i$ is the smallest or the largest element of $T_j(M)$). These conditions may be added to Definition~\ref{def:dominomatrix} without changing the sum. In particular only $s\leq\frac12|\bfc|$ occurs, in accordance with the depth bound~\eqref{eq:depthbound}.

Marking the first nonzero entry of each row of~$M$ gives a matrix with entries in~$\{0,1\}\mspace{1mu}$, with exactly one entry~$1$ in every row, whose $j$th column is the indicator function of~$A_j(M)\mspace{1mu}$, a partial domino tiling of~$T_j(M)\mspace{1mu}$. Hence the sum in Corollary~\ref{cor:main} can be organised first over these matrices and then over the~$M$ with a given one, which are those with positive entries at the positions of its entries~$1$ and zeros to their left.
\end{remark}

\begin{corollary}\label{cor:depthone} For $\ell,k\geq 1\mspace{1mu}$, we have \[
\EE_{\epsilon}(\ell,k) \= -\sum_{s\geq 0} (-1)^s\!\!\sum_{ \substack{m_i,n_i\geq 1 \\ m_1+\ldots+m_s=\ell\\n_1+\ldots+n_s=k }}\!\! b_{m_1+n_1+\epsilon-1}\,b_{m_2+n_2}\cdots b_{m_s+n_s}\,.
\]
\end{corollary}
\begin{proof}
We have $r=1$ and $\mathcal{T}(1)=\{\emptyset,\{1,2\}\}\mspace{1mu}$, and the polynomials $p_{A,\bfa}$ were computed after Lemma~\ref{lem:pcoeff}. Let $M\in \mathcal{M}_{2\times s}(\ell,k)$ and write $(m_j,n_j)$ for its $j$th column. The set $T_1(M)$ is nonempty of even cardinality, so $T_1(M)=[2]\mspace{1mu}$. Hence $A_1(M)=\{1,2\}\mspace{1mu}$, $A_j(M)=\emptyset$ for $j\geq2$ and $\bfa_j(M)=(m_j,n_j)$ for all~$j\mspace{1mu}$. The summand of $M$ vanishes unless $m_j,n_j\geq1$ for all~$j\mspace{1mu}$, in which case it equals \[
y\,x^{m_1+n_1+\epsilon-2}\; \bigl(-y\,x^{m_2+n_2-1}\bigr)\cdots \bigl(-y\,x^{m_s+n_s-1}\bigr) \= (-1)^{s-1}\,b_{m_1+n_1+\epsilon-1}\,b_{m_2+n_2}\cdots b_{m_s+n_s}\msp. \qedhere
\]
\end{proof}

In particular, $\EE_{\epsilon}(\ell,k)$ is homogeneous of weight $\ell+k+\epsilon-1\mspace{1mu}$, since $b_d$ has weight~$d\mspace{1mu}$, and therefore coincides with the classical word~$\EE^{\mathrm{cl}}_{\epsilon}(\ell,k)\mspace{1mu}$. For $\epsilon=1\mspace{1mu}$, applying~$Z$ as in~\eqref{eq:drop1} gives the explicit formula of Kaneko and Sakata for multiple zeta values of height one \cite[Theorem~1.1]{KS} (see also \cite[Theorem~1.5]{HMSW}).

In depth two, for $\epsilon=0$ and $\bfc=(1,1,1,2)\mspace{1mu}$, three matrices have a nonzero summand, \[
  \begin{pmatrix}1&0\\1&0\\0&1\\0&2\end{pmatrix}\,,\qquad
  \begin{pmatrix}1&0\\0&1\\0&1\\2&0\end{pmatrix}\,,\qquad
  \begin{pmatrix}0&1\\0&1\\1&0\\2&0\end{pmatrix}\,,
\]
with $T_1(M)=\{1,2\}\mspace{1mu}$, $\{1,4\}$ and $\{3,4\}$ respectively. With the coefficients $p_{A,\bfa}$ of Section~\ref{sec:coeffs}, their summands are $y\cdot y\,x=b_1b_2\mspace{1mu}$, $y\,x\cdot(-y)=-b_2b_1$ and $y\,x\cdot y=b_2b_1\mspace{1mu}$. The last two cancel, and $\EE_{0}(1,1,1,2)=b_1b_2\mspace{1mu}$.

\section{\texorpdfstring{Explicit $q$-series identities}{Explicit q-series identities}}\label{sec:qseries}
In this section we apply the maps $Z^\dagger_{q,N}$ of Theorem~\ref{thm:HMW} to the words and their generating series, which turns the results of the preceding sections into identities between $q$-series. In Section~\ref{sec:finite} this proves Proposition~\ref{prop:Wlinear}. Under $Z^\dagger_{q,N}$ the letter~$y$ becomes the passage from one summation index to the next, the letter~$x$ becomes $1/(1-q^{n})\mspace{1mu}$, and the inverse of $1+y\,\mathsf C_r$ becomes an ordered product of matrices whose entries are explicit rational functions of $q^{n}$ and the~$t_s\mspace{1mu}$, as we show in Section~\ref{sec:product}. In Section~\ref{sec:qchains} we expand the product along chains, which gives identities of generating series, and in Section~\ref{sec:partitions} we take coefficients, which gives identities between the $\zeta^\dagger_q$ as well as identities between sums over partitions such as~\eqref{eq:derg1}.

\subsection{Finite truncations}\label{sec:finite}
We prove Proposition~\ref{prop:Wlinear} by passing to a finite level~$N\mspace{1mu}$, where \eqref{eq:HMWfinite} characterises the words. For $N\geq1$ let $\Psi^{0}_N$ be the linear map on $\bigoplus_{r\geq0}\Q[x_1,\ldots,x_r,y_1,\ldots,y_r]$ with $\Psi^{0}_N(p_{\bfc})=\zeta^\dagger_{q,N}(\bfc)\mspace{1mu}$, and let $\Psi^{1}_N$ be the linear map on $\bigoplus_{r\geq0}x_r\Q[x_1,\ldots,x_r]$ with \[
  \Psi^{1}_N\biggl(\prod_{j=1}^{r}\binom{x_j}{k_j-1}\biggr)\=\sum_{A\subseteq\{j\,:\,k_j=1\}}\ \sum_{\substack{0<n_1\leq\dots\leq n_r<N\\ n_j<n_{j+1}\text{ if }j\notin A}}\ \prod_{j\in A}\frac{q^{N-n_j}}{1-q^{N-n_j}}\ \prod_{j\notin A}\frac{q^{n_j(k_j-1)}}{(1-q^{n_j})^{k_j}}
\]
for $\bfk\in\mathbb{I}^{\mathrm{adm}}$, the right-hand side being the finite Bradley--Zhao model $\zeta^{\diamond,\mathrm{BZ}}_{q,N}(\bfk)$ of \cite[Section~1]{HMW26}. On $\bigoplus_{r\geq0}(x_1\cdots x_r)^{\epsilon}\Q[x_1,\ldots,x_r]$ the map~$\Psi^{\epsilon}_N$ is \eqref{eq:phi} with the summation restricted to $n_r<N\mspace{1mu}$, by \eqref{eq:pcmult} and \eqref{eq:BZgeq2}, since the set~$A$ is empty when all $k_j\geq2\mspace{1mu}$. In particular~$\Psi^{1}_N$ and $\Psi^{0}_N$ agree there, and $\Psi^{0}_N\circ\rho=Z^\dagger_{q,N}\mspace{1mu}$. As in Section~\ref{sec:Lmaps}, it follows that $\Psi^{\epsilon}_N=\Psi^{0}_N\circ\mathcal L^{\epsilon}$.

\begin{lemma}\label{lem:finite}
Let $\epsilon\in\{0,1\}$ and $N\geq1\mspace{1mu}$. Then
\begin{equation}\label{eq:PsiNstar}
  \Psi^{\epsilon}_N(p\ast^{\epsilon} f)\=\Psi^{\epsilon}_N(p)\,\Psi^{\epsilon}_N(f)
\end{equation}
for every~$p$ in the domain of~$\Psi^{\epsilon}_N$ and every $f\in\bigoplus_{r\geq0}(x_1\cdots x_r)^{\epsilon}\Q[x_1,\ldots,x_r]\mspace{1mu}$.
\end{lemma}

\begin{proof}For $0<a<N$ put $u_a\df q^{\epsilon(N-a)}/(1-q^{N-a})\mspace{1mu}$, and for $0\leq n'<n<N$ and $L\geq0$ let $H_L(n',n]$ be the sum of $u_{a_1}\cdots u_{a_L}$ over all $n'<a_1\leq\dots\leq a_L\leq n\mspace{1mu}$. For $\bfc=(0^{\ell_1-1},k_1,\dots,0^{\ell_r-1},k_r)$ the definition of $\zeta^\dagger_{q,N}(\bfc)$ and $q^{n}(1-q^{n})^{-k}=\sum_{m\geq1}\binom{m+k-2}{k-1}q^{mn}$ give \[
  \Psi^{0}_N(p_{\bfc})\=\sum_{\substack{0<n_1<\dots<n_r<N\\ m_1,\dots,m_r\geq1}}\ \prod_{j=1}^{r}\binom{m_j+k_j-2}{k_j-1}\,H_{\ell_j-1}(n_{j-1},n_j]\;q^{m_1n_1+\dots+m_rn_r}\qquad(n_0\df0)\mspace{1mu},
\]
the summation indices at the zeros of~$\bfc$ contributing the factors~$u_a$ and no power of~$q\mspace{1mu}$. So $\Psi^{0}_N(p)$ is obtained from \eqref{eq:phi} by writing~$p$ in the basis~\eqref{eq:pc} and replacing each factor $\binom{y_j-y_{j-1}+\ell-2}{\ell-1}\mspace{1mu}$, which is $H_{\ell-1}(n_{j-1},n_j]$ with all $u_a$ equal to~$1\mspace{1mu}$, by $H_{\ell-1}(n_{j-1},n_j]\mspace{1mu}$. For $\epsilon=1\mspace{1mu}$, expanding $p\in x_r\Q[x_1,\ldots,x_r]$ in the basis $\prod_j\binom{x_j}{k_j-1}$ and using \eqref{eq:BZgeom}, whose term $m=0$ is present precisely when $k_j=1\mspace{1mu}$, gives 
\[
  \Psi^{1}_N(p)\=\sum_{A\subseteq[r-1]}\ \sum_{\substack{0<n_1\leq\dots\leq n_r<N\\ n_j<n_{j+1}\text{ if }j\notin A}}\ \sum_{\substack{m_1,\dots,m_r\geq0\\ m_j=0\text{ if }j\in A}}p(m_1,\dots,m_r)\,q^{m_1n_1+\dots+m_rn_r}\prod_{j\in A}u_{n_j}\,,
\]
the terms with $m_r=0$ vanishing as $x_r$ divides~$p\mspace{1mu}$. Here the $n_j$ with $j\in A$ play the role of the indices~$a$ above: they contribute the factor~$u_{n_j}$ and no power of~$q\mspace{1mu}$. In both cases $f$ has no such indices: for $\epsilon=0$ it does not involve the~$y_j\mspace{1mu}$, and for $\epsilon=1$ it vanishes as soon as one of its multiplicities is~$0\mspace{1mu}$, so that only $A=\emptyset$ contributes to $\Psi^{1}_N(f)$ and $\Psi^{1}_N(f)=\Psi^{0}_N(f)\mspace{1mu}$.

Now merge a configuration for~$p$ with one for~$f$ as in the proof of Proposition~\ref{prop:star}: order the parts of both increasingly, a part common to both becoming one part with the multiplicities added, and keep the indices~$a$ of~$p$ in their intervals, an index equal to a part of~$f$ in front of that part. By induction on the number of parts, applied to the parts larger than the smallest one, the parts and multiplicities of the merged configurations are those of $p\ast^{\epsilon}f\mspace{1mu}$, a common part giving the letter product $g\diamond^{\epsilon}g'$ of the two factors. For $\epsilon=1$ its term with the multiplicity of~$f$ equal to~$0$ is absent, as $f$ vanishes there, and its term with the multiplicity of~$p$ equal to~$0$ is a position $j\notin A$ of~$p$ with $m_j=0\mspace{1mu}$. The merged indices~$a$ obey the rule above, and every term of $\Psi^{\epsilon}_N(p\ast^{\epsilon}f)$ arises exactly once, as an index~$a$ may equal the next part and a position $j\notin A$ may not. It remains to compare the factors~$u_a\mspace{1mu}$. For $\epsilon=1$ they are unchanged. For $\epsilon=0$ the parts of~$f$ that fall into an interval $(n_{j-1},n_j)$ subdivide it, and in $p\ast^{0}f$ the factor $h_j(y_j-y_{j-1})$ of~$p$ becomes~$h_j$ of the sum of the consecutive differences of parts. By linearity we may take $h_j=\binom{y+L-1}{L}\mspace{1mu}$. The expansion 
\[
  \binom{y+y'+L-1}{L}\=\sum_{i+i'=L}\binom{y+i-1}{i}\binom{y'+i'-1}{i'}
\]
over the two subintervals created by a part~$n''$ of~$f\mspace{1mu}$, followed by the replacement of each binomial coefficient by the corresponding sum~$H\mspace{1mu}$, gives $\sum_{i+i'=L}H_i(n_{j-1},n'']\,H_{i'}(n'',n_j]=H_L(n_{j-1},n_j]\mspace{1mu}$, since splitting a weakly increasing tuple at~$n''$ is a bijection; likewise for several parts. This is the factor in $\Psi^{0}_N(p)\,\Psi^{0}_N(f)\mspace{1mu}$, and \eqref{eq:PsiNstar} follows.
\end{proof}

\begin{proof}[Proof of Proposition~\ref{prop:Wlinear}]
Let $p$ be in the domain of~$\mathcal L^{\epsilon}$ and $f\in\bigoplus_{r\geq0}(x_1\cdots x_r)^{\epsilon}\Q[x_1,\ldots,x_r]\mspace{1mu}$. By $\Psi^{\epsilon}_N=\Psi^{0}_N\circ\mathcal L^{\epsilon}$ and Lemma~\ref{lem:finite} 
\[
  \Psi^{0}_N\bigl(\mathcal L^{\epsilon}(p\ast^{\epsilon} f)\bigr)\=\Psi^{\epsilon}_N(p\ast^{\epsilon} f)\=\Psi^{\epsilon}_N(p)\,\Psi^{\epsilon}_N(f)\=\Psi^{0}_N(\mathcal L^{\epsilon}p)\,\Psi^{0}_N(f)\=\Psi^{0}_N\bigl(\mathcal L^{\epsilon}(p)\ast^{0} f\bigr)
\]
for all~$N\geq1\mspace{1mu}$, as $\Psi^{\epsilon}_N$ and~$\Psi^{0}_N$ agree on~$f\mspace{1mu}$. In the last term $\ast^{0}$ may be replaced by~$\ast^{\epsilon}$, since for $\epsilon=1$ both $\mathcal L^{1}p$ and~$f$ lie in $\bigoplus_{r\geq0}x_1\cdots x_r\Q[x_1,\ldots,x_r]\mspace{1mu}$, where the two products agree by Proposition~\ref{prop:star}. Both arguments of~$\Psi^{0}_N$ lie in $\bigoplus_{r\geq0}\Q[x_1,\ldots,x_r]\mspace{1mu}$, where $\Psi^{0}_N\circ\rho=Z^\dagger_{q,N}\mspace{1mu}$, and these maps separate points as $N$ varies, by the injectivity recalled after Theorem~\ref{thm:HMW}.
\end{proof}

\subsection{The ordered product}\label{sec:product}
Read letter by letter, the truncated maps $Z^\dagger_{q,N}$ of Section~\ref{sec:HMWthm} have a simple description: for polynomials $f_1,\dots,f_s$ in~$x$ one has
\begin{equation}\label{eq:letterwise}
  Z^\dagger_{q,N}\bigl(y\,f_1(x)\;y\,f_2(x)\cdots y\,f_s(x)\bigr)
  \=\sum_{0<n_1<\dots<n_s<N}\ \prod_{j=1}^{s}\frac{q^{n_j}}{1-q^{n_j}}\,f_j\Bigl(\frac{1}{1-q^{n_j}}\Bigr)\mspace{1mu},
\end{equation}
because $yx^{c-1}$ contributes $q^{n}(1-q^{n})^{-c}$ and no entry~$c_j$ of the word is zero. We apply $Z^\dagger_{q,N}$ to series in the $t_s$ with coefficients in $\fh^1$ coefficientwise. For a finite $S\subset\Zgeq$ of even cardinality let \[
  \gs^{\epsilon}_{S,N}(q)\df\gs_{S,N}(q)\df\sum_{\bfc\in\Zgeq^{S}}\xi^{\epsilon}_{N}(\bfc;q)\,\bft^{\bfc}\,,\qquad
  \gs^{\epsilon}_{S}(q)\df\gs_{S}(q)\df\lim_{N\to\infty}\gs_{S,N}(q)=\sum_{\bfc\in\Zgeq^{S}}\xi^{\epsilon}(\bfc;q)\,\bft^{\bfc}
\]
be the generating series of the $q$-series of Theorem~\ref{thm:HMW}. In this notation \eqref{eq:HMWfinite} says
\begin{equation}\label{eq:HMWfinite-gs}
  \gs_{S,N}(q)\=Z^\dagger_{q,N}\bigl(\gs_S\bigr)\mspace{1mu}.
\end{equation}

Let $\Omega^{\epsilon}_r(q)\df\Omega_r(q)$ be the matrix $1-\frac{q}{1-q}\,\mathsf C_r\bigl(\frac{1}{1-q}\bigr)\mspace{1mu}$, where $\mathsf C_r(\frac1{1-q})$ is $\mathsf C_r$ with $x$ replaced by $1/(1-q)\mspace{1mu}$. Its entries are the power series \[
  \omega^{\epsilon}_{S\setminus A,S}(q)\df\omega_{S\setminus A,S}(q)\=\delta_{A,\emptyset}+w_{A,S}\big|_{x=1/(1-q)}\qquad(A\in\calT(S))
\]
in $q$ and the~$t_s\mspace{1mu}$, and $\omega_{T,S}(q)=0$ for all other~$T\mspace{1mu}$.

\begin{proposition}\label{prop:transfer}
For $N\geq1$ and $U\in\calU(r)\mspace{1mu}$,
\begin{equation}\label{eq:transfer}
  \gs^{\epsilon}_{U,N}(q)\=\bigl[\Omega^{\epsilon}_r(q)\,\Omega^{\epsilon}_r(q^{2})\cdots\Omega^{\epsilon}_r(q^{N-1})\bigr]_{\emptyset,U}\,.
\end{equation}
\end{proposition}

\begin{proof}
By \eqref{eq:HMWfinite-gs} and \eqref{eq:solved}, $\gs_{U,N}(q)$ is $Z^\dagger_{q,N}$ applied to $[(1+y\,\mathsf C_r)^{-1}]_{\emptyset,U}\mspace{1mu}$. By Remark~\ref{rem:otherentries}, this entry is a sum of words $y\,f_1(x)\cdots y\,f_s(x)\mspace{1mu}$, one for each chain $\emptyset=T_0\subseteq T_1\subseteq\dots\subseteq T_s=U\mspace{1mu}$, with $f_j=-\mathsf C_{T_{j-1},T_j}\mspace{1mu}$. Each coefficient of $\bft^{\bfc}$ receives contributions from finitely many chains only, since the entries of $\mathsf C_r$ have positive order in the~$t_s\mspace{1mu}$. By \eqref{eq:letterwise}, applying $Z^\dagger_{q,N}$ gives \[
  \sum_{s\geq0}\ \sum_{0<n_1<\dots<n_s<N}
  \Bigl[\prod_{j=1}^{s}\Bigl(-\frac{q^{n_j}}{1-q^{n_j}}\,\mathsf C_r\Bigl(\frac{1}{1-q^{n_j}}\Bigr)\Bigr)\Bigr]_{\emptyset,U}\,,
\]
which is the expansion of the entry of the product in~\eqref{eq:transfer}.
\end{proof}

Letting $N\to\infty$ gives the $q$-series form of Theorem~\ref{thm:intro}.

\begin{corollary}\label{cor:product}
For $\epsilon\in\{0,1\}$ and $r\geq1\mspace{1mu}$, \[
  \gs^{\epsilon}_{[2r]}(q)\=\bigl[\Omega^{\epsilon}_r(q)\,\Omega^{\epsilon}_r(q^{2})\,\Omega^{\epsilon}_r(q^{3})\cdots\bigr]_{\emptyset,[2r]}\,.
\]
\end{corollary}

\begin{proof}
The $n$th factor of the product is $\equiv1\pmod{q^{n}}\mspace{1mu}$, so the product converges $q$-adically, to the limit of the right-hand side of~\eqref{eq:transfer} as $N\to\infty\mspace{1mu}$.
\end{proof}

\begin{remark}\label{rem:HMWrec}
Proposition~\ref{prop:transfer} says that the row vector $\gs_{\bullet,N}(q)=(\gs_{U,N}(q))_{U\in\calU(r)}$ satisfies \[
  \gs_{\bullet,N+1}(q)\=\gs_{\bullet,N}(q)\,\Omega_r(q^{N})\mspace{1mu}.
\]
This is, after the substitution $t_s=1-q^{W_s}$ with $(W_1,\dots,W_{2r})=(Y_1,X_1,\dots,Y_r,X_r)$ and clearing denominators, the recurrence in~$N$ of \cite[Corollary~3.5]{HMW26}. Hirose, Maesaka and Watanabe obtain it by manipulating the finite sums, and with the computation in the proof of Proposition~\ref{prop:transfer} and Theorem~\ref{thm:intro} it gives back Theorem~\ref{thm:HMW}.
\end{remark}

To write the entries of $\Omega_r(q^{n})\mspace{1mu}$, for $n\geq1\mspace{1mu}$, as rational functions of $q^{n}$ and the~$t_s\mspace{1mu}$, put, for series $u,v$ without constant term, \[
  \langle u\rangle_n\defis 1-q^{n}-u\,,\qquad
  \langle u,v\rangle_n\defis(1-u)(1-v)-q^{n}\,,
\]
so that $\langle u\rangle_n=\langle u,0\rangle_n\mspace{1mu}$, and $\langle u,v\rangle_n$ is, up to the factor $1-q^{n}$, the polynomial $1-x(u+v)+x\,uv$ of Lemma~\ref{lem:alphabeta} at $x=1/(1-q^{n})\mspace{1mu}$. For a domino $D=\{u,v\}$ write $\langle t_D\rangle_n\defis\langle t_u,t_v\rangle_n\mspace{1mu}$, which is symmetric in $u$ and~$v\mspace{1mu}$.

\begin{lemma}\label{lem:omega}
Let $n\geq1\mspace{1mu}$, let $S\in\calU(r)$ be nonempty, $A\in\calT(S)$ and $T=S\setminus A\mspace{1mu}$. Then \[
  \omega^{\epsilon}_{T,S}(q^{n})\=(-1)^{\mathrm{eo}(A)}\,q^{n\kappa(A)}\,(1-q^{n})^{1-\epsilon\kappa(A)}\;\bft^{\boldsymbol e(A)}\,
  \frac{\prod_{i\ \mathrm{odd}}\langle t_{D_i(T)}\rangle_n}
       {\langle t_{\min S}\rangle_n\,\langle t_{\max S}\rangle_n\prod_{i\ \mathrm{even}}\langle t_{D_i(S)}\rangle_n}\,,
\]
and $\omega^{\epsilon}_{\emptyset,\emptyset}(q^{n})=1\mspace{1mu}$.
\end{lemma}

\begin{proof}
At $x=1/(1-q^{n})$ one has $x-1=q^{n}/(1-q^{n})$ and $\hat t_s=t_s(1-q^{n})/\langle t_s\rangle_n\mspace{1mu}$, hence \[
  \eta_{D}\=\frac{q^{n}\,t_ut_v}{\langle t_u\rangle_n\langle t_v\rangle_n}\,,
  \qquad
  1-\eta_{D}\=\frac{(1-q^{n})\,\langle t_D\rangle_n}{\langle t_u\rangle_n\langle t_v\rangle_n}
\]
for every domino $D=\{u,v\}\mspace{1mu}$, the second because $\langle t\rangle_n\langle t'\rangle_n-q^{n}tt'=(1-q^{n})\bigl((1-q^{n})-(t+t')+tt'\bigr)\mspace{1mu}$. Therefore
\begin{align*}
  \Pp_T&\=(1-q^{n})^{|T|/2}\,\frac{\prod_{i\ \mathrm{odd}}\langle t_{D_i(T)}\rangle_n}{\prod_{s\in T}\langle t_s\rangle_n}\,,\qquad
  \Pp^{*}_S\=(1-q^{n})^{|S|/2-1}\,\frac{\prod_{i\ \mathrm{even}}\langle t_{D_i(S)}\rangle_n}{\prod_{s\in S\setminus\{\min S,\max S\}}\langle t_s\rangle_n}\,,\\
  \eta_A&\=q^{n\kappa(A)}\prod_{a\in A}\frac{t_a}{\langle t_a\rangle_n}\,.
\end{align*}
We insert this and $x^{(\epsilon-1)\kappa(A)}=(1-q^{n})^{(1-\epsilon)\kappa(A)}$ into $\delta_{T,S}+w_{A,S}\mspace{1mu}$. Then the factors~$\langle t_s\rangle_n$ with $s\in T\sqcup A=S$ cancel against those with $s\in S\setminus\{\min S,\max S\}\mspace{1mu}$, and the exponent of $1-q^{n}$ is $\frac12|T|-(\frac12|S|-1)+(1-\epsilon)\kappa(A)=1-\epsilon\kappa(A)\mspace{1mu}$, as $|S|-|T|=2\kappa(A)\mspace{1mu}$.
\end{proof}

In particular $\omega^{\epsilon}_{T,S}(q^{n})\equiv0\pmod{q^{n\kappa(A)}}$ for $T\neq S\mspace{1mu}$, and $\omega^{\epsilon}_{S,S}(q^{n})\equiv1\pmod{q^{n}}\mspace{1mu}$. For $\epsilon=1$ every column of $\Omega^{1}_r(q^{n})$ sums to~$1\mspace{1mu}$, by Proposition~\ref{prop:collapse}.

\subsection{Sums over chains}\label{sec:qchains}
Expanding the corner entry of the product in Corollary~\ref{cor:product} along the strict chains of \eqref{eq:chains} gives the $q$-series form of that identity. Here the matrices are triangular with respect to inclusion, and the steps on the diagonal are collected into products. 

\begin{corollary}\label{cor:qchains}
For $\epsilon\in\{0,1\}$ and $r\geq1\mspace{1mu}$, \[
  \gs^{\epsilon}_{[2r]}(q)\=\sum_{s\geq1}\ \sum_{\emptyset=U_0\subsetneq U_1\subsetneq\dots\subsetneq U_s=[2r]}\
  \sum_{0<m_1<\dots<m_s}\ \prod_{i=1}^{s}\Bigl(\omega^{\epsilon}_{U_{i-1},U_i}(q^{m_i})
  \prod_{m_i<n<m_{i+1}}\omega^{\epsilon}_{U_i,U_i}(q^{n})\Bigr)\mspace{1mu},
\]
where we put $m_{s+1}\defis\infty\mspace{1mu}$ and the second sum runs over the strict chains from $\emptyset$ to~$[2r]\mspace{1mu}$.
\end{corollary}

For $r=1$ there is one strict chain, and the corollary reads
\begin{equation}\label{eq:depthone-q}
  \gs^{\epsilon}_{[2]}(q)
  \=\sum_{m\geq1}\frac{q^{m}(1-q^{m})^{1-\epsilon}\,t_1t_2}{\langle t_1\rangle_m\langle t_2\rangle_m}
  \prod_{n>m}\frac{(1-q^{n})\,\langle t_1,t_2\rangle_n}{\langle t_1\rangle_n\langle t_2\rangle_n}\,,
\end{equation}
whose right-hand side is symmetric in $t_1$ and~$t_2\mspace{1mu}$. For $\epsilon=0$ this symmetry is the duality $\zeta^\dagger_q(0^{\ell-1},k)=\zeta^\dagger_q(0^{k-1},\ell)\mspace{1mu}$.

\subsection{\texorpdfstring{Quasimodular forms and $q$-series identities}{Quasimodular forms and q-series identities}}\label{sec:partitions}
Applying~$Z^\dagger_q$ to the words of Section~\ref{sec:words} and using Theorem~\ref{thm:HMW} gives identities between the $\zeta^\dagger_q\mspace{1mu}$. In depth one, Corollary~\ref{cor:depthone} gives the following.

\begin{corollary}\label{cor:depthone-q}
For $\ell,k\geq1\mspace{1mu}$, \[
  \zeta^\dagger_q(0^{\ell-1},k)\=-\sum_{s\geq 0}(-1)^s\sum_{\substack{k_i,\ell_i\geq 1 \\ k_1+\ldots+k_s=k\\\ell_1+\ldots+\ell_s=\ell}}\zeta_q^\dagger(k_1+\ell_1-1,k_2+\ell_2,\ldots,k_s+\ell_s)
\,.
\]
\end{corollary}

\begin{proof}
Apply~$Z^\dagger_q$ to the identity in Corollary~\ref{cor:depthone}. By Theorem~\ref{thm:HMW}, $Z^\dagger_q(\EE_{0}(\ell,k))=\zeta^\dagger_q(0^{\ell-1},k)\mspace{1mu}$, and $Z^\dagger_q(b_{c_1}\cdots b_{c_s})=\zeta^\dagger_q(c_1,\dots,c_s)\mspace{1mu}$.
\end{proof}

In depth two, Corollary~\ref{cor:main} expresses each $\zeta^\dagger_q(0^{\ell_1-1},k_1,0^{\ell_2-1},k_2)$ as a $\ZZ$-linear combination of the $\zeta^\dagger_q(\bfd)\mspace{1mu}$. Unlike in depth one, the combination is in general not homogeneous in the weight and its coefficients are not always $\pm1\mspace{1mu}$, because the inverse of $\Pp^{*}_{[4]}=1-\eta_{\{2,3\}}$ brings in the polynomials~$\gamma$ of Section~\ref{sec:coeffs}. For example,
\begin{align*}
  \zeta^\dagger_q(0,1,0,1)&\=\zeta^\dagger_q(2,2)\mspace{1mu},\\
  \zeta^\dagger_q(0,2,1)&\=\zeta^\dagger_q(3,1)-\zeta^\dagger_q(2,2)-\zeta^\dagger_q(1,2,1)-\zeta^\dagger_q(1,1,2)+\zeta^\dagger_q(1,2)\mspace{1mu},\\
  \zeta^\dagger_q(0,2,0,2)&\=-3\zeta^\dagger_q(3,1,2)+2\zeta^\dagger_q(2,4)+4\zeta^\dagger_q(2,2,2)+\zeta^\dagger_q(1,3,2)+3\zeta^\dagger_q(1,2,1,2)\\
  &\qquad+\zeta^\dagger_q(1,1,4)+2\zeta^\dagger_q(1,1,2,2)+\zeta^\dagger_q(3,2)-2\zeta^\dagger_q(2,3)-\zeta^\dagger_q(1,4)\\
  &\qquad-4\zeta^\dagger_q(1,2,2)-\zeta^\dagger_q(1,1,3)+\zeta^\dagger_q(1,3)\mspace{1mu}.
\end{align*}

Corollary~\ref{cor:main} makes the map~$\mathcal L^{0}$ of the introduction explicit: applying~$\rho$ to the words turns it into an identity between $q$-series in which the parts have disappeared from the right-hand side. For $\bfc\in\Zgeq^{2r}$ let $\lambda_{\bfd}\in\ZZ$ be the coefficients of $\EE_{0}(\bfc)$ in the basis of words, \[
  \EE_{0}(\bfc)\=\sum_{s\geq0}\ \sum_{\bfd\in\Zgeq^{s}}\lambda_{\bfd}\,b_{d_1}\cdots b_{d_s}\,,
\]
which Corollary~\ref{cor:main} computes. Written out, Theorem~\ref{thm:HMW} says that for $\boldsymbol\ell,\bfk\in\Zgeq^{r}$ and $\bfc=(\ell_1,k_1,\dots,\ell_r,k_r)\mspace{1mu}$, with $n_0\defis0\mspace{1mu}$,
\begin{multline}\label{eq:partitions}
\sum_{\substack{0<n_1<\dots<n_r\\ m_1,\dots,m_r\geq1}}\ \prod_{j=1}^{r}\binom{m_j+k_j-2}{k_j-1}\binom{n_j-n_{j-1}+\ell_j-2}{\ell_j-1}q^{m_jn_j}\\
\=\sum_{\bfd}\lambda_{\bfd}\!\!\sum_{\substack{0<n_1<\dots<n_s\\ m_1,\dots,m_s\geq1}}\ \prod_{j=1}^{s}\binom{m_j+d_j-2}{d_j-1}q^{m_jn_j}\,.
\end{multline}
Indeed, by \eqref{eq:pc} the left-hand side is $\Psi^{0}(p_{\bfc'})=\zeta^\dagger_q(\bfc')$ with $\bfc'=(0^{\ell_1-1},k_1,\dots,0^{\ell_r-1},k_r)\mspace{1mu}$, and by \eqref{eq:pcmult} the right-hand side is $\Psi^{0}_X(\rho(\EE_{0}(\bfc)))=Z^\dagger_q(\EE_{0}(\bfc))\mspace{1mu}$.

For $r=1\mspace{1mu}$, \eqref{eq:partitions} with the coefficients in Corollary~\ref{cor:depthone-q} reads
\begin{multline*}
\sum_{m,n\geq 1} \binom{m+k-2}{k-1}\binom{n+\ell-2}{\ell-1} q^{mn}\\
\=-\sum_{s\geq 0} (-1)^s\!\! \sum_{ \substack{0<n_1<\ldots<n_s \\ m_1,\ldots,m_s>0}} \sum_{ \substack{k_i,\ell_i\geq 1 \\ k_1+\ldots+k_s=k\\\ell_1+\ldots+\ell_s=\ell}} \prod_{j=1}^s \binom{m_j +k_j+\ell_j-2-\delta_{j,1}}{m_j-1} q^{n_jm_j}\,.
\end{multline*}
The products of binomial coefficients form a basis of $\Q[x,y]\mspace{1mu}$, so this expresses, for every $p\in\Q[x,y]\mspace{1mu}$, the series $\sum_{m,n\geq1}p(m,n)\,q^{mn}$ as $\Psi^{0}_X$ of an explicit polynomial in the multiplicities alone. This is the map~$\mathcal L^{0}$ of the introduction in depth one. In the monomial basis it takes the following form, in which $\mathcal A_k$ is the Eulerian polynomial
\[
  \mathcal A_0(q)\=\mathcal A_1(q)\=1\,,\quad \mathcal A_2(q)\=1+q\,,\quad \mathcal A_3(q)\=1+4q+q^{2}\,,\quad \mathcal A_4(q)\=1+11q+11q^{2}+q^{3}\,,
\]
given by
\[
  \sum_{m\geq1}m^{k}q^{m}\=\frac{q\mathcal A_k(q)}{(1-q)^{k+1}}\,,
\]
and $p^{(a,b)}$ denotes the derivative $\partial_x^a\partial_y^bp\mspace{1mu}$.

\begin{corollary}\label{cor:eulerian}
For every $p\in\Q[x,y]\mspace{1mu}$, \[
\sum_{m,n\geq 1}p(m,n) q^{mn}\=-\sum_{s\geq 1} (-1)^s\!\! \sum_{ \substack{0<n_1<\ldots<n_s }} \sum_{\substack{\boldsymbol{k},\boldsymbol{\ell}\in \NN^s\\ k_j,\ell_j \geq 1 \text{ if } j>1}}\!\! p^{(|\boldsymbol{k}|,|\boldsymbol{\ell}|)}(0,0) \prod_{j=1}^s \frac{q^{n_j}\mathcal A_{k_j}(q^{n_j})\mathcal A_{\ell_j}(q^{n_j})}{k_j!\,\ell_j!\,(1-q^{n_j})^{k_j+\ell_j+\delta_{j,1}}}\,.
\]
\end{corollary}

\begin{proof}
Both sides are linear in~$p\mspace{1mu}$, so let $p=x^ay^b$. Multiply the display above by $u^{\ell}v^{k}$ and sum over $k,\ell\geq1\mspace{1mu}$. Since $\sum_{k\geq1}\binom{m+k-2}{k-1}v^{k}=v(1-v)^{-m}$ and $\sum_{k\geq1}v^{k}(1-q_j)^{-k}=v/(1-q_j-v)\mspace{1mu}$, this gives \[
  \sum_{m,n\geq1}\frac{uv\,q^{mn}}{(1-v)^{m}(1-u)^{n}}\=-\sum_{s\geq1}(-1)^s\sum_{0<n_1<\dots<n_s}\ \prod_{j=1}^{s}\frac{q_j(1-q_j)^{\delta_{j,1}}\,uv}{(1-q_j-u)(1-q_j-v)}\,,\qquad q_j\defis q^{n_j}\,.
\]
Put $v=1-e^{-X}$ and $u=1-e^{-Y}$, and divide by~$uv\mspace{1mu}$. The left-hand side becomes $\sum_{m,n\geq1}e^{X m+Y n}q^{mn}$, in which the coefficient of $X^aY^b/(a!\,b!)$ is $\sum_{m,n\geq1}m^an^bq^{mn}$. On the right, \[
  \frac{uv}{(1-q_j-u)(1-q_j-v)}\=\frac{(e^{X}-1)(e^{Y}-1)}{(1-q_je^{X})(1-q_je^{Y})}\,,
\]
so by the generating series \[
  \sum_{k\geq0}\frac{X^k}{k!}\,\frac{q_j\mathcal A_k(q_j)}{(1-q_j)^{k+1}}\=\frac{q_je^{X}}{1-q_je^{X}}
\]
the factor with $j>1$ is the sum over $k,\ell\geq1\mspace{1mu}$, and the factor with $j=1\mspace{1mu}$, divided by $uv\mspace{1mu}$, the sum over $k,\ell\geq0\mspace{1mu}$, of $X^kY^\ell\,q_j\mathcal A_k(q_j)\mathcal A_\ell(q_j)(1-q_j)^{-k-\ell-\delta_{j,1}}/(k!\,\ell!)\mspace{1mu}$. Comparing the coefficients of $X^aY^b/(a!\,b!)$ gives the claim.
\end{proof}

For $p=x^{a}y$ the corollary takes a short form, in which $\sum_{j=1}^{x}j^{a}$ denotes the polynomial in~$x$ with these values (see also \cite[Proposition~4.6]{Ba}). 

\begin{corollary}\label{cor:Lxay}
For $a\geq0\mspace{1mu}$, \[
  \mathcal L^{0}(x^{a}y)\=\sum_{j=1}^{x_1}j^{a}\oplus\bigl(x_1^{a}-(x_1+x_2)^{a}\bigr)\mspace{1mu},
\]
that is, \[
  \sum_{m,n\geq1}m^{a}n\,q^{mn}\=\Psi^{0}_X\Bigl(\sum_{j=1}^{x_1}j^{a}\Bigr)-\Psi^{0}_X\bigl((x_1+x_2)^{a}-x_1^{a}\bigr)\mspace{1mu}.
\]
\end{corollary}

\begin{proof}
For $p=x^{a}y$ only the terms with $s\leq2$ occur in Corollary~\ref{cor:eulerian}, as $|\boldsymbol\ell|=1\mspace{1mu}$. The term with $s=1$ is $\sum_{n\geq1}(1-q^{n})^{-1}\sum_{m\geq1}m^{a}q^{mn}=\sum_{m,n\geq1}\bigl(\sum_{j=1}^{m}j^{a}\bigr)q^{mn}$, and the terms with $s=2$ add up to $-\sum\binom{a}{k_1}\Psi^{0}_X(x_1^{k_1}x_2^{k_2})$ over $k_1+k_2=a$ with $k_2\geq1\mspace{1mu}$.
\end{proof}

For $a=1$ this is \eqref{eq:derg1}. For $a=2\mspace{1mu}$, $4$ and~$6\mspace{1mu}$,
\begin{align*}
  \mathcal L^{0}(x^{2}y)&\=\bigl(\tfrac13x_1^{3}+\tfrac12x_1^{2}+\tfrac16x_1\bigr)\oplus\bigl(-2x_1x_2-x_2^{2}\bigr)\mspace{1mu},\\
  \mathcal L^{0}(x^{4}y)&\=\bigl(\tfrac15x_1^{5}+\tfrac12x_1^{4}+\tfrac13x_1^{3}-\tfrac1{30}x_1\bigr)\oplus\bigl(-4x_1^{3}x_2-6x_1^{2}x_2^{2}-4x_1x_2^{3}-x_2^{4}\bigr)\mspace{1mu},\\
  \mathcal L^{0}(x^{6}y)&\=\bigl(\tfrac17x_1^{7}+\tfrac12x_1^{6}+\tfrac12x_1^{5}-\tfrac16x_1^{3}+\tfrac1{42}x_1\bigr)\\&\qquad\oplus\bigl(-6x_1^{5}x_2-15x_1^{4}x_2^{2}-20x_1^{3}x_2^{3}-15x_1^{2}x_2^{4}-6x_1x_2^{5}-x_2^{6}\bigr)\mspace{1mu}.
\end{align*}
The left-hand sides of Ramanujan's equations~\eqref{eq:ramanujandiffeq} are $-24\sum_{m,n\geq1}m^{2}n\,q^{mn}$, $240\sum_{m,n\geq1}m^{4}n\,q^{mn}$ and $-504\sum_{m,n\geq1}m^{6}n\,q^{mn}$. Computing their right-hand sides with the multiplication rule of Proposition~\ref{prop:star}, one obtains the same three series as $\Psi^{0}_X$ of the polynomials
\begin{gather*}
  \bigl(\tfrac12x_1+\tfrac12x_1^{3}\bigr)\oplus\bigl(-4x_1x_2\bigr)\mspace{1mu},\qquad
  \bigl(-\tfrac{3}{10}x_1+x_1^{3}+\tfrac{3}{10}x_1^{5}\bigr)\oplus\bigl(-8x_1^{3}x_2-8x_1x_2^{3}\bigr)\mspace{1mu},\\
  \bigl(-\tfrac{103}{98}x_1+\tfrac37x_1^{3}+\tfrac32x_1^{5}+\tfrac{6}{49}x_1^{7}\bigr)\oplus\bigl(-12x_1^{5}x_2+\tfrac{800}{7}x_1^{3}x_2^{3}-12x_1x_2^{5}\bigr)\mspace{1mu}.
\end{gather*}
These differ from the polynomials~$\mathcal L^{0}(x^{a}y)$ above, that is, the differences
{\allowdisplaybreaks
\begin{gather*}
  \bigl(-\tfrac13x_1+\tfrac12x_1^{2}-\tfrac16x_1^{3}\bigr)\oplus\bigl(2x_1x_2-x_2^{2}\bigr)\mspace{1mu},\\
  \bigl(\tfrac4{15}x_1-\tfrac23x_1^{3}+\tfrac12x_1^{4}-\tfrac1{10}x_1^{5}\bigr)\oplus\bigl(4x_1^{3}x_2-6x_1^{2}x_2^{2}+4x_1x_2^{3}-x_2^{4}\bigr)\mspace{1mu},\\
  \bigl(\tfrac{158}{147}x_1-\tfrac{25}{42}x_1^{3}-x_1^{5}+\tfrac12x_1^{6}+\tfrac1{49}x_1^{7}\bigr)\oplus\bigl(6x_1^{5}x_2-15x_1^{4}x_2^{2}-\tfrac{940}{7}x_1^{3}x_2^{3}-15x_1^{2}x_2^{4}+6x_1x_2^{5}-x_2^{6}\bigr)
\end{gather*}}%
lie in $\ker\Psi^{0}_X\mspace{1mu}$. For $p$ of odd degree the series~$\Psi^{0}(p)$ thus has two expressions as $\Psi^{0}_X$ of a polynomial in the multiplicities, one from its quasimodularity through Ramanujan's equations and Proposition~\ref{prop:star}, the other from Theorem~\ref{thm:intro}, and their difference is a relation between multiple~$q$-zeta values.

\appendix
\section{First values}\label{app:table}
Table~\ref{tab:values} lists~$\EE_{0}(\bfc)$ and $\EE_{1}(\bfc)$ in depths one and two, written in the letters $b_k=yx^{k-1}$ of Section~\ref{sec:setup} and computed from Definition~\ref{def:Eq} by computer. By \eqref{eq:reflection} only one of $\bfc=(c_1,\dots,c_{2r})$ and its reversal $(c_{2r},\dots,c_1)$ is listed. In every entry the largest weight and the largest depth of the words equal the bounds~\eqref{eq:weightbound} and \eqref{eq:depthbound}, so these bounds are sharp. Table~\ref{tab:L} lists the polynomials~$\mathcal L^{0}(x^{a}y^{b})$ of Section~\ref{sec:partitions}, written by depth. By Lemma~\ref{lem:Lprops}\ref{it:Lprops-iota}, $\mathcal L^{0}(x^{a}y^{b})=\mathcal L^{0}(x^{b}y^{a})\mspace{1mu}$.

{\footnotesize
\begin{longtable}{@{}l@{\quad}>{\raggedright\arraybackslash}p{0.405\textwidth}@{\quad}>{\raggedright\arraybackslash}p{0.405\textwidth}@{}}
\caption{The words $\EE_{0}(\bfc)$ and $\EE_{1}(\bfc)$ in depth at most two.}\label{tab:values}\\[-3pt]
\toprule
$\bfc$ & $\EE_{0}(\bfc)$ & $\EE_{1}(\bfc)$\\
\midrule
\endfirsthead
\toprule
$\bfc$ & $\EE_{0}(\bfc)$ & $\EE_{1}(\bfc)$\\
\midrule
\endhead
\bottomrule
\endlastfoot
\multicolumn{3}{@{}l}{\emph{Depth one}, $|\bfc|\leq8$}\\[2pt]
$(1,1)$ & $b_{1}$ & $b_{2}$ \\[1pt]
$(1,2)$ & $b_{2}$ & $b_{3}$ \\[1pt]
$(1,3)$ & $b_{3}$ & $b_{4}$ \\[1pt]
$(2,2)$ & $b_{3}-b_{1}b_{2}$ & $b_{4}-b_{2}^{2}$ \\[1pt]
$(1,4)$ & $b_{4}$ & $b_{5}$ \\[1pt]
$(2,3)$ & $b_{4}-b_{1}b_{3}-b_{2}^{2}$ & $b_{5}-b_{2}b_{3}-b_{3}b_{2}$ \\[1pt]
$(1,5)$ & $b_{5}$ & $b_{6}$ \\[1pt]
$(2,4)$ & $b_{5}-b_{1}b_{4}-b_{2}b_{3}-b_{3}b_{2}$ & $b_{6}-b_{2}b_{4}-b_{3}^{2}-b_{4}b_{2}$ \\[1pt]
$(3,3)$ & $b_{5}-b_{1}b_{4}-2b_{2}b_{3}-b_{3}b_{2}+b_{1}b_{2}^{2}$ & $b_{6}-b_{2}b_{4}-2b_{3}^{2}-b_{4}b_{2}+b_{2}^{3}$ \\[1pt]
$(1,6)$ & $b_{6}$ & $b_{7}$ \\[1pt]
$(2,5)$ & $b_{6}-b_{1}b_{5}-b_{2}b_{4}-b_{3}^{2}-b_{4}b_{2}$ & $b_{7}-b_{2}b_{5}-b_{3}b_{4}-b_{4}b_{3}-b_{5}b_{2}$ \\[1pt]
$(3,4)$ & $b_{6}-b_{1}b_{5}-2b_{2}b_{4}-2b_{3}^{2}-b_{4}b_{2}+b_{1}b_{2}b_{3}+b_{1}b_{3}b_{2}+b_{2}^{3}$ & $b_{7}-b_{2}b_{5}-2b_{3}b_{4}-2b_{4}b_{3}-b_{5}b_{2}+b_{2}^{2}b_{3}+b_{2}b_{3}b_{2}+b_{3}b_{2}^{2}$ \\[1pt]
$(1,7)$ & $b_{7}$ & $b_{8}$ \\[1pt]
$(2,6)$ & $b_{7}-b_{1}b_{6}-b_{2}b_{5}-b_{3}b_{4}-b_{4}b_{3}-b_{5}b_{2}$ & $b_{8}-b_{2}b_{6}-b_{3}b_{5}-b_{4}^{2}-b_{5}b_{3}-b_{6}b_{2}$ \\[1pt]
$(3,5)$ & $b_{7}-b_{1}b_{6}-2b_{2}b_{5}-2b_{3}b_{4}-2b_{4}b_{3}-b_{5}b_{2}+b_{1}b_{2}b_{4}+b_{1}b_{3}^{2}+b_{1}b_{4}b_{2}+b_{2}^{2}b_{3}+b_{2}b_{3}b_{2}+b_{3}b_{2}^{2}$ & $b_{8}-b_{2}b_{6}-2b_{3}b_{5}-2b_{4}^{2}-2b_{5}b_{3}-b_{6}b_{2}+b_{2}^{2}b_{4}+b_{2}b_{3}^{2}+b_{2}b_{4}b_{2}+b_{3}b_{2}b_{3}+b_{3}^{2}b_{2}+b_{4}b_{2}^{2}$ \\[1pt]
$(4,4)$ & $b_{7}-b_{1}b_{6}-2b_{2}b_{5}-3b_{3}b_{4}-2b_{4}b_{3}-b_{5}b_{2}+b_{1}b_{2}b_{4}+2b_{1}b_{3}^{2}+b_{1}b_{4}b_{2}+2b_{2}^{2}b_{3}+2b_{2}b_{3}b_{2}+b_{3}b_{2}^{2}-b_{1}b_{2}^{3}$ & $b_{8}-b_{2}b_{6}-2b_{3}b_{5}-3b_{4}^{2}-2b_{5}b_{3}-b_{6}b_{2}+b_{2}^{2}b_{4}+2b_{2}b_{3}^{2}+b_{2}b_{4}b_{2}+2b_{3}b_{2}b_{3}+2b_{3}^{2}b_{2}+b_{4}b_{2}^{2}-b_{2}^{4}$ \\[1pt]
\midrule
\multicolumn{3}{@{}l}{\emph{Depth two}, $|\bfc|\leq8$}\\[2pt]
$(1,1,1,1)$ & $b_{1}^{2}$ & $b_{2}^{2}$ \\[1pt]
$(1,1,1,2)$ & $b_{1}b_{2}$ & $b_{2}b_{3}$ \\[1pt]
$(1,1,2,1)$ & $b_{2}b_{1}$ & $b_{3}b_{2}$ \\[1pt]
$(1,1,1,3)$ & $b_{1}b_{3}$ & $b_{2}b_{4}$ \\[1pt]
$(1,1,2,2)$ & $-b_{2}^{2}+b_{3}b_{1}+b_{1}b_{2}-b_{1}^{2}b_{2}-b_{1}b_{2}b_{1}$ & $-b_{3}^{2}+b_{4}b_{2}+b_{2}b_{3}-2b_{2}^{3}$ \\[1pt]
$(1,1,3,1)$ & $b_{3}b_{1}$ & $b_{4}b_{2}$ \\[1pt]
$(1,2,1,2)$ & $b_{2}^{2}$ & $b_{3}^{2}$ \\[1pt]
$(1,2,2,1)$ & $2b_{2}^{2}-b_{1}b_{2}+b_{1}^{2}b_{2}$ & $2b_{3}^{2}-b_{2}b_{3}+b_{2}^{3}$ \\[1pt]
$(2,1,1,2)$ & $b_{1}b_{3}+2b_{2}^{2}-b_{3}b_{1}-b_{1}b_{2}+b_{1}b_{2}b_{1}$ & $b_{2}b_{4}+2b_{3}^{2}-b_{4}b_{2}-b_{2}b_{3}+b_{2}^{3}$ \\[1pt]
$(1,1,1,4)$ & $b_{1}b_{4}$ & $b_{2}b_{5}$ \\[1pt]
$(1,1,2,3)$ & $-b_{2}b_{3}-b_{3}b_{2}+b_{4}b_{1}+b_{1}b_{3}+b_{2}^{2}-b_{1}^{2}b_{3}-b_{1}b_{2}^{2}-b_{1}b_{3}b_{1}-b_{2}^{2}b_{1}$ & $-b_{3}b_{4}-b_{4}b_{3}+b_{5}b_{2}+b_{2}b_{4}+b_{3}^{2}-b_{2}^{2}b_{3}-2b_{2}b_{3}b_{2}-b_{3}b_{2}^{2}$ \\[1pt]
$(1,1,3,2)$ & $-2b_{2}b_{3}+b_{4}b_{1}+b_{1}b_{3}+b_{2}^{2}-b_{1}^{2}b_{3}-b_{1}b_{3}b_{1}-b_{2}b_{1}b_{2}-b_{2}^{2}b_{1}$ & $-2b_{3}b_{4}+b_{5}b_{2}+b_{2}b_{4}+b_{3}^{2}-b_{2}^{2}b_{3}-b_{2}b_{3}b_{2}-2b_{3}b_{2}^{2}$ \\[1pt]
$(1,1,4,1)$ & $b_{4}b_{1}$ & $b_{5}b_{2}$ \\[1pt]
$(1,2,1,3)$ & $b_{2}b_{3}$ & $b_{3}b_{4}$ \\[1pt]
$(1,2,2,2)$ & $b_{3}b_{2}-b_{2}b_{1}b_{2}$ & $b_{4}b_{3}-b_{2}^{2}b_{3}+b_{2}b_{3}b_{2}-b_{3}b_{2}^{2}$ \\[1pt]
$(1,2,3,1)$ & $2b_{2}b_{3}+b_{3}b_{2}-b_{1}b_{3}-b_{2}^{2}+b_{1}^{2}b_{3}+b_{2}b_{1}b_{2}$ & $2b_{3}b_{4}+b_{4}b_{3}-b_{2}b_{4}-b_{3}^{2}+b_{2}^{2}b_{3}+b_{3}b_{2}^{2}$ \\[1pt]
$(1,3,1,2)$ & $b_{3}b_{2}$ & $b_{4}b_{3}$ \\[1pt]
$(2,1,1,3)$ & $b_{1}b_{4}+2b_{2}b_{3}+b_{3}b_{2}-b_{4}b_{1}-b_{1}b_{3}-b_{2}^{2}+b_{1}b_{3}b_{1}+b_{2}^{2}b_{1}$ & $b_{2}b_{5}+2b_{3}b_{4}+b_{4}b_{3}-b_{5}b_{2}-b_{2}b_{4}-b_{3}^{2}+b_{2}b_{3}b_{2}+b_{3}b_{2}^{2}$ \\[1pt]
$(2,1,2,2)$ & $b_{2}b_{3}-b_{1}b_{2}^{2}$ & $b_{3}b_{4}-b_{2}b_{3}b_{2}$ \\[1pt]
$(1,1,1,5)$ & $b_{1}b_{5}$ & $b_{2}b_{6}$ \\[1pt]
$(1,1,2,4)$ & $-b_{2}b_{4}-b_{3}^{2}-b_{4}b_{2}+b_{5}b_{1}+b_{1}b_{4}+b_{2}b_{3}+b_{3}b_{2}-b_{1}^{2}b_{4}-b_{1}b_{2}b_{3}-b_{1}b_{3}b_{2}-b_{1}b_{4}b_{1}-b_{2}b_{3}b_{1}-b_{3}b_{2}b_{1}$ & $-b_{3}b_{5}-b_{4}^{2}-b_{5}b_{3}+b_{6}b_{2}+b_{2}b_{5}+b_{3}b_{4}+b_{4}b_{3}-b_{2}^{2}b_{4}-b_{2}b_{3}^{2}-2b_{2}b_{4}b_{2}-b_{3}^{2}b_{2}-b_{4}b_{2}^{2}$ \\[1pt]
$(1,1,3,3)$ & $-2b_{2}b_{4}-2b_{3}^{2}+b_{5}b_{1}+b_{1}b_{4}+2b_{2}b_{3}+b_{3}b_{2}-b_{1}^{2}b_{4}-b_{1}b_{4}b_{1}-b_{2}b_{1}b_{3}+b_{2}^{3}-2b_{2}b_{3}b_{1}-b_{3}b_{1}b_{2}-b_{3}b_{2}b_{1}-2b_{1}b_{2}^{2}+b_{1}^{2}b_{2}^{2}+b_{1}b_{2}b_{1}b_{2}+b_{1}b_{2}^{2}b_{1}$ & $-2b_{3}b_{5}-2b_{4}^{2}+b_{6}b_{2}+b_{2}b_{5}+2b_{3}b_{4}+b_{4}b_{3}-b_{2}b_{3}^{2}-b_{2}b_{4}b_{2}-b_{3}b_{2}b_{3}-b_{3}^{2}b_{2}-2b_{4}b_{2}^{2}-b_{2}^{2}b_{3}-b_{2}b_{3}b_{2}+3b_{2}^{4}$ \\[1pt]
$(1,1,4,2)$ & $-2b_{2}b_{4}-b_{3}^{2}+b_{5}b_{1}+b_{1}b_{4}+b_{2}b_{3}+b_{3}b_{2}-b_{1}^{2}b_{4}-b_{1}b_{4}b_{1}-b_{2}b_{1}b_{3}-b_{2}b_{3}b_{1}-b_{3}b_{1}b_{2}-b_{3}b_{2}b_{1}$ & $-2b_{3}b_{5}-b_{4}^{2}+b_{6}b_{2}+b_{2}b_{5}+b_{3}b_{4}+b_{4}b_{3}-b_{2}^{2}b_{4}-b_{2}b_{4}b_{2}-b_{3}b_{2}b_{3}-b_{3}^{2}b_{2}-2b_{4}b_{2}^{2}$ \\[1pt]
$(1,1,5,1)$ & $b_{5}b_{1}$ & $b_{6}b_{2}$ \\[1pt]
$(1,2,1,4)$ & $b_{2}b_{4}$ & $b_{3}b_{5}$ \\[1pt]
$(1,2,2,3)$ & $-b_{3}^{2}+b_{4}b_{2}+b_{2}b_{3}-b_{2}b_{1}b_{3}-2b_{2}^{3}$ & $-b_{4}^{2}+b_{5}b_{3}+b_{3}b_{4}-b_{2}b_{3}^{2}+b_{2}b_{4}b_{2}-2b_{3}b_{2}b_{3}-b_{3}^{2}b_{2}$ \\[1pt]
$(1,2,3,2)$ & $-2b_{2}b_{4}+b_{3}^{2}+b_{4}b_{2}+b_{1}b_{4}+2b_{2}b_{3}-b_{3}b_{2}-b_{1}b_{3}-b_{1}^{2}b_{4}-b_{1}b_{3}b_{2}-b_{2}b_{1}b_{3}-4b_{2}^{3}+b_{3}b_{1}b_{2}+b_{1}^{2}b_{3}+3b_{1}b_{2}^{2}-2b_{1}^{2}b_{2}^{2}-b_{1}b_{2}b_{1}b_{2}$ & $-2b_{3}b_{5}+b_{4}^{2}+b_{5}b_{3}+b_{2}b_{5}+2b_{3}b_{4}-b_{4}b_{3}-b_{2}b_{4}-2b_{2}^{2}b_{4}-2b_{3}b_{2}b_{3}-3b_{3}^{2}b_{2}+b_{4}b_{2}^{2}+2b_{2}^{2}b_{3}+2b_{2}b_{3}b_{2}-3b_{2}^{4}$ \\[1pt]
$(1,2,4,1)$ & $2b_{2}b_{4}+b_{3}^{2}+b_{4}b_{2}-b_{1}b_{4}-b_{2}b_{3}-b_{3}b_{2}+b_{1}^{2}b_{4}+b_{2}b_{1}b_{3}+b_{3}b_{1}b_{2}$ & $2b_{3}b_{5}+b_{4}^{2}+b_{5}b_{3}-b_{2}b_{5}-b_{3}b_{4}-b_{4}b_{3}+b_{2}^{2}b_{4}+b_{3}b_{2}b_{3}+b_{4}b_{2}^{2}$ \\[1pt]
$(1,3,1,3)$ & $b_{3}^{2}$ & $b_{4}^{2}$ \\[1pt]
$(1,3,2,2)$ & $2b_{3}^{2}-b_{2}b_{3}+b_{2}^{3}-b_{3}b_{1}b_{2}$ & $2b_{4}^{2}-b_{3}b_{4}-b_{2}^{2}b_{4}+b_{2}b_{3}^{2}+b_{3}^{2}b_{2}-b_{4}b_{2}^{2}$ \\[1pt]
$(1,3,3,1)$ & $4b_{2}b_{4}+2b_{3}^{2}-2b_{1}b_{4}-4b_{2}b_{3}+b_{1}b_{3}+2b_{1}^{2}b_{4}+2b_{2}b_{1}b_{3}+2b_{2}^{3}-b_{1}^{2}b_{3}-b_{1}b_{2}^{2}+b_{1}^{2}b_{2}^{2}$ & $4b_{3}b_{5}+2b_{4}^{2}-2b_{2}b_{5}-4b_{3}b_{4}+b_{2}b_{4}+2b_{2}^{2}b_{4}+2b_{3}b_{2}b_{3}+2b_{3}^{2}b_{2}-b_{2}^{2}b_{3}-b_{2}b_{3}b_{2}+b_{2}^{4}$ \\[1pt]
$(1,4,1,2)$ & $b_{4}b_{2}$ & $b_{5}b_{3}$ \\[1pt]
$(2,1,1,4)$ & $b_{1}b_{5}+2b_{2}b_{4}+b_{3}^{2}+b_{4}b_{2}-b_{5}b_{1}-b_{1}b_{4}-b_{2}b_{3}-b_{3}b_{2}+b_{1}b_{4}b_{1}+b_{2}b_{3}b_{1}+b_{3}b_{2}b_{1}$ & $b_{2}b_{6}+2b_{3}b_{5}+b_{4}^{2}+b_{5}b_{3}-b_{6}b_{2}-b_{2}b_{5}-b_{3}b_{4}-b_{4}b_{3}+b_{2}b_{4}b_{2}+b_{3}^{2}b_{2}+b_{4}b_{2}^{2}$ \\[1pt]
$(2,1,2,3)$ & $b_{2}b_{4}-b_{1}b_{2}b_{3}-b_{1}b_{3}b_{2}-2b_{2}^{3}+b_{3}b_{1}b_{2}+b_{1}b_{2}^{2}-b_{1}b_{2}b_{1}b_{2}$ & $b_{3}b_{5}-b_{2}b_{3}^{2}-b_{2}b_{4}b_{2}-2b_{3}^{2}b_{2}+b_{4}b_{2}^{2}+b_{2}b_{3}b_{2}-b_{2}^{4}$ \\[1pt]
$(2,1,3,2)$ & $-b_{3}^{2}+b_{4}b_{2}+b_{2}b_{3}-b_{1}b_{3}b_{2}-2b_{2}^{3}$ & $-b_{4}^{2}+b_{5}b_{3}+b_{3}b_{4}+b_{2}^{2}b_{4}-2b_{2}b_{3}^{2}-b_{3}b_{2}b_{3}-b_{3}^{2}b_{2}$ \\[1pt]
$(2,2,1,3)$ & $b_{2}b_{4}+2b_{3}^{2}-b_{4}b_{2}-b_{2}b_{3}-b_{1}b_{2}b_{3}+b_{2}^{3}$ & $b_{3}b_{5}+2b_{4}^{2}-b_{5}b_{3}-b_{3}b_{4}-b_{2}^{2}b_{4}+b_{2}b_{3}^{2}-b_{2}b_{4}b_{2}+b_{3}b_{2}b_{3}$ \\[1pt]
$(2,2,2,2)$ & $2b_{2}b_{4}-b_{1}b_{4}-2b_{2}b_{3}+b_{3}b_{2}+b_{1}b_{3}+b_{1}^{2}b_{4}+b_{1}b_{3}b_{2}+4b_{2}^{3}-3b_{3}b_{1}b_{2}-b_{1}^{2}b_{3}-4b_{1}b_{2}^{2}+2b_{1}^{2}b_{2}^{2}+3b_{1}b_{2}b_{1}b_{2}$ & $2b_{3}b_{5}-b_{2}b_{5}-2b_{3}b_{4}+b_{4}b_{3}+b_{2}b_{4}+b_{2}^{2}b_{4}+b_{2}b_{4}b_{2}+4b_{3}^{2}b_{2}-3b_{4}b_{2}^{2}-2b_{2}^{2}b_{3}-3b_{2}b_{3}b_{2}+5b_{2}^{4}$ \\[1pt]
$(3,1,1,3)$ & $b_{1}b_{5}+2b_{2}b_{4}+3b_{3}^{2}-b_{5}b_{1}-b_{1}b_{4}-2b_{2}b_{3}-b_{3}b_{2}-b_{1}b_{2}b_{3}+b_{1}b_{4}b_{1}+2b_{2}b_{3}b_{1}+b_{3}b_{2}b_{1}+b_{1}b_{2}^{2}-b_{1}b_{2}^{2}b_{1}$ & $b_{2}b_{6}+2b_{3}b_{5}+3b_{4}^{2}-b_{6}b_{2}-b_{2}b_{5}-2b_{3}b_{4}-b_{4}b_{3}-b_{2}^{2}b_{4}+b_{2}b_{4}b_{2}+2b_{3}^{2}b_{2}+b_{4}b_{2}^{2}+b_{2}^{2}b_{3}-b_{2}^{4}$ \\[1pt]
\end{longtable}}

{\footnotesize
\begin{longtable}{@{}l@{\quad}>{\raggedright\arraybackslash}p{0.87\textwidth}@{}}
\caption{The polynomials $\mathcal L^{0}(x^{a}y^{b})$ for $a\geq b\geq1$ and $a+b\leq6\mspace{1mu}$.}\label{tab:L}\\
\toprule
$p$ & $\mathcal L^{0}(p)$\\
\midrule
\endfirsthead
\toprule
$p$ & $\mathcal L^{0}(p)$\\
\midrule
\endhead
\bottomrule
\endlastfoot
$x^{1}y^{1}$ & $\bigl(\tfrac{1}{2}x_1+\tfrac{1}{2}x_1^{2}\bigr)\oplus \bigl(-x_2\bigr)$ \\[1pt]
$x^{2}y^{1}$ & $\bigl(\tfrac{1}{6}x_1+\tfrac{1}{2}x_1^{2}+\tfrac{1}{3}x_1^{3}\bigr)\oplus \bigl(-x_2^{2}-2x_1x_2\bigr)$ \\[1pt]
$x^{3}y^{1}$ & $\bigl(\tfrac{1}{4}x_1^{2}+\tfrac{1}{2}x_1^{3}+\tfrac{1}{4}x_1^{4}\bigr)\oplus \bigl(-x_2^{3}-3x_1x_2^{2}-3x_1^{2}x_2\bigr)$ \\[1pt]
$x^{2}y^{2}$ & $\bigl(\tfrac{1}{6}x_1+\tfrac{1}{3}x_1^{2}+\tfrac{1}{3}x_1^{3}+\tfrac{1}{6}x_1^{4}\bigr)\oplus \bigl(-\tfrac{1}{3}x_2-2x_1x_2-\tfrac{2}{3}x_2^{3}-4x_1x_2^{2}-2x_1^{2}x_2\bigr)\oplus 4x_2x_3$ \\[1pt]
$x^{4}y^{1}$ & $\bigl(-\tfrac{1}{30}x_1+\tfrac{1}{3}x_1^{3}+\tfrac{1}{2}x_1^{4}+\tfrac{1}{5}x_1^{5}\bigr)\oplus \bigl(-x_2^{4}-4x_1x_2^{3}-6x_1^{2}x_2^{2}-4x_1^{3}x_2\bigr)$ \\[1pt]
$x^{3}y^{2}$ & $\bigl(\tfrac{1}{15}x_1+\tfrac{1}{4}x_1^{2}+\tfrac{1}{3}x_1^{3}+\tfrac{1}{4}x_1^{4}+\tfrac{1}{10}x_1^{5}\bigr)\oplus \bigl(-\tfrac{1}{2}x_2^{2}-2x_1x_2-3x_1x_2^{2}-3x_1^{2}x_2-\tfrac{1}{2}x_2^{4}-4x_1x_2^{3}-6x_1^{2}x_2^{2}-2x_1^{3}x_2\bigr)\oplus \bigl(6x_2x_3^{2}+6x_2^{2}x_3+12x_1x_2x_3\bigr)$ \\[1pt]
$x^{5}y^{1}$ & $\bigl(-\tfrac{1}{12}x_1^{2}+\tfrac{5}{12}x_1^{4}+\tfrac{1}{2}x_1^{5}+\tfrac{1}{6}x_1^{6}\bigr)\oplus \bigl(-x_2^{5}-5x_1x_2^{4}-10x_1^{2}x_2^{3}-10x_1^{3}x_2^{2}-5x_1^{4}x_2\bigr)$ \\[1pt]
$x^{4}y^{2}$ & $\bigl(-\tfrac{1}{30}x_1+\tfrac{1}{10}x_1^{2}+\tfrac{1}{3}x_1^{3}+\tfrac{1}{3}x_1^{4}+\tfrac{1}{5}x_1^{5}+\tfrac{1}{15}x_1^{6}\bigr)\oplus \bigl(\tfrac{1}{15}x_2-\tfrac{2}{3}x_2^{3}-4x_1x_2^{2}-4x_1^{2}x_2-4x_1x_2^{3}-6x_1^{2}x_2^{2}-4x_1^{3}x_2-\tfrac{2}{5}x_2^{5}-4x_1x_2^{4}-8x_1^{2}x_2^{3}-8x_1^{3}x_2^{2}-2x_1^{4}x_2\bigr)\oplus \bigl(8x_2x_3^{3}+12x_2^{2}x_3^{2}+8x_2^{3}x_3+24x_1x_2x_3^{2}+24x_1x_2^{2}x_3+24x_1^{2}x_2x_3\bigr)$ \\[1pt]
$x^{3}y^{3}$ & $\bigl(\tfrac{1}{10}x_1+\tfrac{1}{5}x_1^{2}+\tfrac{1}{4}x_1^{3}+\tfrac{1}{4}x_1^{4}+\tfrac{3}{20}x_1^{5}+\tfrac{1}{20}x_1^{6}\bigr)\oplus \bigl(-\tfrac{1}{5}x_2-3x_1x_2-\tfrac{1}{2}x_2^{3}-6x_1x_2^{2}-\tfrac{9}{2}x_1^{2}x_2-3x_1x_2^{3}-9x_1^{2}x_2^{2}-3x_1^{3}x_2-\tfrac{3}{10}x_2^{5}-3x_1x_2^{4}-9x_1^{2}x_2^{3}-6x_1^{3}x_2^{2}-\tfrac{3}{2}x_1^{4}x_2\bigr)\oplus \bigl(6x_2x_3+18x_1x_2x_3+6x_2x_3^{3}+18x_2^{2}x_3^{2}+6x_2^{3}x_3+36x_1x_2x_3^{2}+36x_1x_2^{2}x_3+18x_1^{2}x_2x_3\bigr)\oplus \bigl(-36x_2x_3x_4\bigr)$ \\[1pt]
\end{longtable}}

\end{document}